\documentclass[3p,times]{elsarticle}

\usepackage{ulem}

\usepackage{caption}
\usepackage{subcaption}
\usepackage{xcolor}

\usepackage{amssymb}
\usepackage[figuresright]{rotating}

\begin{document}

\begin{frontmatter}

\title{Nonlinear dynamics of asymmetric bistable energy harvesters}

\author[unesp]{Jo\~ao Pedro~Norenberg}
\author[uerj]{Roberto Luo}
\author[uerj]{Vinicius Gon\c{c}alves Lopes}
\author[uerj]{\\Jo\~{a}o Victor~L.~L.~Peterson}
\author[uerj]{Americo~Cunha~Jr\corref{cor}}
\cortext[cor]{Corresponding author.}

\address[unesp]{S\~ao Paulo State University, Ilha Solteira, SP, Brazil}
\address[uerj]{Rio de Janeiro State University, Rio de Janeiro, RJ, Brazil}

\begin{abstract}
The paper investigates asymmetries effects over a nonlinear vibration energy harvester dynamics. The asymmetric system performance is compared with symmetric ones. Different asymmetry levels on restoring force and gravity action are investigated from a system-sloping angle variation. Bifurcation diagrams and basins of attraction are used to examine the local and global characteristics underlying dynamical systems under different excitation energy. The results show the adverse effects of asymmetries on system dynamics. They also reveal ways to overcome them by canceling asymmetric influence from optimal sloping angle values and improving asymmetric system performance over symmetrical ones. This comprehensive numerical study provides novel valuable insights into asymmetrical energy harvester dynamics, a wide and still less explored topic.
\end{abstract}


\end{frontmatter}

\section{Introduction}

From recent years, the increasing Internet of Things (IoT) technology diffusion has been claiming the need for alternatives on the power supplying embedded and portable electronic devices \cite{NORD201997}. Energy harvesting systems have been emerging as promising solutions in this area, with potential applications ranging from large-scale energy supply to power small devices, such as sensors, medical implants, and micro/nano-electromechanical systems (MEMS/NEMS)~\cite{Koka_2014,SEOL_2013}.

Due to its wide potential applications, energy harvesters have been extensively investigated. Recent papers of Abdelkareem et al. \cite{ABDELKAREEM_2018} presented a regenerative energy harvesting method from the automotive suspension; Phillips \cite{Phillips_2021} investigated nanoscale sensors, Karami and Inman \cite{Karami_2012} proposed harvesting energy from heart vibrations to power pacemakers and for biomedical applications \cite{panda2022}. Catacuzzeno et al. \cite{Catacuzzeno_2019} proposed energy harvesting from cell membrane on powering wireless communication devices. The design and implementation of novel energy harvesting devices have also been proposed and discussed in recent works as \cite{Kazemi_2021,Zhong_2022,Zhou_2021}. The use of metamaterials has also been explored as a potential environmental energy harvester alternative \cite{Lee_2022,XIAO2023109808,tan2019, sun2022, Lu2020p105532,Miranda_2020,Nie2019p287,Rezaei2021p106618,Tao_2021,Vijayan2015p101, cao2021}.

This paper investigates the nonlinear bistable piezo-magneto-elastic energy harvester device by Cottone et al. \cite{cottone2009p080601}, which harvests electrical power from a cantilever ferromagnetic beam mechanical vibration. The magnetic forces provided by a pair of permanent magnets placed under the beam's free edge confer its nonlinear dynamic characteristics to the device. Unlike linear energy harvesting devices, which non-negligible electrical output power is limited on resonant excitation frequency, nonlinear systems, as those proposed by Erturk et al. \cite{erturk2009p254102}, provide a wide excitation frequency range for useful electrical power harvesting. For these devices, energy harvesting improvements from system dynamics have been widely explored, as in the papers of Daqaq et al. \cite{19_dAQAQ}, where the nonlinearities role on energy transduction is reviewed, Mann and Owens \cite{Mann_2010}, where the dynamics mathematical solutions co-existence is investigated for broader frequency response, Tan et al. \cite{Tan_2020}, where optimal performance is analyzed, and Khovanov et al. \cite{Khovanov_2021} and Litak et al. \cite{Litak_2010} which deals with system stochastic excitation. Other studies have investigated the nonlinear vibration energy harvesting from various perspectives, including polynomial nonlinearities \cite{Lua_2020,Erturk_2011,Lopes,Peterson_2017}, nonlinear piezoelectric coupling \cite{13_triplett,Leadenham_2015}, variable potentials \cite{Chen2023p107997}, rectifier circuits \cite{Basquerotto_2015, xu2022}, shape memory alloy \cite{Adeodato2021p106206}, rotating elements \cite{Wang2023p108033}, combining with vibration suppression \cite{REZAEI2021106618,CHEN2023107997} and using time-frequency analysis \cite{Varanis_2020} and global sensitivity analysis \cite{norenbergNoDy_2022,Norenberg2021_cobem}. Some authors have also been proposing system improvements, as
\cite{CUNHA_2021,Firoozy2017p227,Huang2022p107609,deLaRoca_2019} or exploring tri-stable, multi-stable or triboelectric configurations \cite{Egbe2022p107722,Kim_2014,Litak_2021,Zhou_2014}.

Despite its relevance and such huge potential applications, testified by much promising research, much attention has been paid to symmetric systems. However, symmetric systems manufacturing is usually hard due to a wide of disturbances, such as those from its manufacturing process and environmental or operational conditions, like humidity, temperature, and nearby magnetic fields. Recent studies \cite{Cao_2015,Halvorsen,He_asymmetric,Wang_asymmetric_2018,Wang_asymmetric_2019} have recognized the importance of asymmetric potentials and characterized them as nonlinear restoring forces using a quadratic term. In many real bistable harvesters, this asymmetric behavior can result from using non-identical magnets. Further, gravity states as another effect that can introduce asymmetric forces resulting from geometric imperfections and operating conditions. For example,  in \cite{Wang2018}, authors explore the harvesting process under asymmetric excitation, which resembles human lower-limb motion; Wang et al. \cite{Wang_bias_2018} uses from sloping the plane where the system is attached to improve harvesting efficiency. The understanding of the impacts of asymmetric layouts on bistable energy harvesters may lead to improvements to devices' design and performance improvements, even reducing manufacturing costs and limitations, which may pave the way for new potential unexplored applications.

However, the current technical literature lacks a comprehensive analysis of the dynamic behavior of asymmetric bistable energy harvesters. Most of the available studies have concentrated on investigating specific situations by fixing all parameters except for one parameter of interest. For example, Balakrishnan et al. \cite{Balakrishnan_2019} briefly investigated the bifurcation diagrams from varying excitation frequency and the nonlinear term. Meanwhile, works in \cite{Wang_bias_2018,Dong_2021} calculated bifurcations for both frequency and amplitude of excitation. Wang et al. \cite{Wang_asymmetric_2018} also studied basins of attraction, but only for different values of the nonlinear term and excitation amplitude. This work addresses this relevant topic by comprehensively analyzing asymmetry on bistable energy harvesting devices.
To introduce asymmetric features, a quadratic term is added to restoring force to mimic non-identical magnets. A sloping angle is introduced as the asymmetric gravity force.

This paper investigates dynamical system for different excitation amplitude and frequencies for different initial conditions and asymmetry levels. The symmetric energy harvester dynamics is compared with  devices on different asymmetry conditions. The analysis focuses on a high asymmetry level model, where a single quadratic term refers to its asymmetry on large sloping angle values and a special asymmetry condition where such an angle suppresses the quadratic term influence. Bifurcation diagrams and basin of attraction analyses are held to characterize the dynamic behavior and identify any instances of chaos or multiple solutions. For electrical output voltage dynamics, basins of attraction are computed from the 0-1 test for chaos \cite{gottwald2016}, which offers a more efficient computational approach when compared to the usual Lyapunov exponent method.

As its main contributions, the paper addresses a comprehensive investigation of energy harvester dynamics for different operational conditions, even considering asymmetric configurations. It provides a deep dynamics analysis while showing a wide oversight of its best output excitation conditions. A nonlinear dynamics analysis and solutions co-existence is also contemplated for asymmetry conditions. Asymmetry energy harvester dynamics is compared with symmetric devices' behavior, evidencing asymmetry impacts on bistable energy harvester performance. The paper's scope matches no similar precedent analysis in the current literature to the best of the author's knowledge.

Section 2 presents symmetric and asymmetric dynamical systems investigated; Section 3 covers the bifurcation analysis for different excitation amplitude and frequencies and system sloping angle ; Section 4 provides an analysis on basins of attraction ; Finally, the paper main conclusions are presented in Section 5.

\section{Nonlinear dynamical systems}
\label{nonlinear_syst_sect}

This section describes the physical system and presents the governing equations of motion for the nonlinear vibration bistable energy harvester both on perfect symmetric and asymmetric conditions. It also outlines a methodology for defining the optimal sloping angle value on compensating the asymmetry effects.

\subsection{Physical system of interest}

\begin{figure}
\centering
    \begin{subfigure}[b]{0.48\textwidth}
         \centering
         \includegraphics[width=0.75\textwidth]{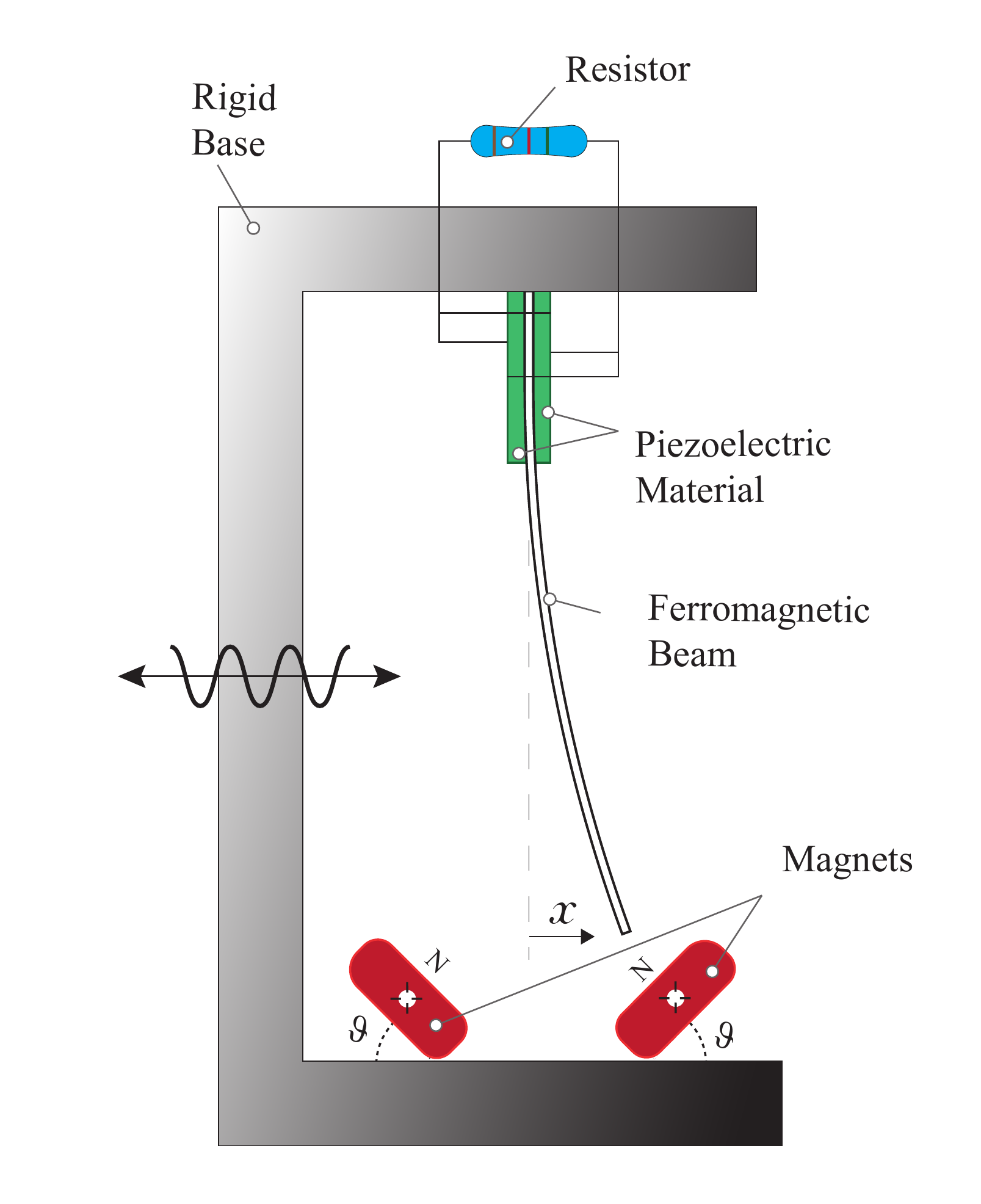}
         \caption{symmetric system}
     \end{subfigure}
     \begin{subfigure}[b]{0.48\textwidth}
         \centering
         \includegraphics[width=0.75\textwidth]{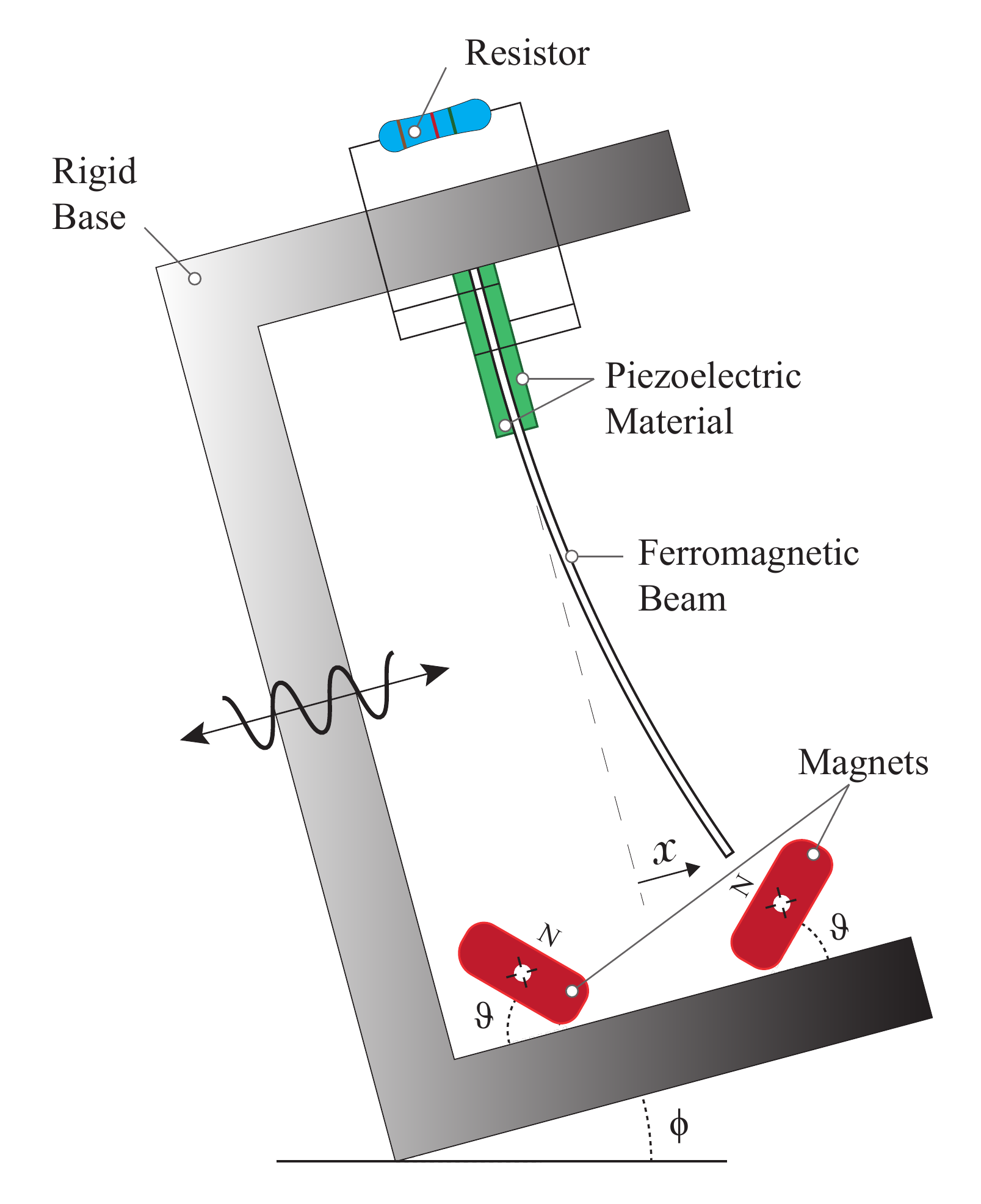}
         \caption{asymmetric system}
     \end{subfigure}    
    \caption{The illustration depicts the two configurations of the piezo-magneto-elastic energy harvesters: (a) the symmetric configuration and (b) the asymmetric configuration.}
    \label{harvesting_device_fig}
\end{figure}

The electromechanincal system of interest is shown in Fig.\ref{harvesting_device_fig}. The system consists in a cantilever ferromagnetic elastic beam attached to a rigid base. A pair of identical magnets are placed in the bottom of the structure, under the ferromagnetic beam's free edge. A pair of piezoelectric plates are placed close to the beam's fixed edge, powering a purely resistive load from harvesting mechanical vibration into electrical power. The system mechanical excitation is provided by an external oscillation amplitude and frequency-controlled source. Fig.\ref{harvesting_device_fig}a presents the symmetric bistable energy harvester. In contrast, Fig.~\ref{harvesting_device_fig}b shows the asymmetric version, which includes a sloping plane angle $\phi$ and non-identical magnets for emulating asymmetries from manufacturing and operation processes.

Erturk et al. \cite{erturk2009p254102} presented the symmetric bistable energy harvester described above, which evolves according to the following initial value problem
\begin{equation}
	\ddot{x} + 2 \, \xi \, \dot{x} - \frac{1}{2} x \left(1-x^2 \right) - \chi \, v 
	= f \, \cos{\left( \Omega \, t \right)} \, ,
	\label{ode_mechanical}
\end{equation}
\begin{equation}
	\dot{v} + \lambda \, v + \kappa \, \dot{x} = 0 \, ,
	\label{ode_electrical}
\end{equation}
\begin{equation}
	x(0) = x_0, ~ \dot{x}(0) = \dot{x}_0, ~ v(0) = v_0 \, ,
	\label{initial_conditions}
\end{equation}
where $t$ denotes time; $x$ is the modal amplitude of oscillation; $v$ is the voltage in the resistor; $\xi$ is the damping ratio; $f$ is the rigid base oscillation amplitude; $\Omega$ is the external excitation frequency; $\lambda$ is a reciprocal time constant; the piezoelectric coupling terms are represented by $\chi$, in the mechanical equation, and by $\kappa$ in the electrical one; $x_{0}$ represents the modal initial position;  $\dot{x}_{0}$ is the modal initial velocity; $v_{0}$ denotes the initial voltage over the resistor. The upper dot is an abbreviation for time derivative. These parameters are all dimensionless.

In their study \cite{Cao_2015}, Cao et al. first introduced an asymmetry on a bistable energy harvester by adding a quadratic term to the nonlinear restoring force. This term accounted for the asymmetry caused by the magnets' eccentricity, previously verified by Halvorsen \cite{Halvorsen} and He and Daqaq \cite{He_asymmetric}. In the sequence, a self-gravity force acting on the harvester was also modeled to characterize the influence of a sloping plane. It is worth noting that this force can be considered an asymmetric external excitation with its intensity depending on the slope angle $\phi$. The resulting lumped-parameter equations of motion that govern the behavior of the asymmetric bistable energy harvester are given by
\begin{equation}
        \ddot{\mathnormal{x}}+2\;\xi\;\dot{\mathnormal{x}}-\frac{1}{2}\;\mathnormal{x}\;(1+2\delta\mathnormal{x}-\mathnormal{x}^2)-\chi\; \mathnormal{v}  = 
        \mathnormal{f}\;\cos{\left(\Omega\; t\right)} + \mathnormal{p}\;\sin{\left(\phi\right)} \, ,
    \label{ode_mechanical_asymmetric}
\end{equation}
\begin{equation}
    \dot{\mathnormal{v}} + \lambda\;\mathnormal{v}+\kappa\; \dot x = 0 \, ,
    \label{ode_eletrical_asymmetric}
\end{equation}
\begin{equation}
    \mathnormal{x}(0) = \mathnormal{x}_0, \;\dot{\mathnormal{x}}(0) = \dot{\mathnormal{x}}_0, \; \mathnormal{v}(0) = \mathnormal{v}_0 \, ,
\end{equation}
where $\delta$ is a coefficient of the quadratic nonlinearity, $\mathnormal{p}$ is the equivalent dimensionless gravity of ferromagnetic beam, and $\phi$ is the sloping angle. 


\subsection{Optimal sloping angle}

The authors in \cite{Wang_asymmetric_2018,Wang_asymmetric_2019} have highlighted that asymmetry can negatively impact the process of energy harvesting. To counteract this, Wang et al. \cite{Wang_bias_2018} verified that an optimal sloping angle can be used to compensate for the asymmetry in the potential functions. In this way, the asymmetric excitation can mitigate the effect of the asymmetric quadratic term. Considering the stable solutions of the model's nonlinear restoring force
\begin{equation}
    F_r(x) = -\frac{1}{2}x(1+2\delta~x - x^2) - p\sin{\left(\phi\right)} \,
    \label{eq:restoring}
\end{equation}
can be found by assuming that its derivative is equal to zero, such that its stable solutions are given by
\begin{equation}
    x_{1,2} = \frac{2\delta\pm\sqrt{4\delta^2+3}}{3} \, .
    \label{eq:x12}
\end{equation}

To ensure that the restoring force has the same intensity on equilibrium points, the sum of the nonlinear restoring force at the stable solutions of Eq. (\ref{eq:x12}) must be zero, i.e., $F_r(x_1) + F_r(x_2) = 0$. This condition enables the calculation of the optimal sloping angle value
\begin{equation}
    \phi_{opt} = \sin^{-1}{\left(\frac{8\delta^3+9\delta}{27 p}\right)} \, .
    \label{eq:optimal_angle}
\end{equation}

\begin{figure}[!]
    \centering
    \includegraphics[width=0.45\textwidth]{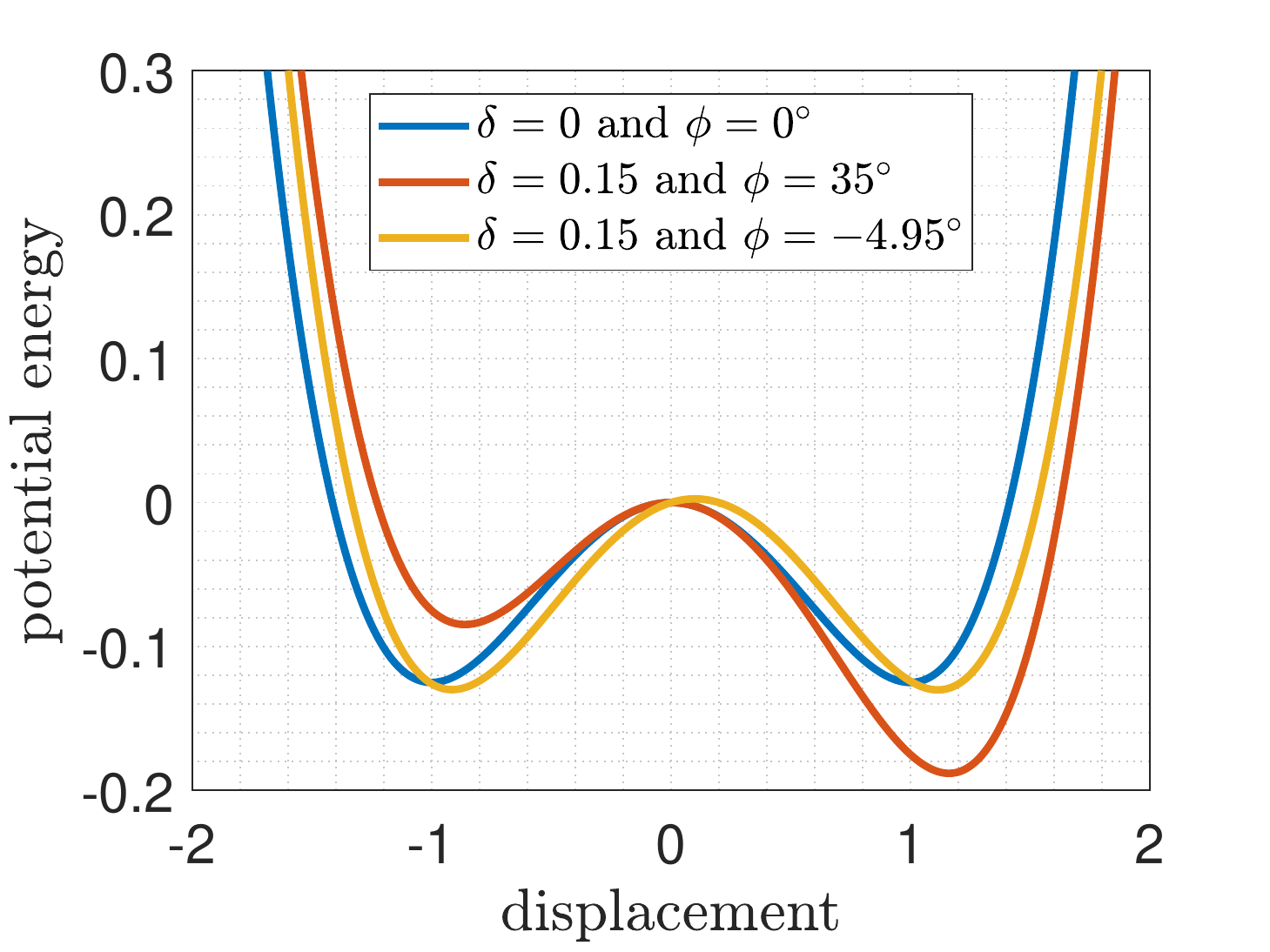}
    \caption{Potential energy of the symmetric ($\delta=0$ and $\phi=0^\circ$) , asymmetric ($\delta=0.15$ and $\phi=35^\circ$) and optimal asymmetric ($\delta=0.15$ and $\phi_{opt}=-4.95^\circ$) models. Optimal asymmetric and symmetric potentials energy are symmetrical deep well, but they present different equilibrium points. The asymmetric potential energy has the right well deeper than the left well.}
    \label{fig:potential}
\end{figure}

Fig.\ref{fig:potential} shows three types of potential energy: symmetric, asymmetric, and optimal asymmetric (when the sloping angle is optimal for compensating the asymmetry effect). The potential energy is obtained by integrating the Eq.~(\ref{eq:restoring}). The asymmetric potential has a deep and a shallow well. The symmetric and compensated potentials have similar shapes, with equally deep wells but with different equilibrium points.


\section{Bifurcation analysis}
\label{sec_bifurcation_analysis}

This section focuses on numerical simulations of bifurcation diagrams for the symmetric and asymmetric bistable energy harvesters. The dynamic behavior of increasing and decreasing amplitude and frequency values, referred to as forward and backward, respectively, and the slope angle are investigated. The methodology of bifurcation calculus is also approached.

The frequency and amplitude sampling intervals are set to $0.1 \leq \Omega \leq 1.4$ and $0.01 \leq f \leq 0.3$, respectively. The sloping angle intervals are defined from $-35^\circ$ to $35^\circ$. The remaining model parameters are fixed at $\xi = 0.01$, $\chi = 0.05$, $\lambda = 0.05$, $\kappa = 0.5$, $p=0.59$, and $\delta = 0.15$. The initial conditions are $(x_{0},\dot{x}_{0},\upsilon_{0}) = (1,0,0)$.
For the amplitude and frequency bifurcation diagram, an asymmetric system with $\phi=35^\circ$ is analyzed to investigate the effect of high-intensity of asymmetry excitation. 
The special condition when the asymmetry of excitation mitigates the asymmetry of the quadratic term is also explored. For the given $\delta$, the optimal angle value $\phi_{opt}$ is equal to $-4.95^\circ$, obtained from Eq.~(\ref{eq:optimal_angle}). To integrate the dynamics, we used the fourth-order Runge-Kutta method with a relative tolerance of $10^{-6}$ and an absolute tolerance of $10^{-9}$.

\begin{figure}
	\centering
    \includegraphics[width=0.85\textwidth]{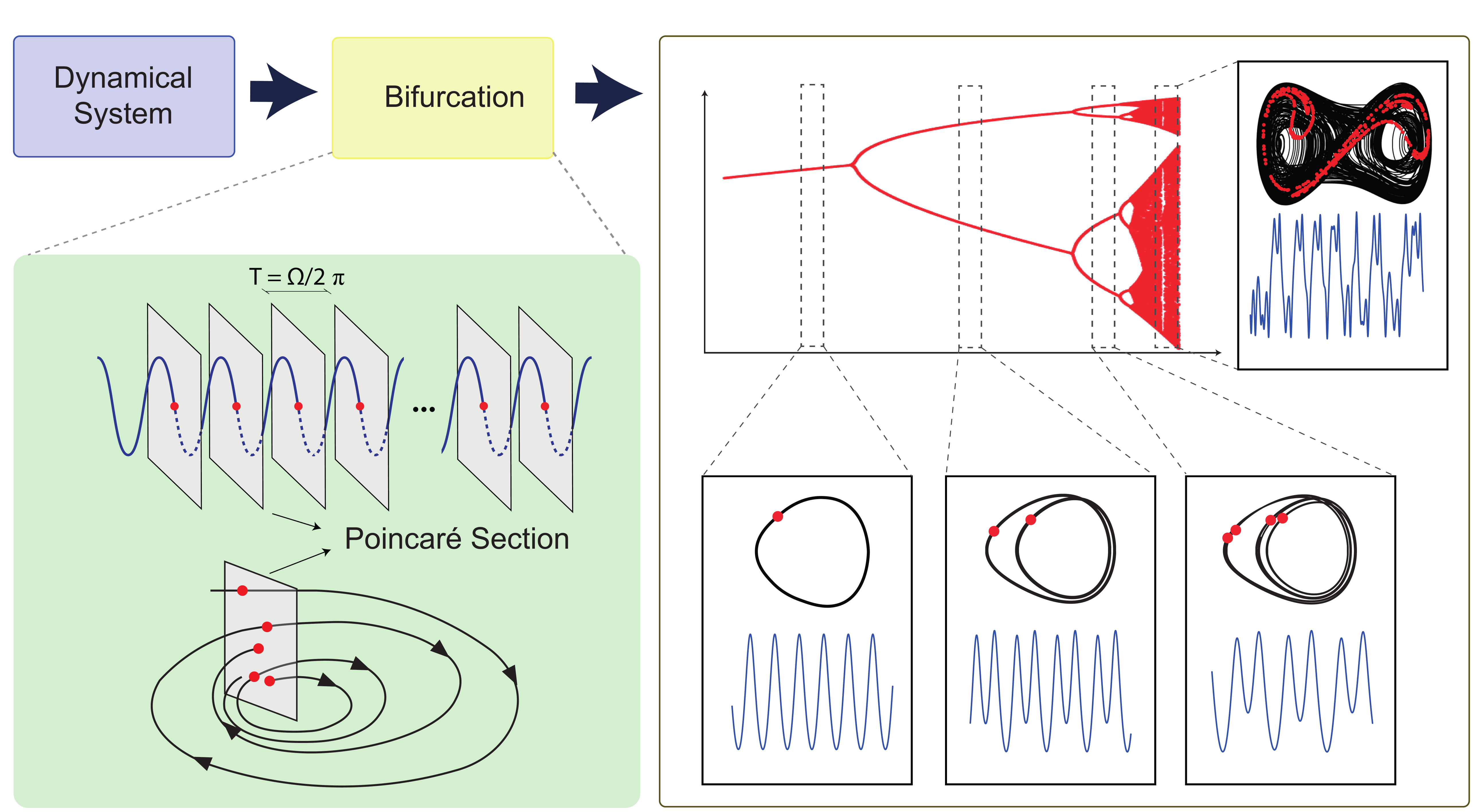}
	\caption{Visual depiction of the procedure for creating bifurcation diagrams. The dynamical system is integrated using a fixed time step $T=\frac{\Omega}{2\pi}$ to obtain the Poincaré section of the time series. This process is repeated while varying the control parameter. Finally, the bifurcation diagram is plotted from the Poincaré section of the time series according to their control parameter. An intersection of the bifurcation diagram indicates a one-period response while doubling periods leads to a chaotic response.}
	\label{fig:bifur_illustration}
\end{figure}

The bifurcation diagrams are created by varying a parameter of interest, such as frequency, amplitude, or sloping angle, in either a forward or backward direction. The diagrams are generated using the system's steady-state condition, defined as the last $10\%$ of the response time series. For each excitation amplitude or frequency, the initial conditions are updated with the last system state and sampled from a set of $1200$ regularly spaced values within the analysis intervals. The simulation is run for $1000$ forcing cycles. Fig.\ref{fig:bifur_illustration} illustrates the proceeding to obtain the bifurcation diagrams using the Poincaré section.

The frequency observation window is set from $\Omega = 0.1$ to $\Omega = 0.9$, with a step of $0.1$. In contrast, the excitation amplitude, with frequency as the parameter of interest, is set from $f = 0.019$ to $f = 0.275$, with a step of $0.032$. The defined range of values corresponds to feasible environmental conditions that, according to our previous work \cite{norenbergNoDy_2022}, at this range of values can significantly affect the motion of the system. The nonlinear dynamics effects are studied by comparing the symmetric and asymmetric models. The forward and backward bifurcation diagrams are presented, with voltage-time series sampled from regions of interest selected by visual inspection of the diagrams. 
The voltage-time series in the Supplementary Material provides additional information for readers to understand the dynamic behavior for several points of bifurcation diagrams.

\subsection{Amplitude analysis}
\label{subsec_amplitude_analysis}

\begin{figure}[!]
\centering
\includegraphics[width = 1\textwidth]{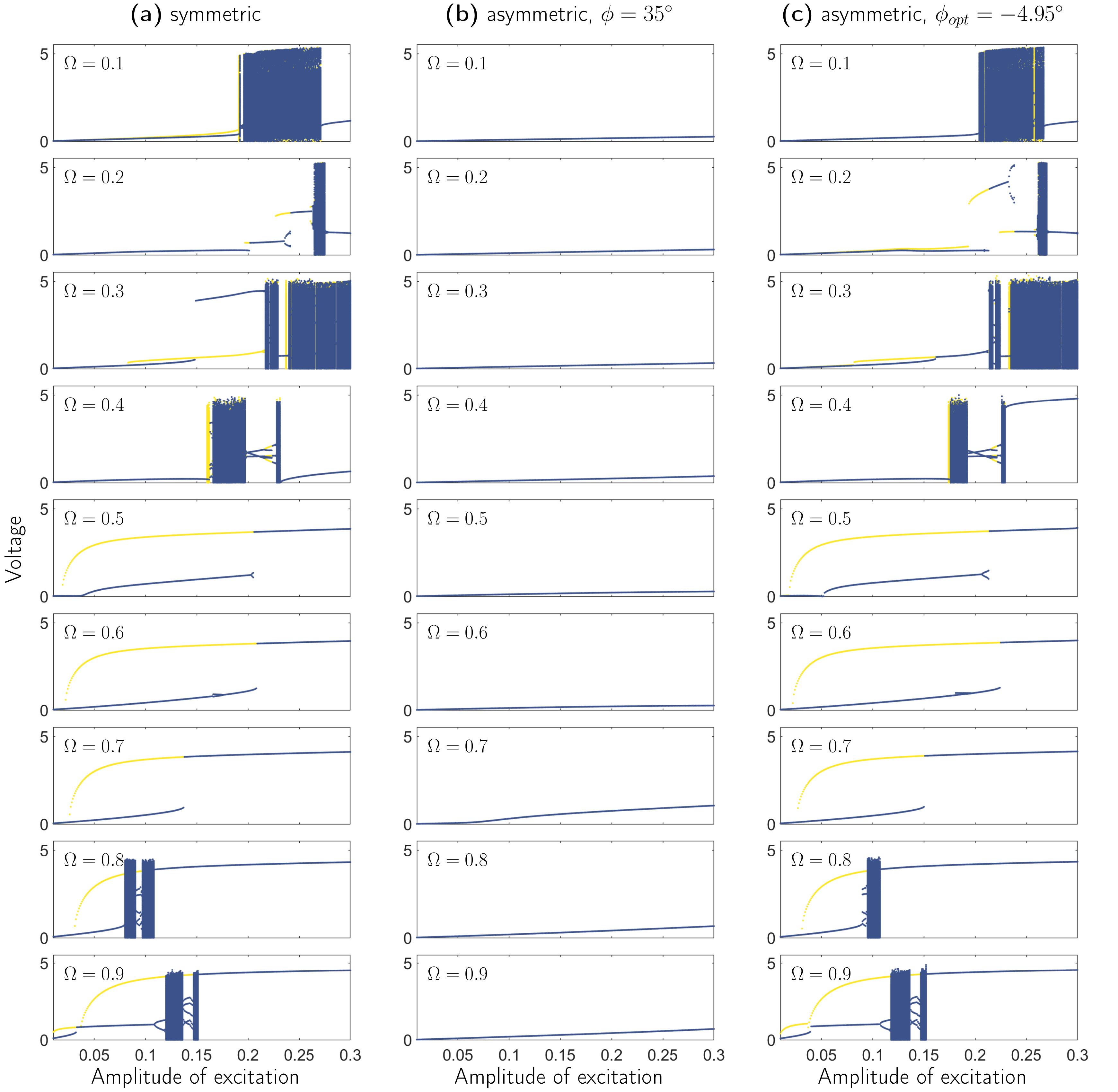}
\caption{Bifurcation diagram of output voltage in bistable energy harvesters as a function of excitation amplitude for different frequency values. The forward sweep is displayed in blue, and the backward sweep in yellow. (a) symmetric harvester, (b) asymmetric harvester with $\phi = 35^\circ$, (c) asymmetric harvester with $\phi_{opt} = -4.95^\circ$.}
\label{fig:plot3_amplitude}
\end{figure}

The voltage bifurcation diagrams of the symmetric and asymmetric (high asymmetry, $\phi=35^\circ$, and optimal asymmetry, $\phi_{opt}=-4.95^\circ$, conditions) energy harvesters as a function of the excitation amplitude under different frequencies are displayed in Fig.~\ref{fig:plot3_amplitude}. The results show the dynamic behavior for all three cases, revealing a wealth of chaotic and multi-periodicity responses.


Fig.~\ref{fig:plot3_amplitude}a displays the symmetrical nonlinear voltage forward and backward bifurcation diagrams, displaying regular multiple-period regions and discontinuities. High frequencies result in better energy recovery, with regular dynamics observed for $\Omega\leq0.7$ at the beginning of the amplitude interval. The backward diagrams show an increase in voltage output for frequencies from $0.5$ onwards and a smooth increase starting at $f=0.015$. The highest amplitudes at $\Omega=0.5$ result in multiple period regions on the backward diagrams, while higher frequencies show regular and unchanged diagrams on both forward and backward sweeping. The voltage reaches a steady-state value for $f\geq0.2$ and $\Omega=0.6$ onwards. However, lower frequencies lead to chaotic behavior. The forward diagrams show massive regions of chaotic behavior for $\Omega=0.8$ and $\Omega=0.9$ for $f<0.2$.

Fig.~\ref{fig:plot3_amplitude}b displays the asymmetrical nonlinear voltage forward and backward bifurcation diagrams for $\phi=35^\circ$. The bifurcation diagrams did not reveal a wealthy dynamic behavior. The system exhibits a unique periodic solution for all amplitude and frequency ranges, but its severe asymmetry results in low voltage output amplitude. It occurs when the external excitation cannot overcome the potential barrier, causing the system to remain in intra-well motion. To transition from one stable point to another and perform inter-well motion, a significant amount of external energy is necessary.

Fig.~\ref{fig:plot3_amplitude}c shows the bifurcation diagram for the asymmetric system with the optimal sloping angle, $\phi_{opt}=-4.95^\circ$. The bifurcation diagram reveals multiple periodic regions for the highest amplitudes when $\Omega\leq0.4$. For $\Omega\geq0.8$, chaotic regions also appear for $f\leq0.2$. The forward and backward bifurcation diagrams demonstrate the coexistence of solutions when $\Omega\geq0.5$ and $f\leq0.2$. However, for this range of frequency and $f\geq0.2$, they indicate the same response. The dynamic behavior under different excitations is similar to the symmetric dynamic behavior. Chaotic behavior occurs under almost the same $f$ and $\Omega$ conditions. Both forward and backward bifurcations reveal the same multiple solutions. Moreover, the voltage output for the optimal slope angle value resulted in a higher amplitude than $\phi=35^\circ$.

\subsection{Frequency analysis}
\label{subsec_frequency_analysis}

\begin{figure}
\centering
\includegraphics[width = 1\textwidth]{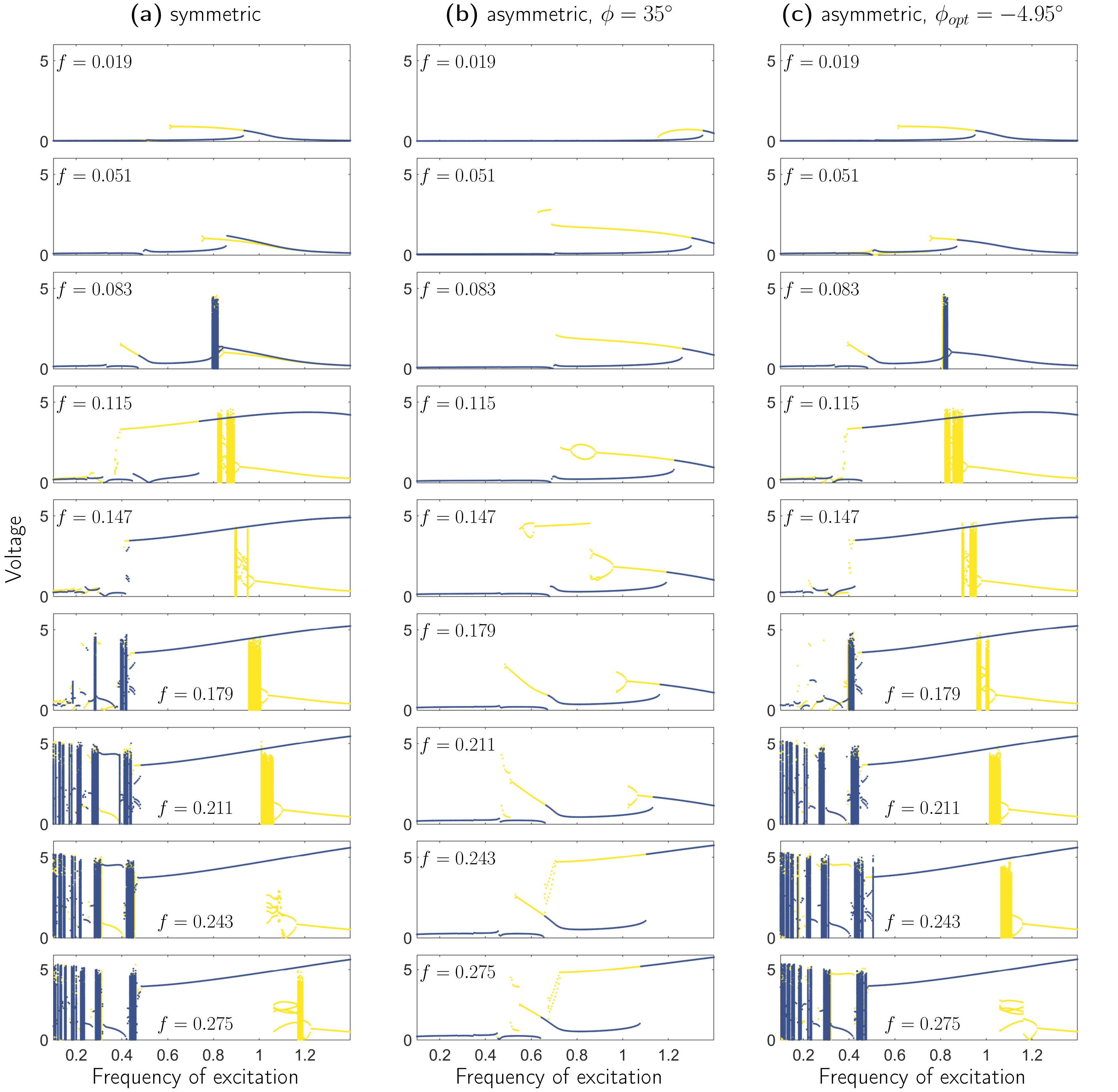}
\caption{Bifurcation diagram of output voltage in bistable energy harvesters as a function of the excitation frequency. The forward sweep is displayed in blue and the backward sweep in yellow. (a) symmetric harvester, (b) asymmetric harvester with $\phi = 35^\circ$, (c) asymmetric harvester with $\phi_{opt} = -4.95^\circ$.}
\label{fig:plot3_frequency}
\end{figure}

Bifurcation diagrams of bistable energy harvesters by examining the voltage response as the frequency is swept for both symmetric and asymmetric models (with high asymmetry, $\phi=35^\circ$, and optimal asymmetry, $\phi_{opt}=-4.95^\circ$) is investigated. The results are presented in Fig.~\ref{fig:plot3_frequency}, with the forward sweeping depicted in blue and the backward sweeping shown in yellow.

In terms of chaos and period multiplicity, the symmetric model (see Fig.~\ref{fig:plot3_frequency}a) exhibits a rich behavior in both forward and backward sweeping. Both sweepings show periodic responses that nearly overlap when f is less than 0.051. However, for f = 0.083 and higher, both bifurcation diagrams reveal a progression to chaotic patterns. For frequencies lower than 0.4, the forward and backward bifurcations occur simultaneously. The forward sweeping of the diagrams indicates higher output regular voltages for higher frequencies, whereas the downward sweeping reveals chaotic zones. Moreover, at amplitudes around f = 0.115 and 0.147, both forward and backward diagrams show regions with discontinuities in the early frequencies. For $\mathnormal{f}\geq0.179$, there are regions with chaos incidence and narrow regions with multiple periods, both in forward and backward sweeping.

The asymmetric model with a slope angle of $\phi = 35^\circ$ (see Fig.~\ref{fig:plot3_frequency}b) exhibited a significant impact of the strong asymmetry on the energy harvesting process. Unlike the symmetric model, the asymmetric model did not display chaotic regions, and the dynamic behavior became periodic for all excitation values. The voltage output increased with increasing frequencies. The forward and backward bifurcations coincided for $\Omega$ values below 0.4 and above 1.2. The backward curve showed bifurcation points within this frequency range only for periodicity responses with higher voltage output. This suggests that the system can change its motion to higher energy orbits. However, it does not prove that the system performs an inter-well motion as the symmetric system.

The asymmetric model with the optimal sloping angle, $\phi_{opt} = -4.95^\circ$ (Fig.\ref{fig:plot3_frequency}c) showed that the chaotic regions re-appeared for frequencies near $0.8$ and low frequencies with high amplitudes. The high voltage output was seen for amplitudes higher than $0.115$ and frequencies higher than $0.4$. Backward and forward bifurcations indicate chaotic motion for frequencies lower than $0.5$ simultaneously. Chaotic motions and multi-periodicity response can reappear at a higher frequency for the backward curve. The asymmetric behavior with optimal angle is similar to the symmetric case's behavior presented in Fig.\ref{fig:plot3_frequency}a. The chaoticity and bifurcation points occur at the same frequency values. Therefore these results demonstrate that the optimal angle can compensate for the influence of asymmetry and behave as a symmetric system.

\subsection{Slope angle analysis}
\label{subsec_angle_analysis}

\begin{figure}
\centering
\includegraphics[width = 1\textwidth]{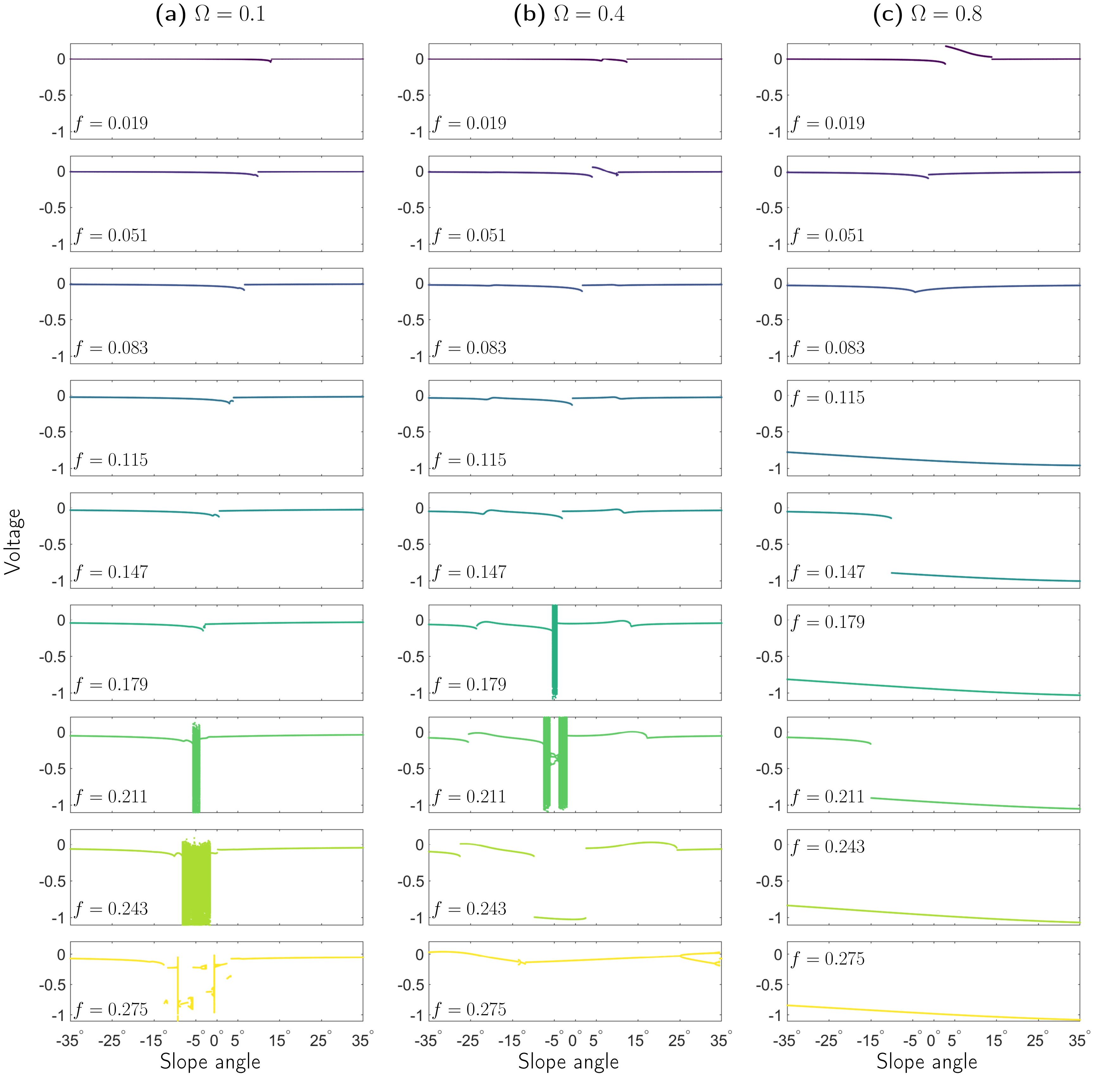}
\caption{Bifurcation diagram of output voltage in asymmetric bistable energy harvesters as a function of slope angle under different amplitude of excitation for (a) $\Omega=0.1$, (b) $\Omega=0.4$, (c) $\Omega=0.8$.}
\label{fig:plot3_angle}
\end{figure}

Bifurcation diagrams for the asymmetric bistable energy harvesters as a function of the slope angle while keeping the asymmetric quadratic term fixed at $\delta=0.15$ are investigated. The study examines the effect of sweeping the slope angle from $-35^\circ$ to $35^\circ$ on the voltage under different excitation conditions, ranging from low to high input energy. The study results are presented in the bifurcation diagrams in Fig.~\ref{fig:plot3_angle}.

Fig.\ref{fig:plot3_angle}a displays the bifurcation diagrams when the excitation frequency is 0.1 under various amplitudes of excitation (from low energy, $f=0.019$, to high energy, $f=0.275$). In most cases, the response is periodic with a single period. However, when the amplitude exceeds $0.211$, a chaotic region emerges around the angle of $-5^\circ$. A similar behavior is observed when the frequency is increased to 0.4, Fig.\ref{fig:plot3_angle}b, but the chaotic region appears at a lower amplitude. Nonetheless, when the amplitude reaches $0.243$, the chaotic motion disappears, and the response becomes periodic. Finally, when the excitation frequency is $0.8$, Fig.~\ref{fig:plot3_angle}c, the system exhibits a periodic behavior for all tested excitation amplitudes.

It is worth noting that the slope angle significantly affects the energy harvesting process, especially when it is close to the optimal slope angle. This observation is consistent with the findings of \cite{norenbergNoDy_2022}. Therefore, careful consideration of the slope angle is necessary when studying the behavior of asymmetric bistable energy harvesting systems.

\section{Basins of attraction}
\label{basins_sect}

This section presents a numerical study of the harvester's dynamics to identify the different attractors that can coexist with the underlying system. The study integrates the dynamics for various initial conditions and maps the possible attractors that the phase space trajectory can accumulate using the 0-1 test for chaos. Firstly we described the methodology to classify the dynamic behavior using the 0-1 test for chaos. Then, we present the results for symmetric and asymmetric models. Finally, a comparison of the relative area of the basins is investigated.

We propose a novel approach for determining the basins of attraction using the 0-1 test for chaos \cite{Gottwald_test}. Since its introduction, the 0-1 test for chaos has been widely used in various fields, such as biology \cite{Toker_2020}, engineering \cite{SAVI2017307}, and finance \cite{Xin_2013}. This tool is based on time series analysis and has a superior computational performance compared to the Lyapunov exponent method. The 0-1 test is computationally more efficient as it focuses on statistical divergences, while the Lyapunov exponent requires phase space reconstruction and eigenvalue calculation, which can be a computationally intensive task. Although the Lyapunov exponent provides a more informative analysis by defining the strength of chaos, determining bifurcation points and dimensionality of the attractor, the 0-1 test is more general and cost-effective. Several works have compared these methods \cite{Lok_2016,Bernardini_2016}, and some have critically studied them \cite{Hu_2005,Romero-Bastida_2009}.

\begin{figure}
    \centering
    \begin{subfigure}[b]{0.48\textwidth}
         \centering
         \includegraphics[width=0.95\textwidth]{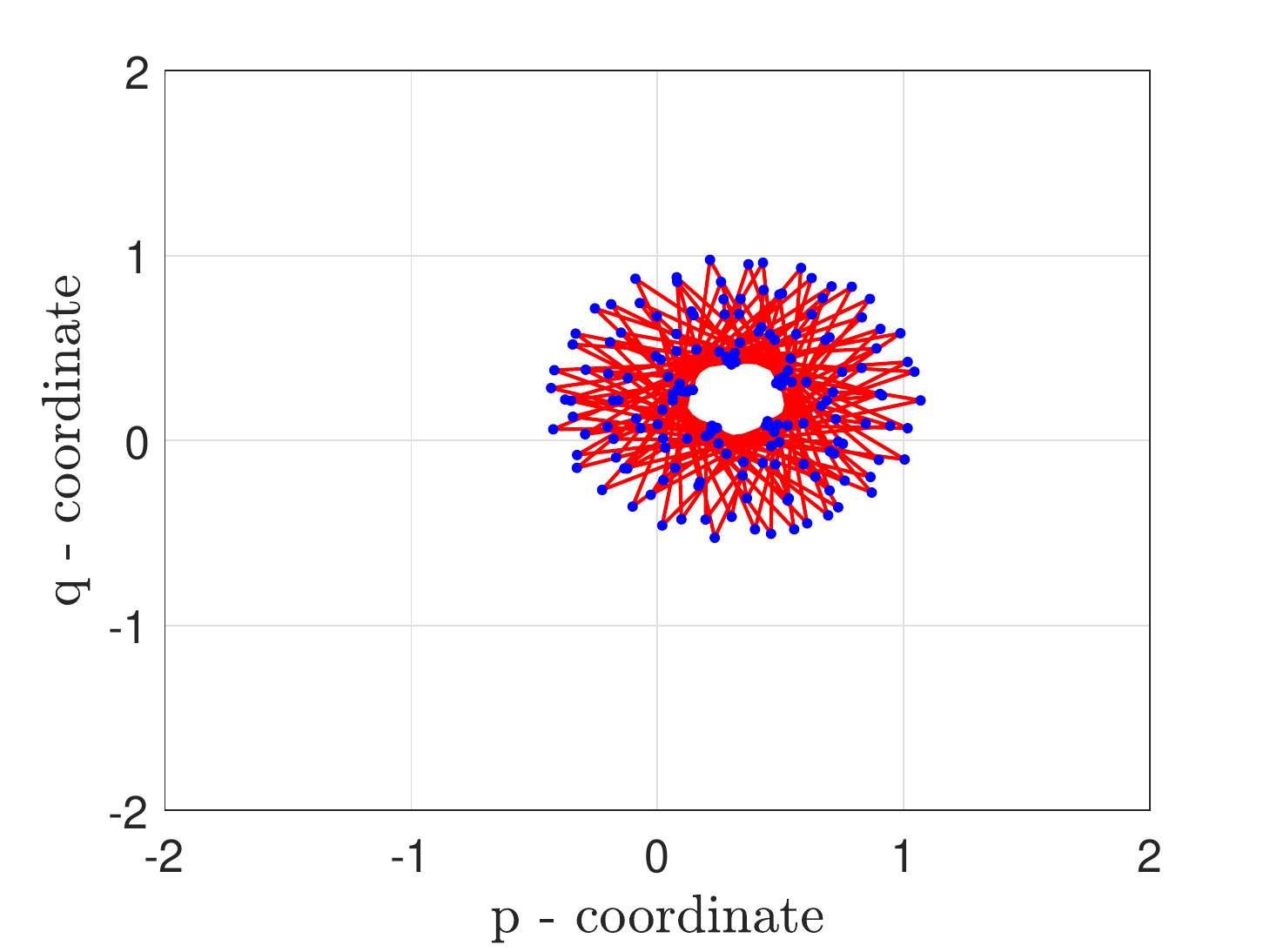}
         \caption{bounded projection}
     \end{subfigure}
     \begin{subfigure}[b]{0.48\textwidth}
         \centering
         \includegraphics[width=0.95\textwidth]{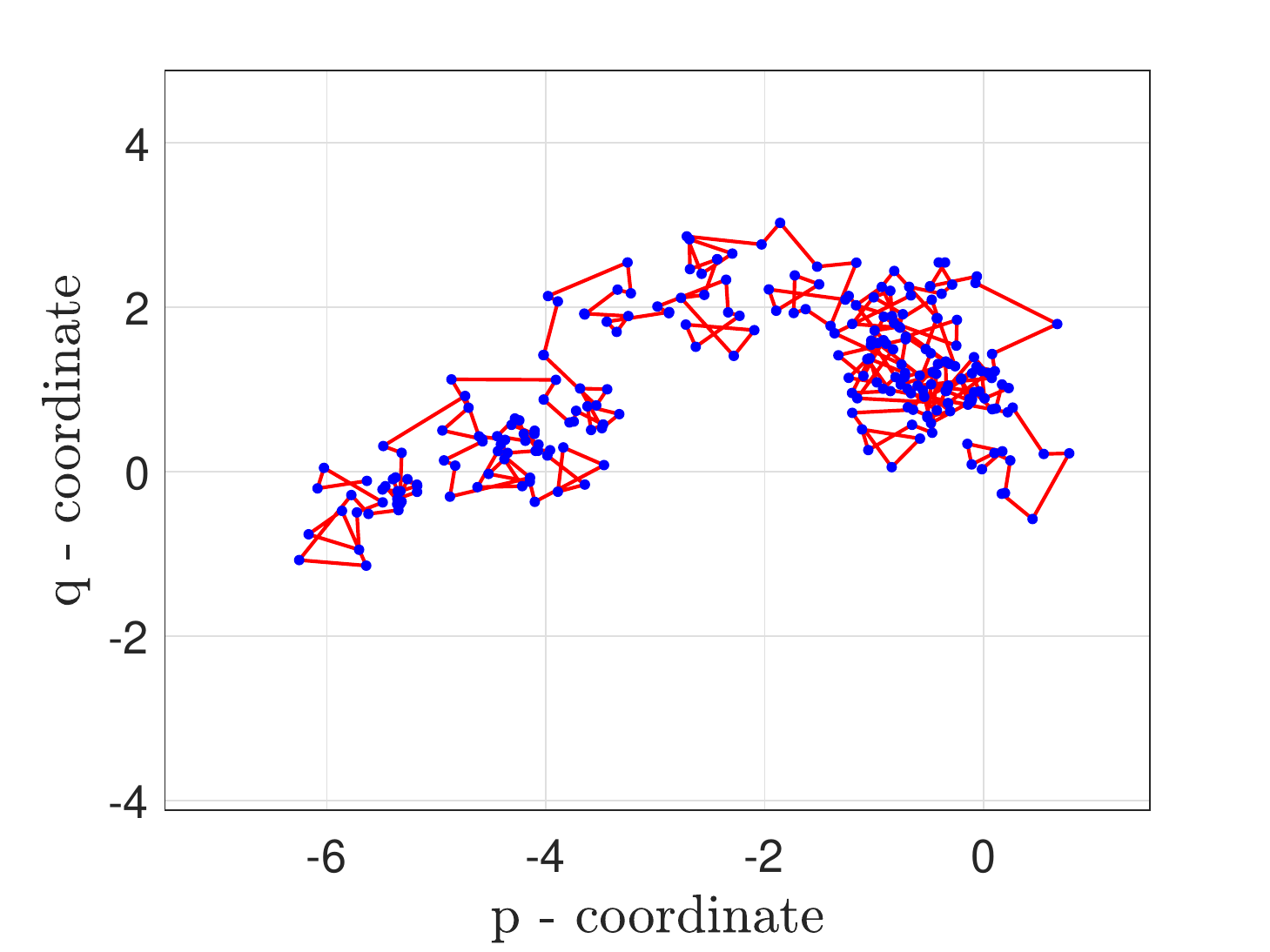}
         \caption{diffuse projection}
     \end{subfigure} 
    \caption{Dynamic extension in the $p \times q$ phase space. The bounded projection displays for periodic response under an excitation amplitude of $\mathnormal{f}=0.115$, while the diffuse projection corresponds to a chaotic response under an excitation amplitude of $\mathnormal{f}=0.083$.} 
    \label{fig:pxq_space}
\end{figure}

The 0-1 test involves mapping the time series $x(t)$ into a pair of transformed coordinates 
\begin{equation}
    p_n(c) = \sum_{j = 1}^n x(t_j) \, \cos{(j\,c)} \, ,
\end{equation}
\begin{equation}
    q_n(c) = \sum_{j = 1}^n x(t_j) \, \sin{(j\,c)}\, ,
\end{equation}
which $c$ is a random value uniformly drawn in the support $[0,2\pi)$ and $n=1,2,\dots,N$. For regular behavior, the coordinates $p_n$ and $q_n$ exhibit bounded motion, while for chaotic behavior, they asymptotically behave as Brownian motion. The time series projections in the $p \times q$ space are shown in Fig.~\ref{fig:pxq_space}a (periodic response) and \ref{fig:pxq_space}b (chaotic response) for two different values of the excitation amplitude. 

\begin{figure}
    \centering
    \includegraphics[width=0.85\textwidth]{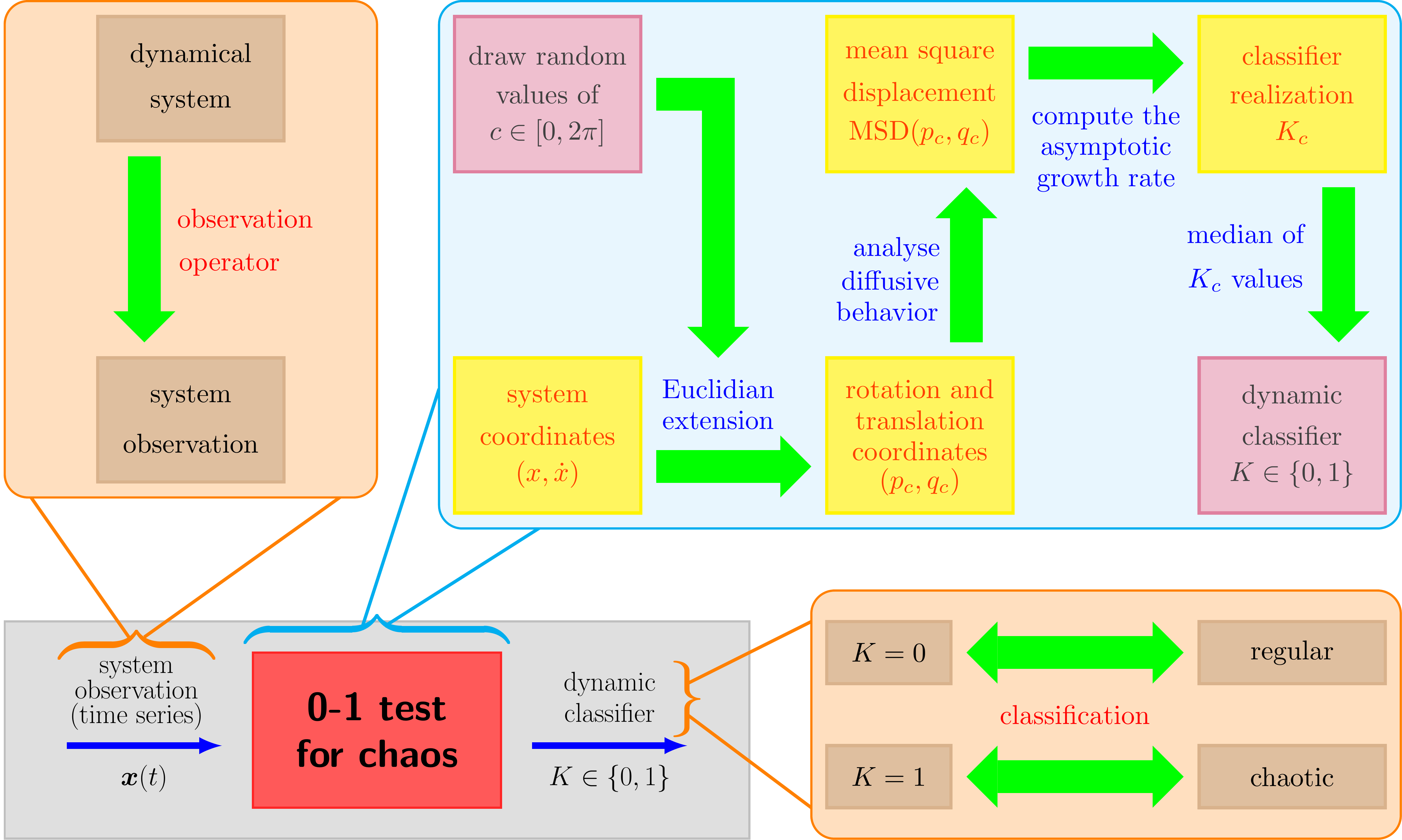}
    \caption{Schematic illustration of the 0-1 test for chaos classification. From the time series, the test consists of transforming the original dynamics ($x$ and $\dot{x}$) into Euclidean extension coordinates ($p_c$ and $q_c$) to calculate the mean square displacement ($MSD(p_c,q_c)$). Then, the statistical classifier $K_c$ is obtained. The median of $K_c$ for several realizations of $c$ random values, named $K$, is computed. Finally, if $K=0$, the dynamic behavior is classified as regular, whereas if $K=1$, it is classified as chaotic.}
    \label{fig:test01}
\end{figure}

The mean square deviation
\begin{equation}
        M_n(c) = \lim_{N\rightarrow \infty} \frac{1}{N} \sum_{j = 1}^N \left( \left[p_{j+n}(c) - p_j(c) \right]^2 +  \right. \\ \left. \left[q_{j+n}(c) - q_j(c)\right]^2 \right) \,
\end{equation}
is used to analyze diffusive or non-diffusive behavior. The classifier
\begin{equation}
    K_c = \lim_{N\rightarrow \infty} \frac{Cov(\overline{t_n},\overline{M_n})}{\sqrt{Var(\overline{t_n})\,Var(\overline{M_n})}},
\end{equation}
where $\overline{M_n} = (M_1,M_2,\dots,M_n)$, $\overline{t_n} = (t_1,t_2,\dots,t_n)$, $Cov$ and $Var$ are covariance and variance operators, is then calculated for several values of $c$, and the median is obtained. The dynamic behavior is classified as chaotic ($K > 0.8$), regular ($K < 0.2$), or inconclusive ($0.2 < K < 0.8$) as shown in Figure~\ref{fig:test01}. The convergence of the 0-1 test over the c sampling process and its effectiveness compared to the Lyapunov exponents are proved and demonstrated in the Supplementary Material.

\begin{figure}
	\centering
    \includegraphics[width=0.85\textwidth]{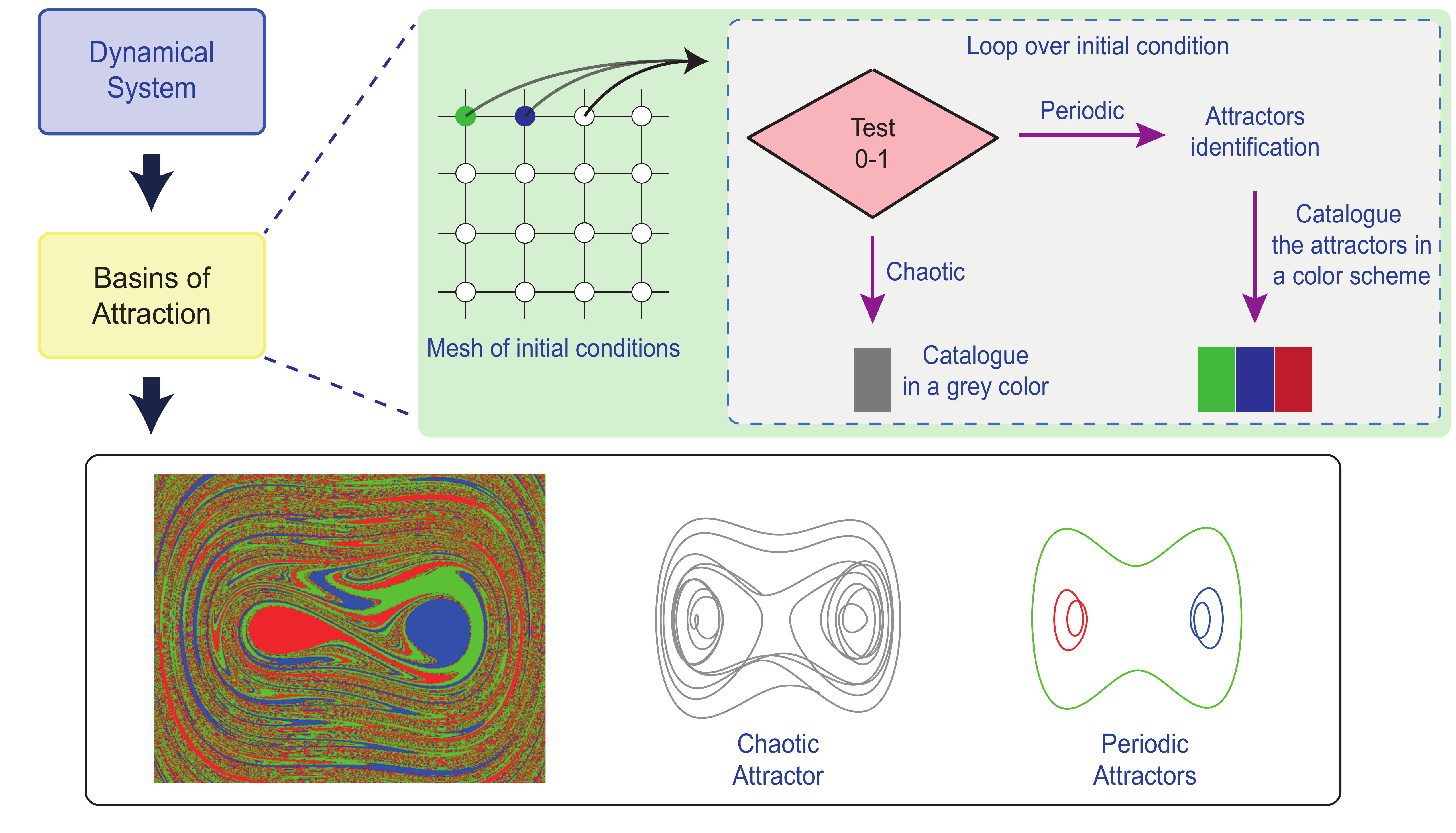}
	\caption{Diagram of the process for computing the basins of attraction. The 0-1 test for chaos is applied to each initial condition. If the motion is chaotic, the grid point is assigned a gray color. If the motion is regular, the next step is to identify the attractor, and the grid point is assigned a color representing a different attractor.}
	\label{fig:basins_cal_fig}
\end{figure}

This process is repeated for a grid of initial conditions to obtain the basins of attraction, shown in Figure~\ref{fig:basins_cal_fig}. If the dynamic behavior is classified as chaotic, the node is colored grey. If the dynamic behavior is classified as regular, the attractors are further identified and distinguished using different colors. The green color indicates a regular high-energy orbit, which is the optimal scenario for harvesting, while the red and blue colors indicate low-energy regular orbits. Other colors, such as magenta and cyan, indicate different regular behaviors.

\subsection{Symmetric bistable energy harvester}


A symmetric bistable energy harvester is studied. The investigation examines two excitation amplitudes, $f = 0.083$ and $f=0.115$, at varying excitation frequencies. Next, the study investigates various excitation amplitudes for a forcing frequency of $0.8$.

\begin{figure}
    \centering
    \begin{subfigure}[b]{0.31\textwidth}
         \centering
         \includegraphics[width=1\textwidth]{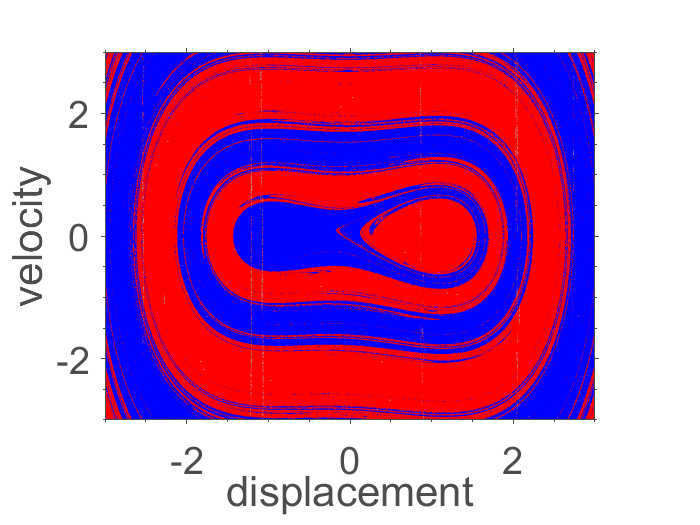}
         \caption{$\Omega = 0.1$}
     \end{subfigure}
     \begin{subfigure}[b]{0.31\textwidth}
         \centering
         \includegraphics[width=\textwidth]{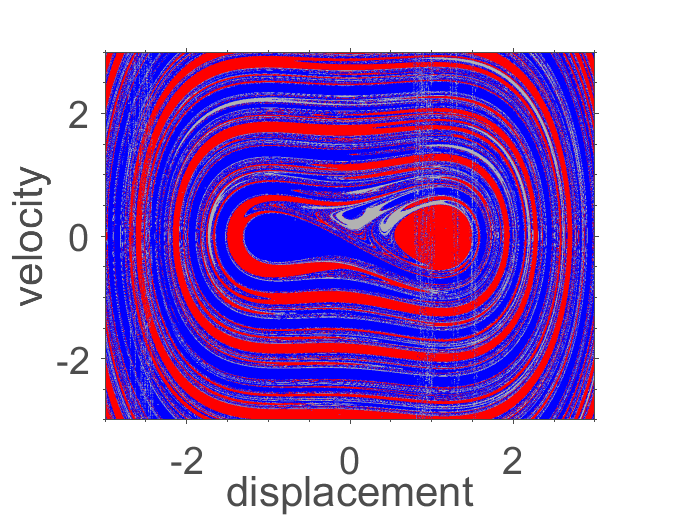}
         \caption{$\Omega = 0.2$}
     \end{subfigure}
     \begin{subfigure}[b]{0.31\textwidth}
         \centering
         \includegraphics[width=\textwidth]{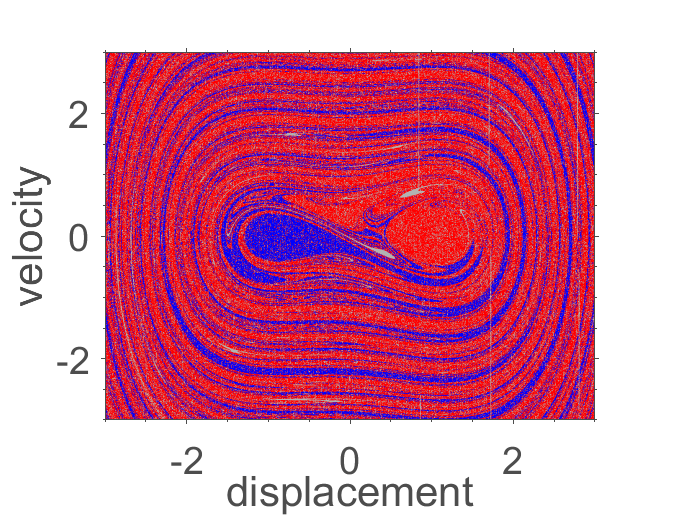}
         \caption{$\Omega = 0.3$}
     \end{subfigure}
     \begin{subfigure}[b]{0.31\textwidth}
         \centering
         \includegraphics[width=\textwidth]{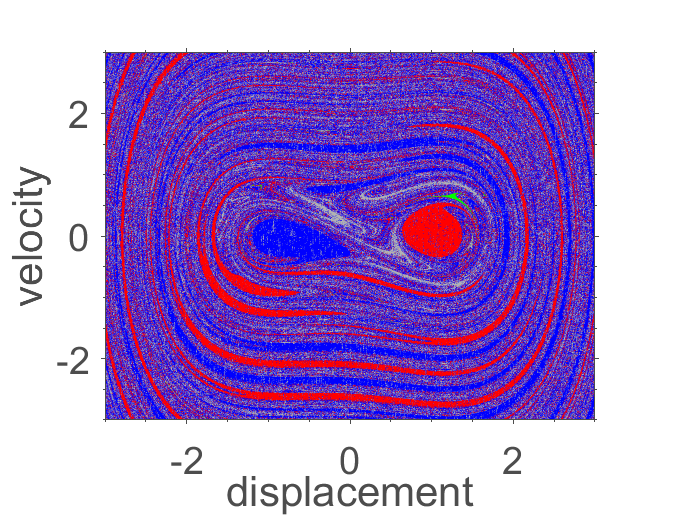}
         \caption{$\Omega = 0.4$}
     \end{subfigure}
     \begin{subfigure}[b]{0.31\textwidth}
         \centering
         \includegraphics[width=\textwidth]{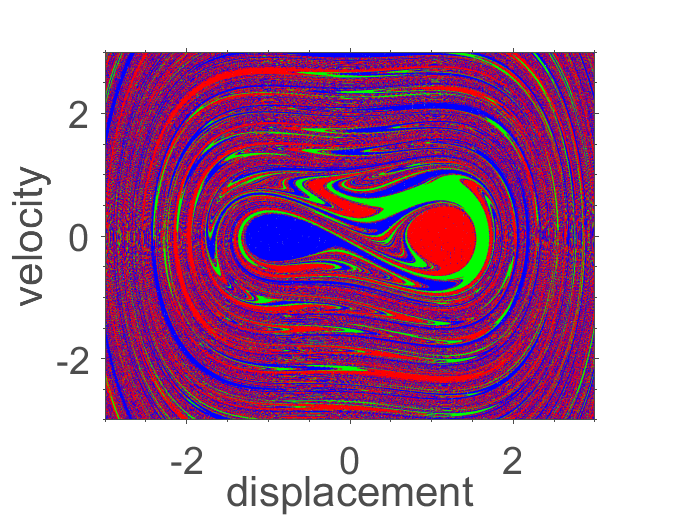}
         \caption{$\Omega = 0.5$}
     \end{subfigure}
     \begin{subfigure}[b]{0.31\textwidth}
         \centering
         \includegraphics[width=\textwidth]{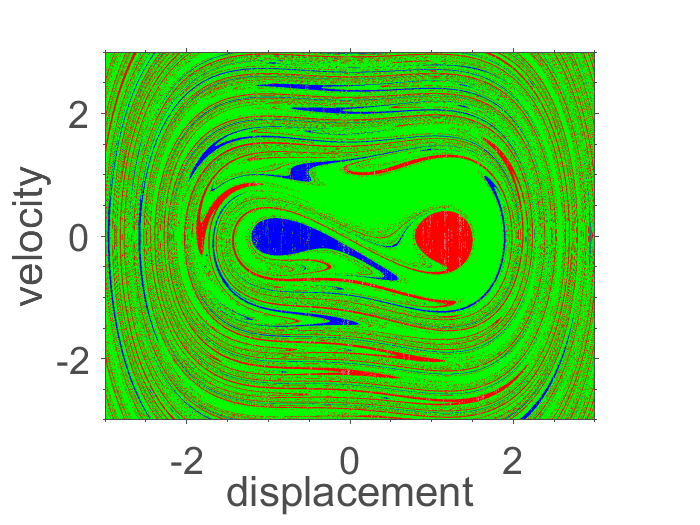}
         \caption{$\Omega = 0.6$}
     \end{subfigure}
     \begin{subfigure}[b]{0.31\textwidth}
         \centering
         \includegraphics[width=\textwidth]{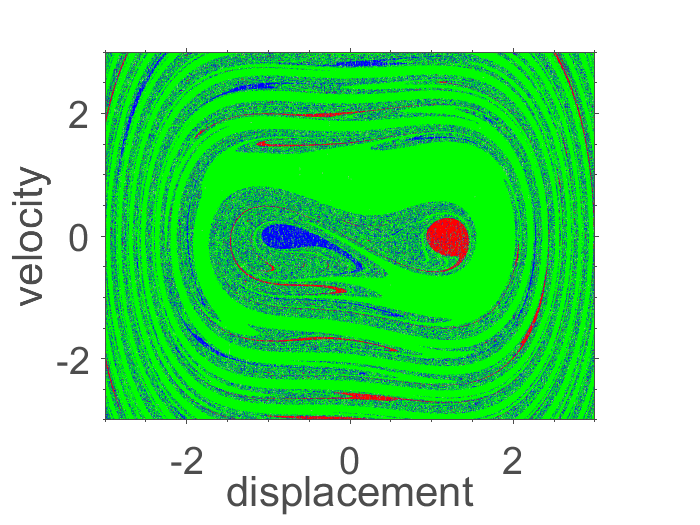}
         \caption{$\Omega = 0.7$}
     \end{subfigure}
     \begin{subfigure}[b]{0.31\textwidth}
         \centering
         \includegraphics[width=\textwidth]{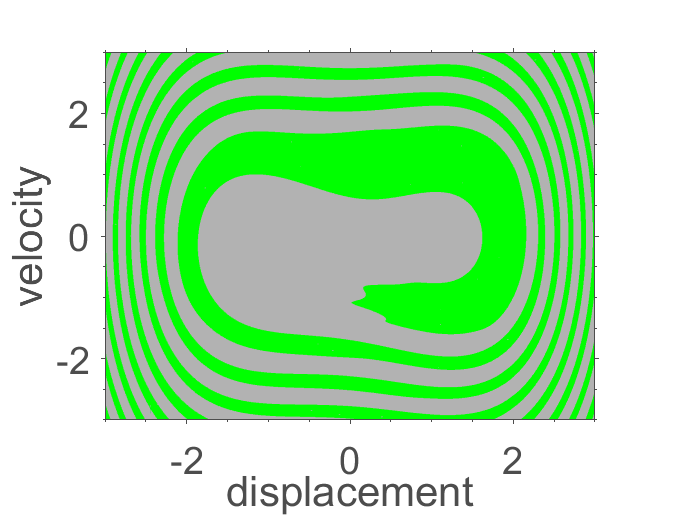}
         \caption{$\Omega = 0.8$}
     \end{subfigure}
     \begin{subfigure}[b]{0.31\textwidth}
         \centering
         \includegraphics[width=\textwidth]{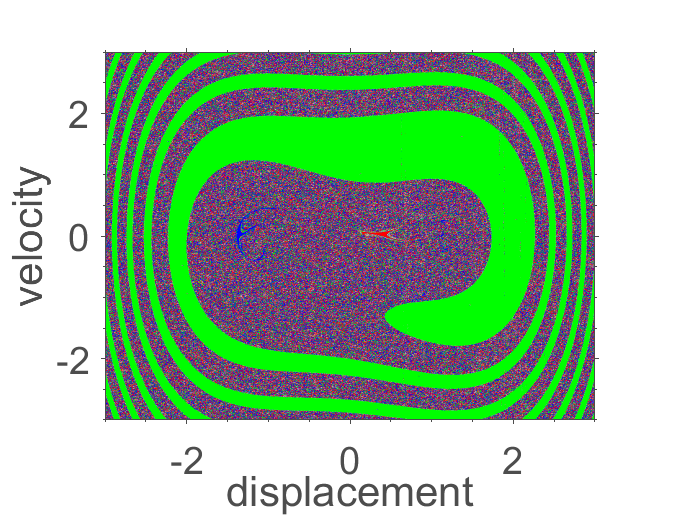}
         \caption{$\Omega = 0.9$}
     \end{subfigure}

    \caption{Basins of attraction intersection with the plane $v=0$ for the symmetric model with a forcing amplitude of $f=0.083$ under forcing frequency from 0.1 to 0.8, within the bounds of $-3 \leq x_0 \leq 3$ and $-3 \leq \dot{x}_0 \leq 3$. The corresponding attractor colors are shown in Fig.\ref{fig:atractors_fig1}.}
    \label{fig:basins_f83}
\end{figure}
\begin{figure}[h]
	\centering
    
    \begin{subfigure}[b]{0.31\textwidth}
         \centering
         \includegraphics[width=\textwidth]{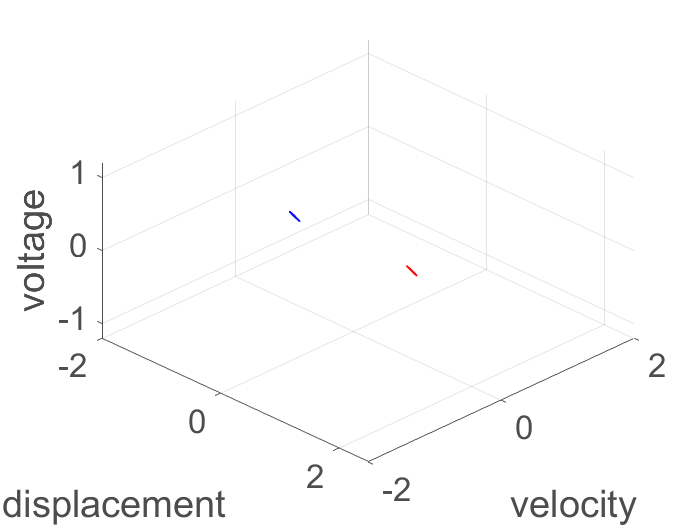}
         \caption{$\Omega = 0.1$}
     \end{subfigure}
     \begin{subfigure}[b]{0.31\textwidth}
         \centering
         \includegraphics[width=\textwidth]{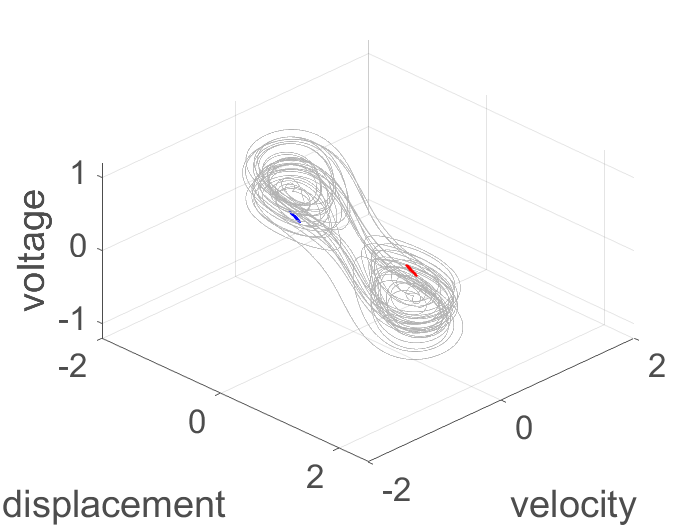}
         \caption{$\Omega = 0.2$}
     \end{subfigure}
     \begin{subfigure}[b]{0.31\textwidth}
         \centering
         \includegraphics[width=\textwidth]{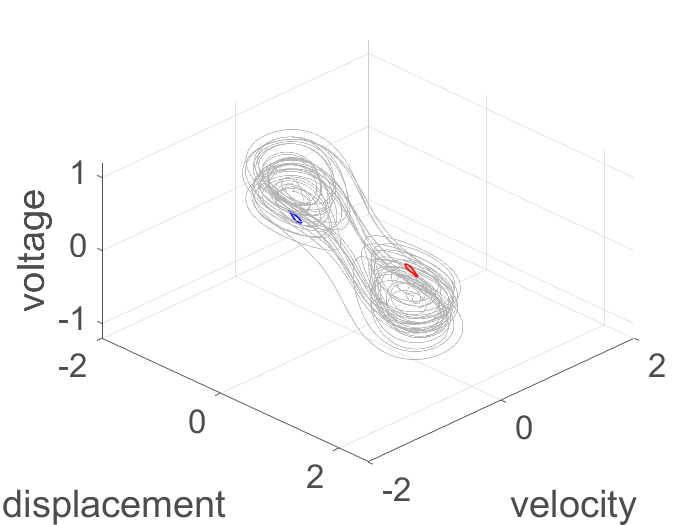}
         \caption{$\Omega = 0.3$}
     \end{subfigure}
     \begin{subfigure}[b]{0.31\textwidth}
         \centering
         \includegraphics[width=\textwidth]{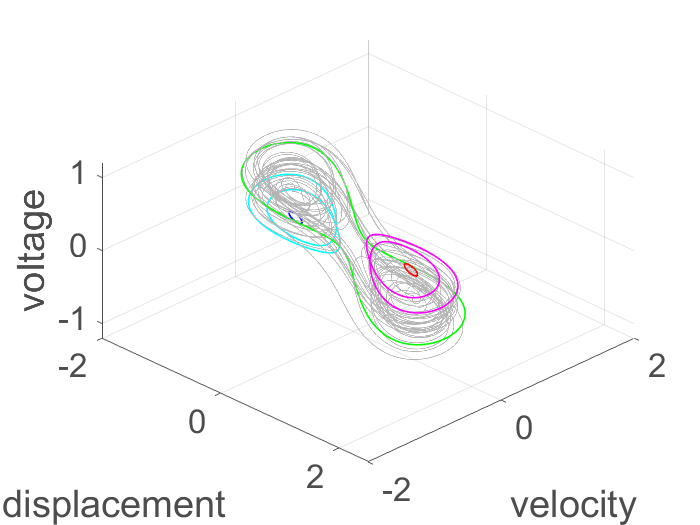}
         \caption{$\Omega = 0.4$}
     \end{subfigure}
     \begin{subfigure}[b]{0.31\textwidth}
         \centering
         \includegraphics[width=\textwidth]{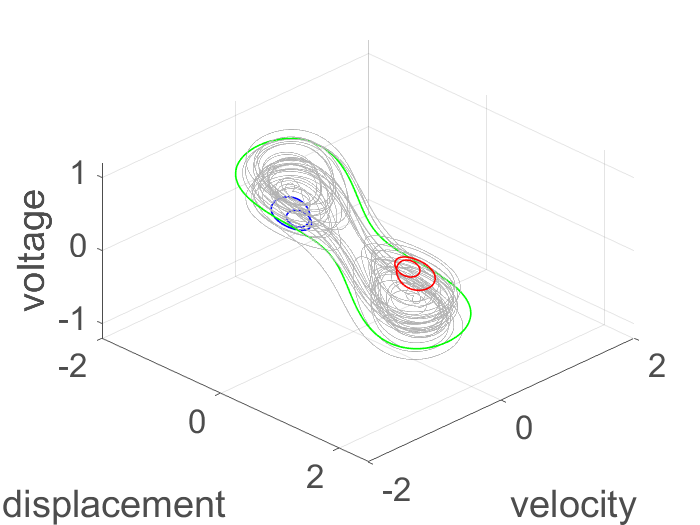}
         \caption{$\Omega = 0.5$}
     \end{subfigure}
     \begin{subfigure}[b]{0.31\textwidth}
         \centering
         \includegraphics[width=\textwidth]{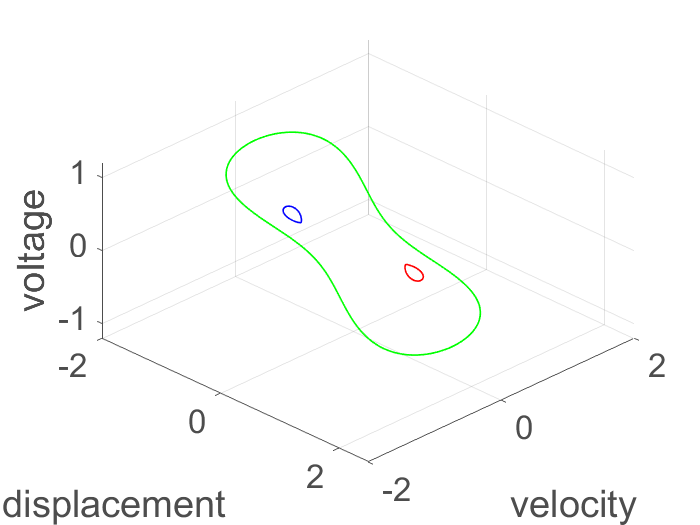}
         \caption{$\Omega = 0.6$}
     \end{subfigure}
     \begin{subfigure}[b]{0.31\textwidth}
         \centering
         \includegraphics[width=\textwidth]{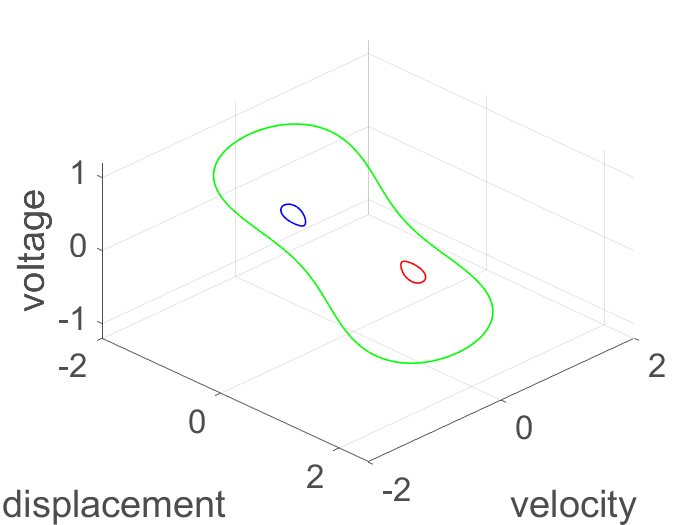}
         \caption{$\Omega = 0.7$}
     \end{subfigure}
     \begin{subfigure}[b]{0.31\textwidth}
         \centering
         \includegraphics[width=\textwidth]{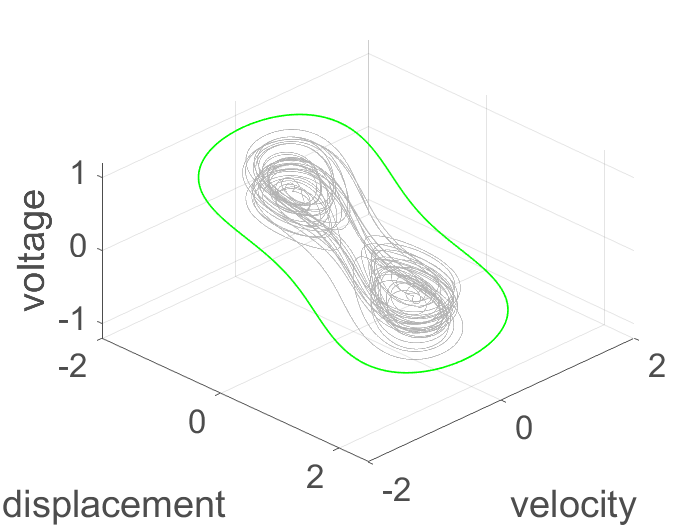}
         \caption{$\Omega = 0.8$}
     \end{subfigure}
     \begin{subfigure}[b]{0.31\textwidth}
         \centering
         \includegraphics[width=\textwidth]{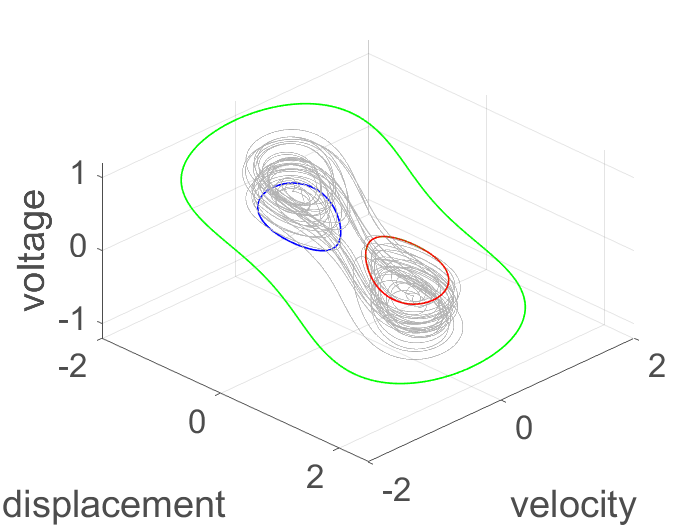}
         \caption{$\Omega = 0.9$}
     \end{subfigure}
	\caption{Attractors of the symmetric model for $f=0.083$ under different forcing frequencies. Attractors correspond to basins of attraction presented in Fig.\ref{fig:basins_f83}.}
	\label{fig:atractors_fig1}
\end{figure}

The Fig. \ref{fig:basins_f83} and \ref{fig:atractors_fig1} present the basins of attraction and the attractors for a fixed excitation amplitude of $\mathnormal{f} = 0.083$, and a range of excitation frequencies from $\Omega=0.1$ to $\Omega=0.9$. The attractors of a dynamical system exhibit different behavior as the excitation frequency $\Omega$ changes. At $\Omega=0.1$, two regular low-energy attractors (blue and red basins) are observed with small displacement and velocity orbits. As $\Omega$ increases, the basins of attraction narrow, showing a more complicated boundary of basins, and chaotic behavior becomes evident. Predominant attractors are identified as red ($\Omega \in {0.3, 0.5}$), blue ($\Omega = 0.4$), and green ($\Omega \in {0.6, 0.7, 0.8}$), with erosions in boundary regions shown as fractals. At $\Omega = 0.4$, there are also cyan and magenta attractors, but they are not prominent in the basins of attraction. Near $\Omega = 0.8$, the basins narrow and move away from equilibrium points, resulting in only high energy, periodic, and chaotic orbits. The system exhibits two regular attractors with low or high-energy orbits and a chaotic attractor.


\begin{figure}
    \centering
    
    \begin{subfigure}[b]{0.31\textwidth}
         \centering
         \includegraphics[width=\textwidth]{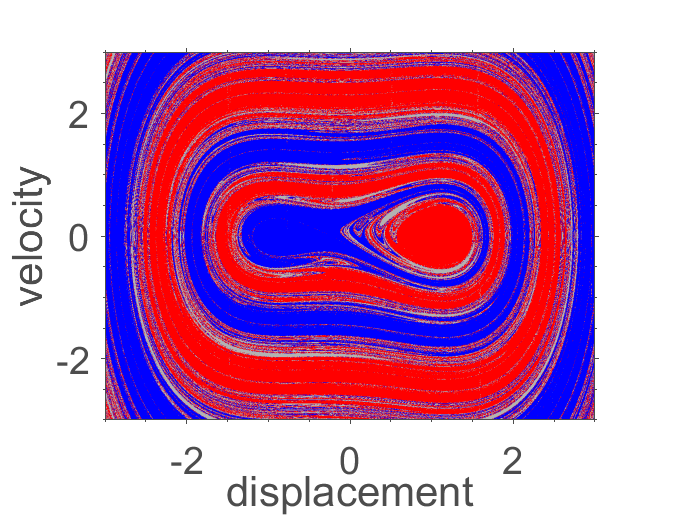}
         \caption{$\Omega = 0.1$}
     \end{subfigure}
     \begin{subfigure}[b]{0.31\textwidth}
         \centering
         \includegraphics[width=\textwidth]{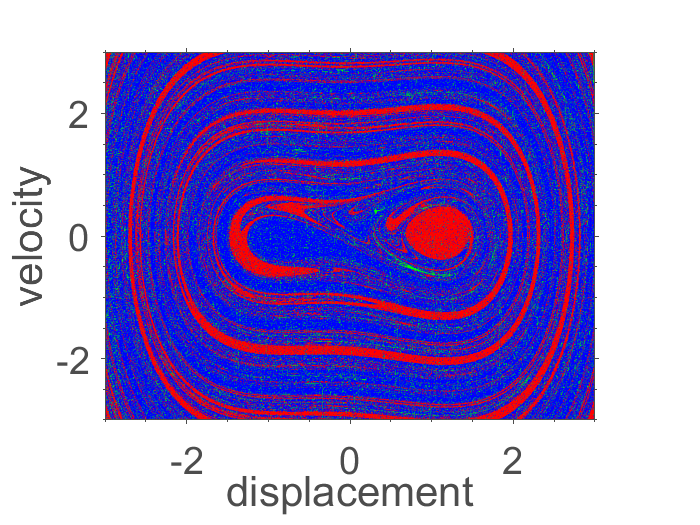}
         \caption{$\Omega = 0.2$}
     \end{subfigure}
     \begin{subfigure}[b]{0.31\textwidth}
         \centering
         \includegraphics[width=\textwidth]{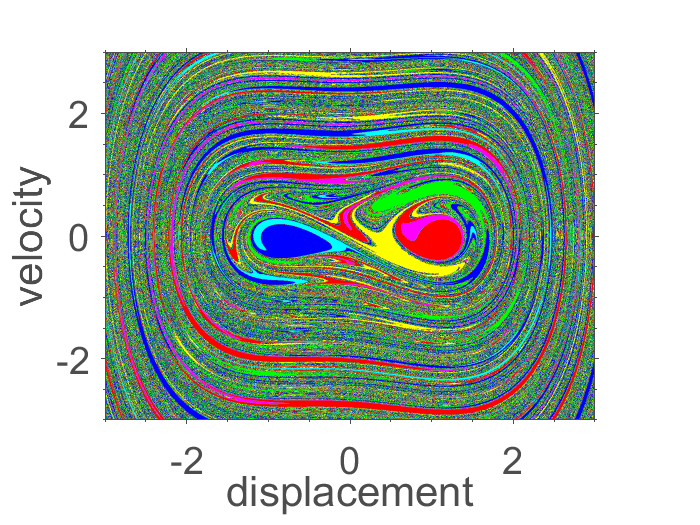}
         \caption{$\Omega = 0.3$}
     \end{subfigure}
     \begin{subfigure}[b]{0.31\textwidth}
         \centering
         \includegraphics[width=\textwidth]{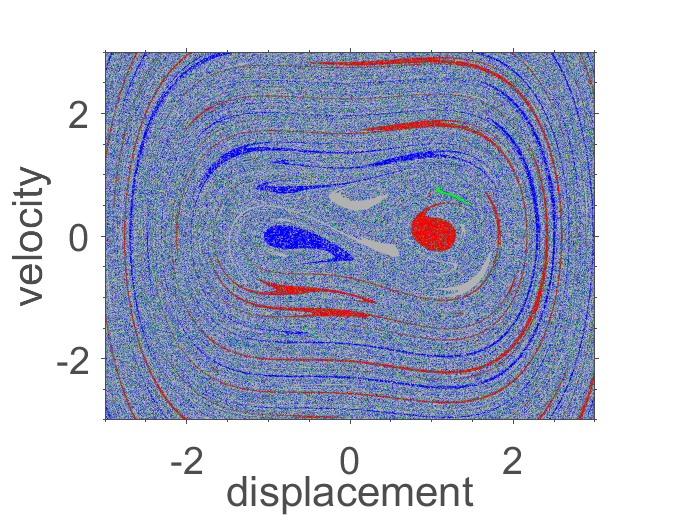}
         \caption{$\Omega = 0.4$}
     \end{subfigure}
     \begin{subfigure}[b]{0.31\textwidth}
         \centering
         \includegraphics[width=\textwidth]{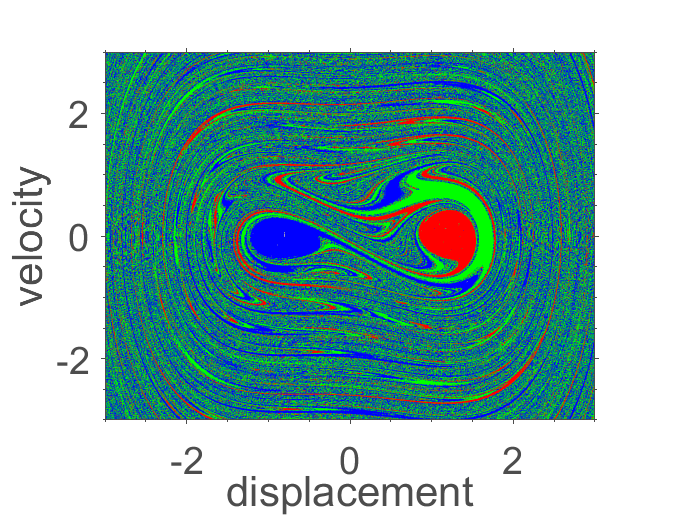}
         \caption{$\Omega = 0.5$}
     \end{subfigure}
     \begin{subfigure}[b]{0.31\textwidth}
         \centering
         \includegraphics[width=\textwidth]{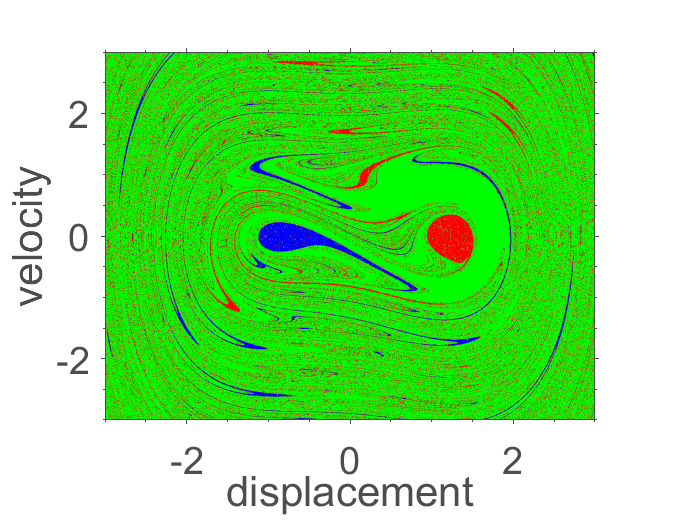}
         \caption{$\Omega = 0.6$}
     \end{subfigure}
     \begin{subfigure}[b]{0.31\textwidth}
         \centering
         \includegraphics[width=\textwidth]{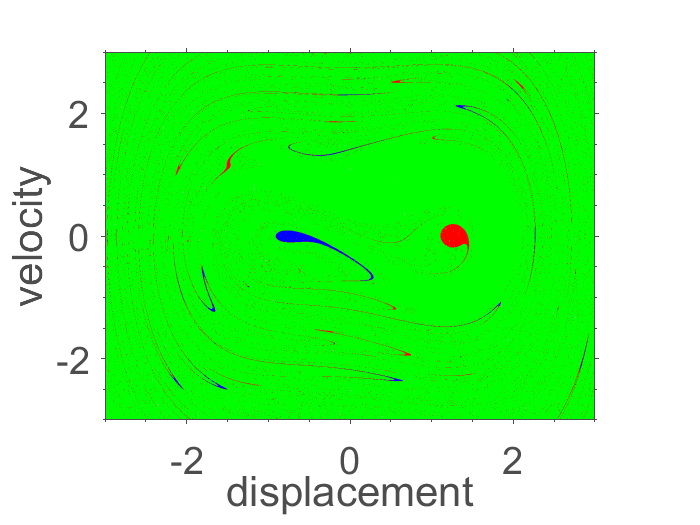}
         \caption{$\Omega = 0.7$}
     \end{subfigure}
     \begin{subfigure}[b]{0.31\textwidth}
         \centering
         \includegraphics[width=\textwidth]{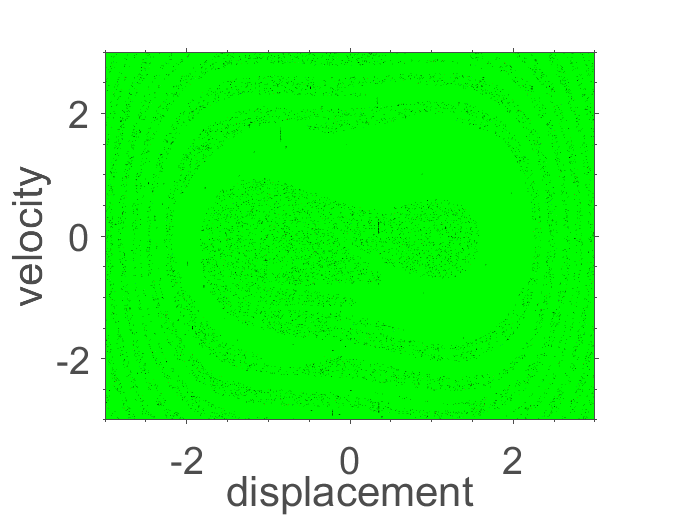}
         \caption{$\Omega = 0.8$}
     \end{subfigure}
     \begin{subfigure}[b]{0.31\textwidth}
         \centering
         \includegraphics[width=\textwidth]{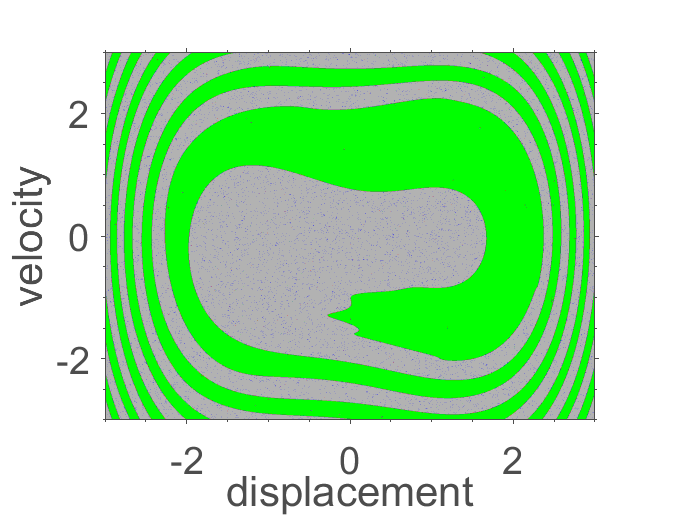}
         \caption{$\Omega = 0.9$}
     \end{subfigure}
    \caption{Basins of attraction intersection with the plane $v=0$ for the symmetric model with a forcing amplitude of $f=0.115$ under forcing frequency from 0.1 to 0.8, within the bounds of $-3 \leq x_0 \leq 3$ and $-3 \leq \dot{x}_0 \leq 3$. The corresponding attractor colors are shown in Fig.\ref{fig:atractors_fig2}.}
    \label{fig:basins_f115}
\end{figure}

\begin{figure}
	\centering
	
    \begin{subfigure}[b]{0.31\textwidth}
         \centering
         \includegraphics[width=\textwidth]{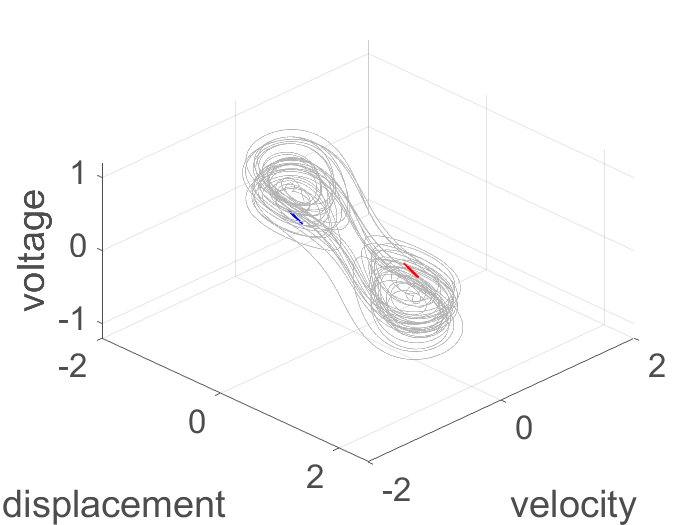}
         \caption{$\Omega = 0.1$}
     \end{subfigure}
     \begin{subfigure}[b]{0.31\textwidth}
         \centering
         \includegraphics[width=\textwidth]{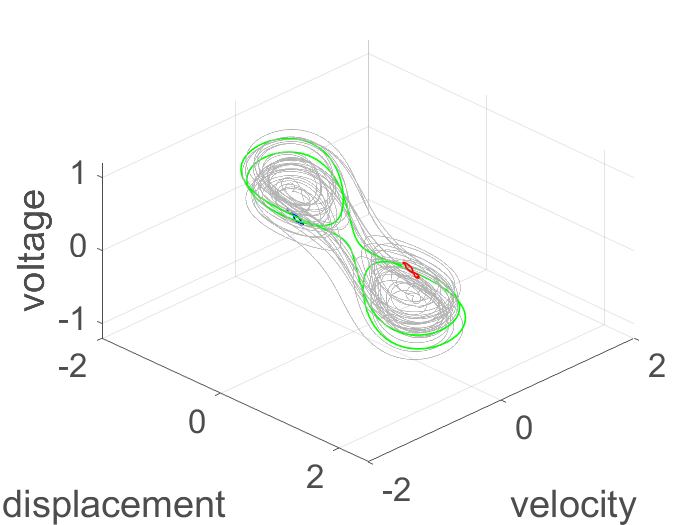}
         \caption{$\Omega = 0.2$}
     \end{subfigure}
     \begin{subfigure}[b]{0.31\textwidth}
         \centering
         \includegraphics[width=\textwidth]{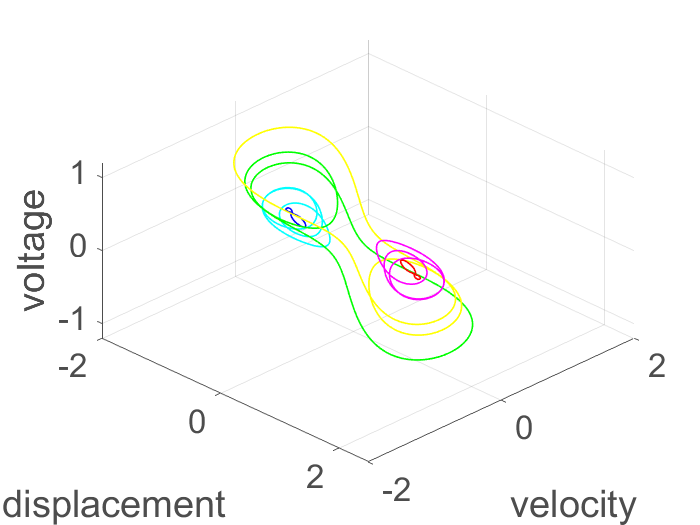}
         \caption{$\Omega = 0.3$}
     \end{subfigure}
     \begin{subfigure}[b]{0.31\textwidth}
         \centering
         \includegraphics[width=\textwidth]{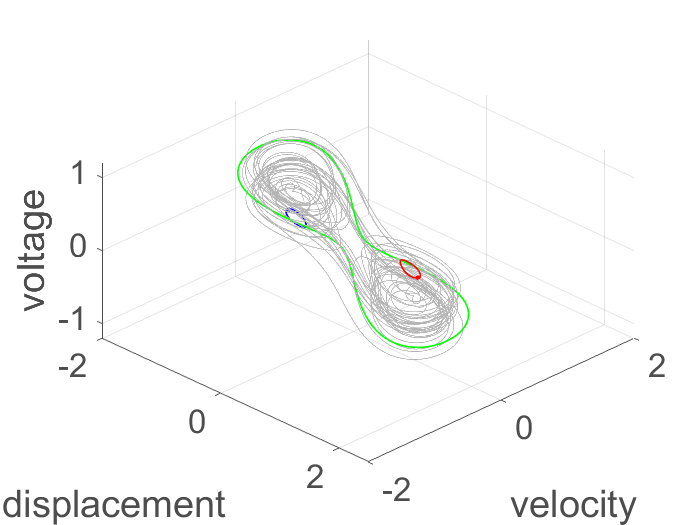}
         \caption{$\Omega = 0.4$}
     \end{subfigure}
     \begin{subfigure}[b]{0.31\textwidth}
         \centering
         \includegraphics[width=\textwidth]{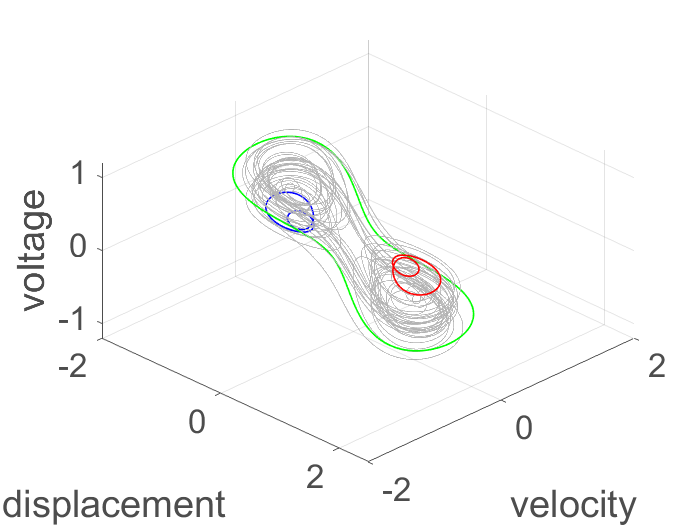}
         \caption{$\Omega = 0.5$}
     \end{subfigure}
     \begin{subfigure}[b]{0.31\textwidth}
         \centering
         \includegraphics[width=\textwidth]{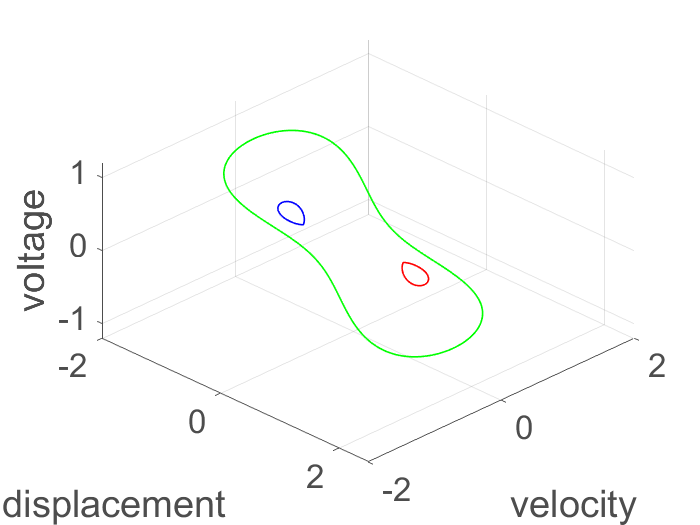}
         \caption{$\Omega = 0.6$}
     \end{subfigure}
     \begin{subfigure}[b]{0.31\textwidth}
         \centering
         \includegraphics[width=\textwidth]{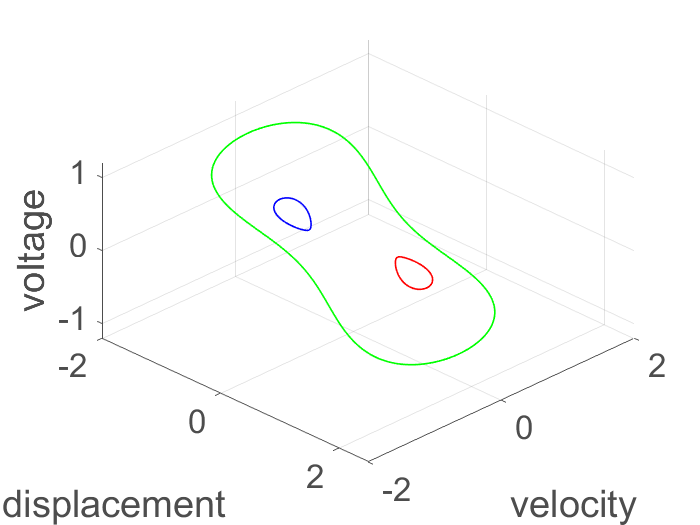}
         \caption{$\Omega = 0.7$}
     \end{subfigure}
     \begin{subfigure}[b]{0.31\textwidth}
         \centering
         \includegraphics[width=\textwidth]{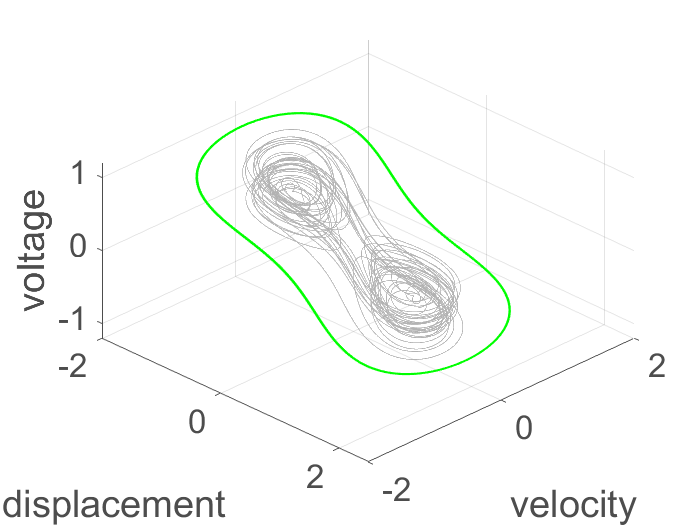}
         \caption{$\Omega = 0.8$}
     \end{subfigure}
     \begin{subfigure}[b]{0.31\textwidth}
         \centering
         \includegraphics[width=\textwidth]{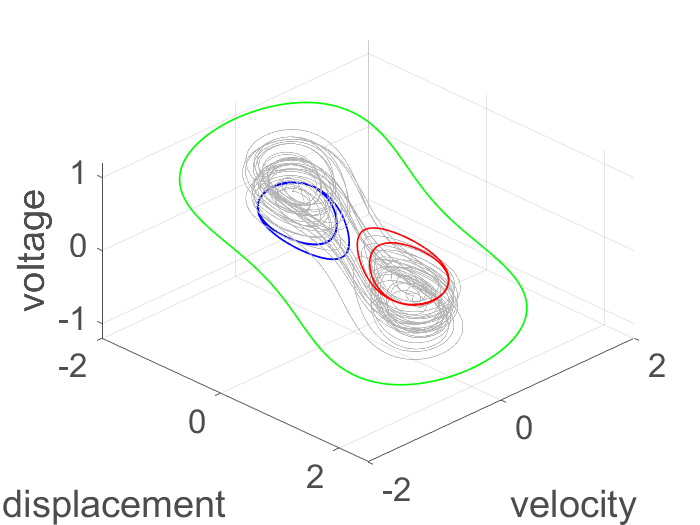}
         \caption{$\Omega = 0.9$}
     \end{subfigure}

	\caption{Attractors of the symmetric model for $f=0.115$ under different forcing frequencies. Attractors correspond to basins of attraction presented in Fig.\ref{fig:basins_f115}.}
	\label{fig:atractors_fig2}
\end{figure}

Fig.~\ref{fig:basins_f115} illustrate the basins of attraction, and Fig.\ref{fig:atractors_fig2} their corresponding attractors for an amplitude of excitation $\mathnormal{f} = 0.115$ and a frequency of excitation $\Omega$ ranging from 0.1 to 0.9. The graphs reveal the presence of chaotic behavior across most excitation frequencies, except for $\Omega = 0.3$, $\Omega = 0.6$, and $0.7$. The low energy orbits, represented by the blue and red attractors, are mainly found at low frequencies. At $\Omega = 0.3$, a yellow basin with high energy orbits and multiple frequencies appears. Additionally, the green attractor emerges at $\Omega = 0.5$, and the blue and red orbits display two frequencies. The green attractor, which features high-energy orbits, becomes increasingly dominant as the frequency increases. At $\Omega = 0.8$, the basins are fully covered by high-energy orbits, dominant by green basins and some scattered chaotic conditions. This frequency achieves the desired outcome for energy harvesting. Finally, at $\Omega = 0.9$, chaos prevails with the green basin, while the blue and red basins with two frequencies are scattered throughout the chaos.

\begin{figure}
    \centering    
    \begin{subfigure}[b]{0.31\textwidth}
         \centering
         \includegraphics[width=\textwidth]{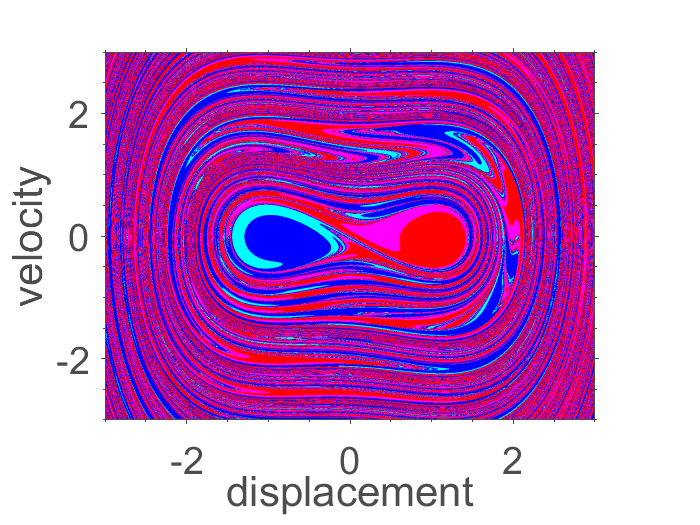}
         \caption{$f = 0.019$}
     \end{subfigure}
     \begin{subfigure}[b]{0.31\textwidth}
         \centering
         \includegraphics[width=\textwidth]{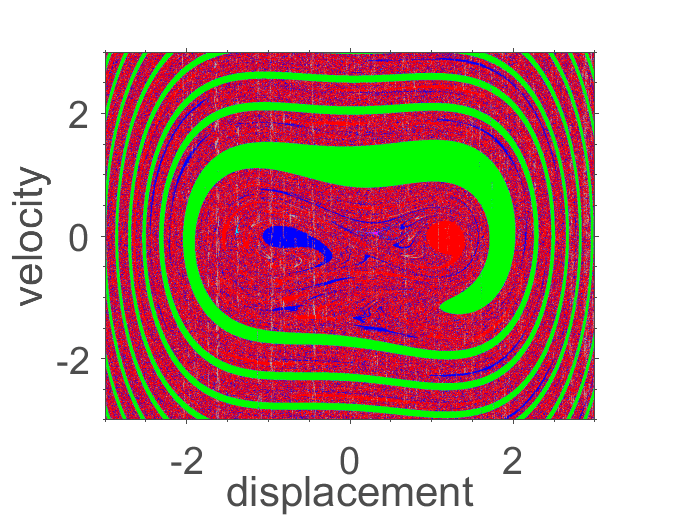}
         \caption{$f = 0.051$}
     \end{subfigure}
     \begin{subfigure}[b]{0.31\textwidth}
         \centering
         \includegraphics[width=\textwidth]{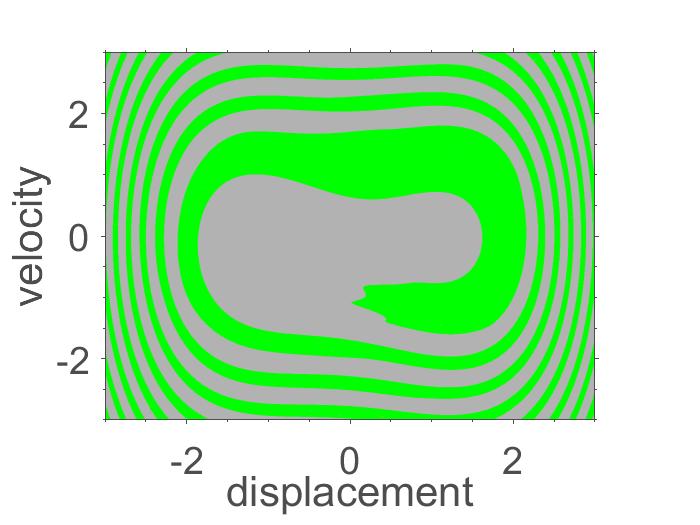}
         \caption{$f = 0.083$}
     \end{subfigure}
     \begin{subfigure}[b]{0.31\textwidth}
         \centering
         \includegraphics[width=\textwidth]{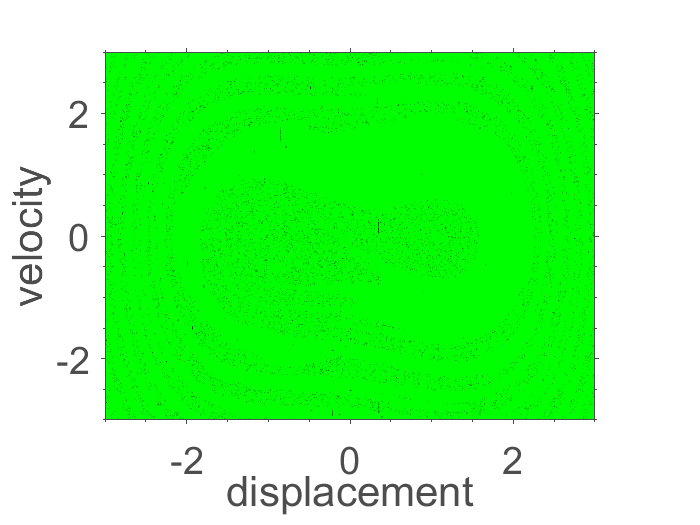}
         \caption{$f = 0.115$}
     \end{subfigure}
     \begin{subfigure}[b]{0.31\textwidth}
         \centering
         \includegraphics[width=\textwidth]{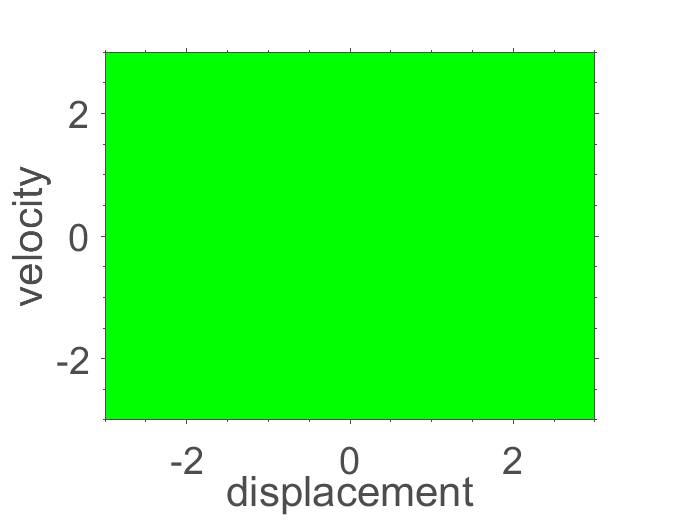}
         \caption{$f = 0.147$}
     \end{subfigure}
     \begin{subfigure}[b]{0.31\textwidth}
         \centering
         \includegraphics[width=\textwidth]{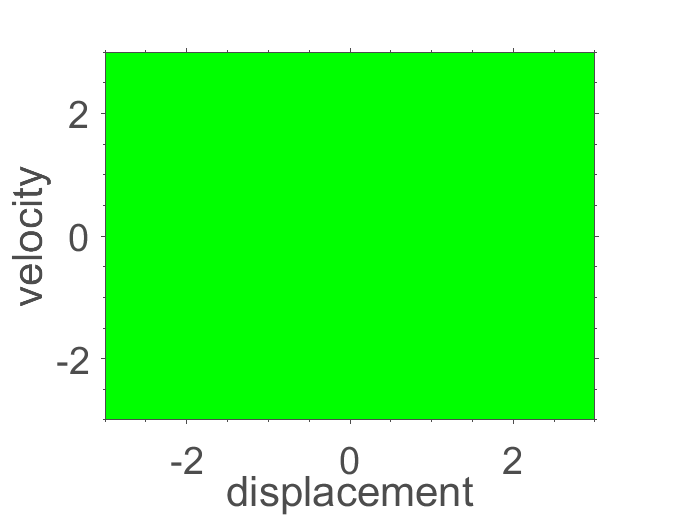}
         \caption{$f = 0.179$}
     \end{subfigure}
     \begin{subfigure}[b]{0.31\textwidth}
         \centering
         \includegraphics[width=\textwidth]{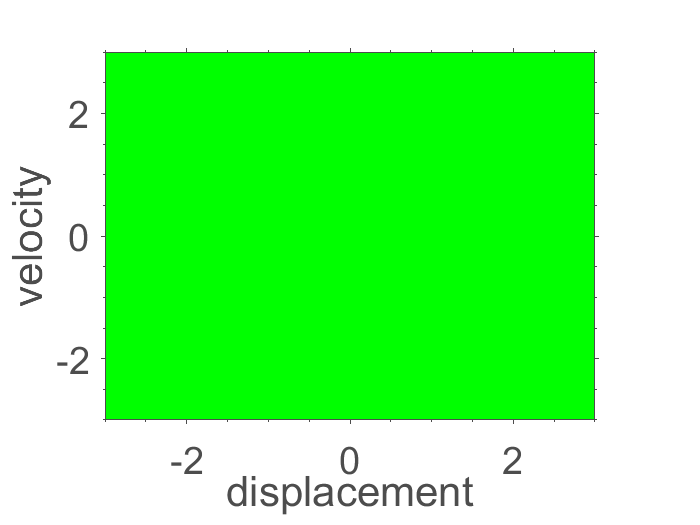}
         \caption{$f = 0.211$}
     \end{subfigure}
     \begin{subfigure}[b]{0.31\textwidth}
         \centering
         \includegraphics[width=\textwidth]{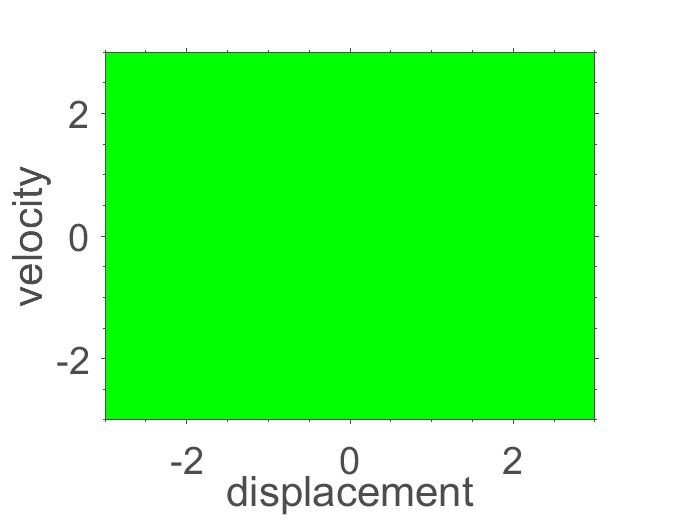}
         \caption{$f = 0.243$}
     \end{subfigure}
     \begin{subfigure}[b]{0.31\textwidth}
         \centering
         \includegraphics[width=\textwidth]{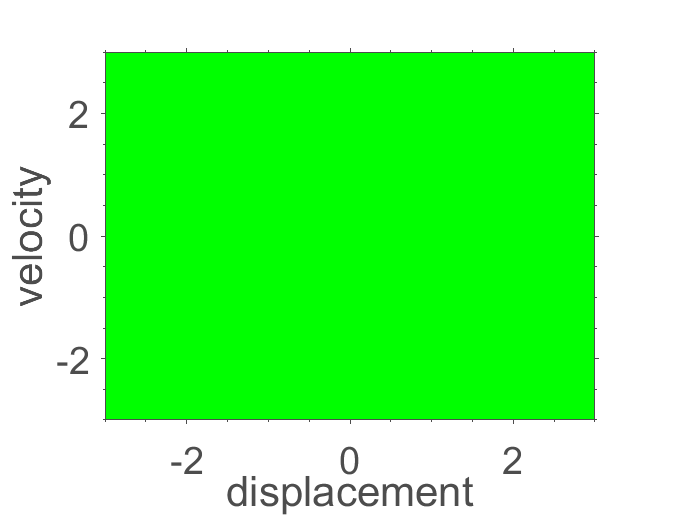}
         \caption{$f = 0.275$}
     \end{subfigure}

    \caption{Basins of attraction intersection with the plane $v=0$ for the symmetric model with a forcing frequency of $\Omega=0.8$ under different forcing amplitude values, within the bounds of $-3 \leq x_0 \leq 3$ and $-3 \leq \dot{x}_0 \leq 3$. The corresponding attractor colors are shown in Fig.\ref{fig:atractors_O08}.}
    \label{fig:basins_O08}
\end{figure}

\begin{figure}
    \centering
    \begin{subfigure}[b]{0.31\textwidth}
         \centering
         \includegraphics[width=\textwidth]{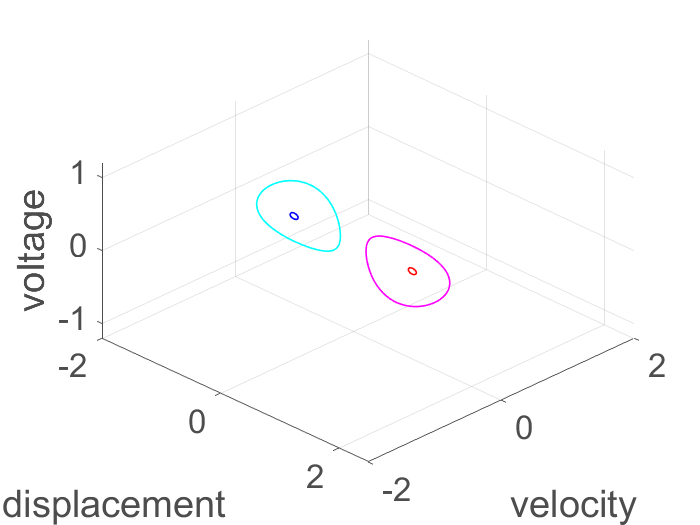}
         \caption{$f = 0.019$}
     \end{subfigure}
     \begin{subfigure}[b]{0.31\textwidth}
         \centering
         \includegraphics[width=\textwidth]{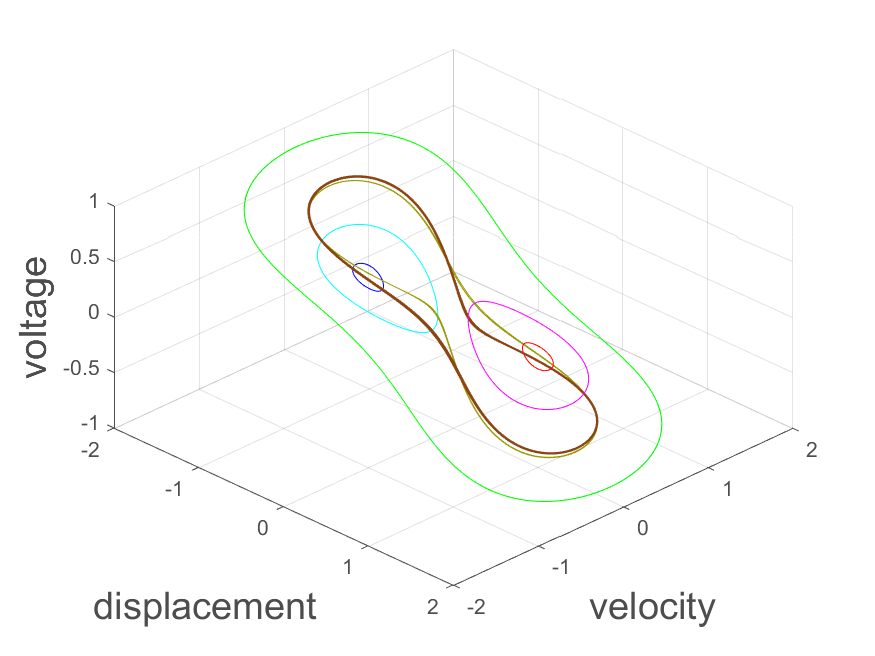}
         \caption{$f = 0.051$}
     \end{subfigure}
     \begin{subfigure}[b]{0.31\textwidth}
         \centering
         \includegraphics[width=\textwidth]{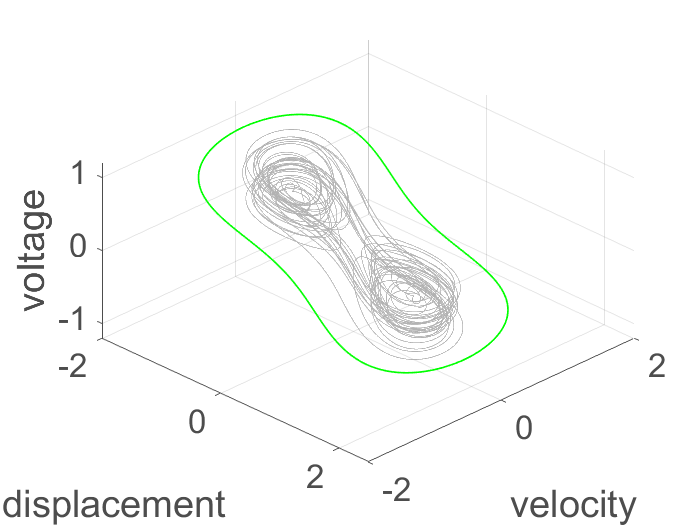}
         \caption{$f = 0.083$}
     \end{subfigure}
     \begin{subfigure}[b]{0.31\textwidth}
         \centering
         \includegraphics[width=\textwidth]{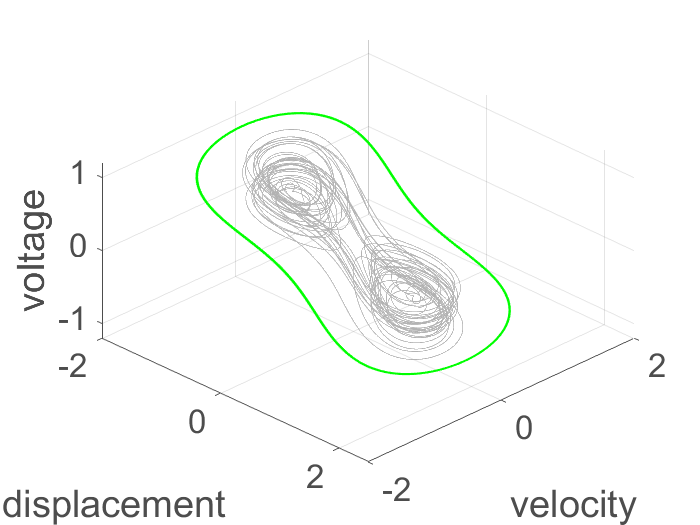}
         \caption{$f = 0.115$}
     \end{subfigure}
     \begin{subfigure}[b]{0.31\textwidth}
         \centering
         \includegraphics[width=\textwidth]{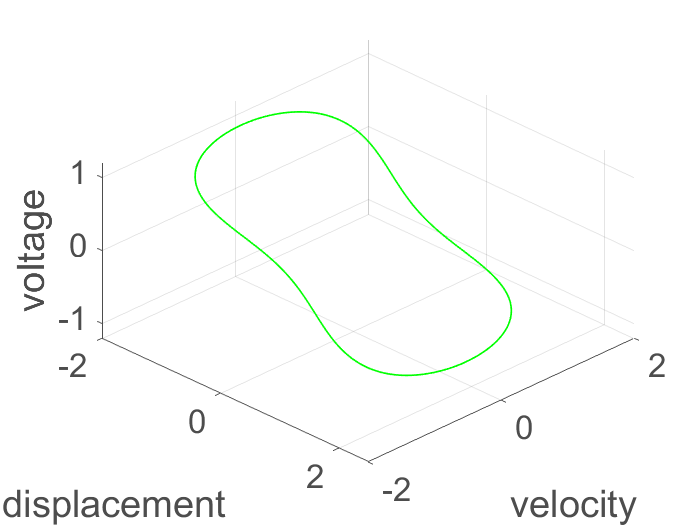}
         \caption{$f = 0.147$}
     \end{subfigure}
     \begin{subfigure}[b]{0.31\textwidth}
         \centering
         \includegraphics[width=\textwidth]{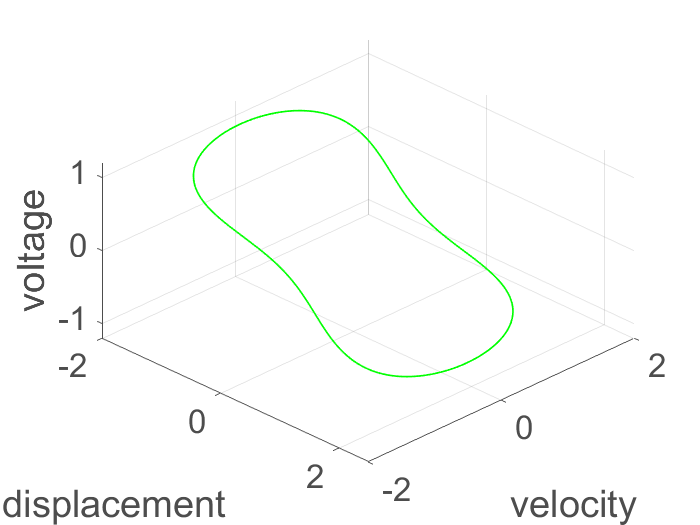}
         \caption{$f = 0.179$}
     \end{subfigure}
     \begin{subfigure}[b]{0.31\textwidth}
         \centering
         \includegraphics[width=\textwidth]{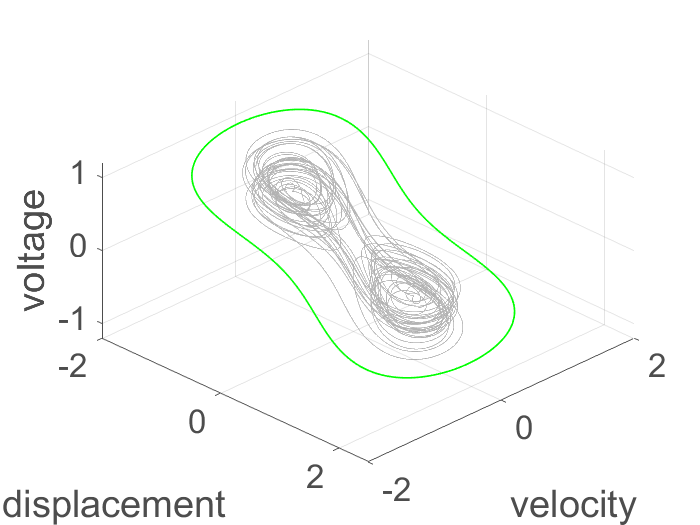}
         \caption{$f = 0.211$}
     \end{subfigure}
     \begin{subfigure}[b]{0.31\textwidth}
         \centering
         \includegraphics[width=\textwidth]{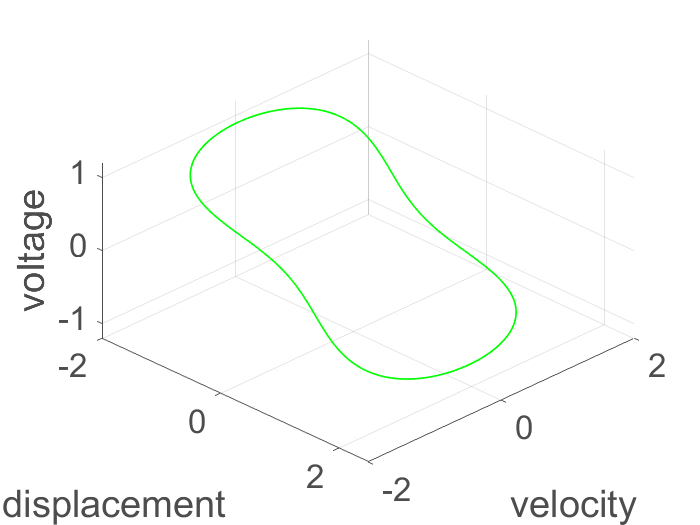}
         \caption{$f = 0.243$}
     \end{subfigure}
     \begin{subfigure}[b]{0.31\textwidth}
         \centering
         \includegraphics[width=\textwidth]{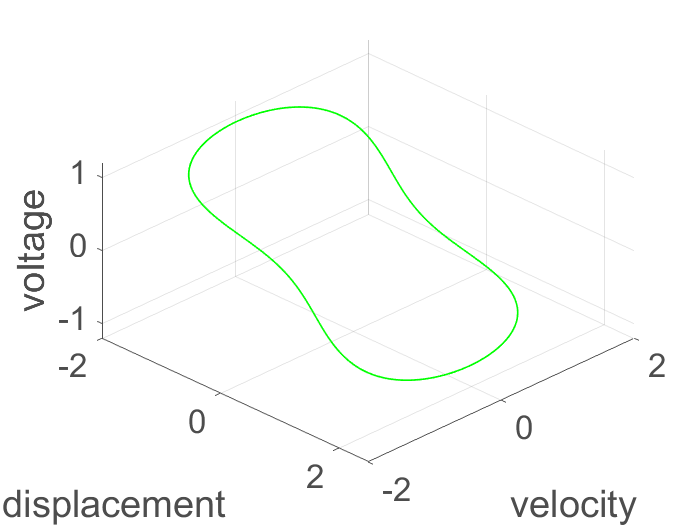}
         \caption{$f = 0.275$}
     \end{subfigure}

	\caption{Attractors of the symmetric model for $\Omega=0.8$ under different forcing amplitude values. Attractors correspond to basins of attraction presented in Fig.\ref{fig:basins_O08}.}
	\label{fig:atractors_O08}
\end{figure}

Now the study analyzes the system's response when the excitation amplitude varies from $0.019$ to $0.275$ with intervals of $0.032$, and the frequency is fixed at $\Omega = 0.8$. This frequency value was chosen for its optimal performance in creating green basins, as shown in Fig.~\ref{fig:basins_O08} and \ref{fig:atractors_O08} that depict the basins of attraction and their corresponding attractors. At $\mathnormal{f} = 0.019$, four regular attractors with fragmented regions illustrate the system's sensitivity to initial conditions. Red and blues basins dominate, but the cyan and magenta basins appear around the equilibrium points. The green basin appears as the excitation amplitude increases to $\mathnormal{f} = 0.051$, but the purple attractor dominates. The purple attractor seems to be a homoclinic orbit. Blue and red areas decrease from this condition. With higher amplitudes, the attractors become more energetic, as seen in $\mathnormal{f} = 0.115$, where the only possible basin is the most energetic or chaotic attractor. This trend continues for higher excitation amplitudes and the attraction basins are covered by green area.

\subsection{Asymmetric bistable energy harvester}

The basins of attraction for asymmetric bistable energy harvesters are investigated. Firstly, we analyze the asymmetric model with a sloping angle of 0, where the only source of asymmetry comes from the quadratic term ($\delta$). Different values of the quadratic term are studied. Then, we introduce a high angle to study the system under strong asymmetry. For these cases, the frequency of excitation is fixed at $\Omega = 0.8$ to compare to the best performance for the symmetric model presented previously, and different forcing amplitude values are tested.

\begin{figure}
    \centering
    \begin{subfigure}[b]{0.31\textwidth}
         \centering
         \includegraphics[width=\textwidth]{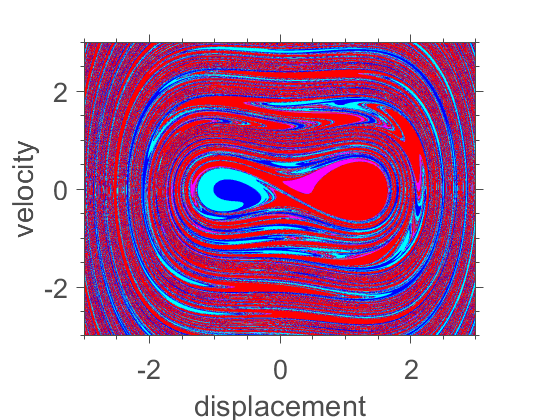}
         \caption{$f = 0.019$ and $\delta = 0.15$}
     \end{subfigure}
     \begin{subfigure}[b]{0.31\textwidth}
         \centering
         \includegraphics[width=\textwidth]{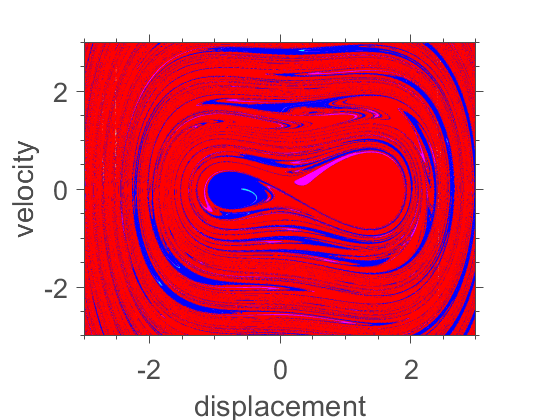}
         \caption{$f = 0.019$ and $\delta = 0.30$}
     \end{subfigure}
     \begin{subfigure}[b]{0.31\textwidth}
         \centering
         \includegraphics[width=\textwidth]{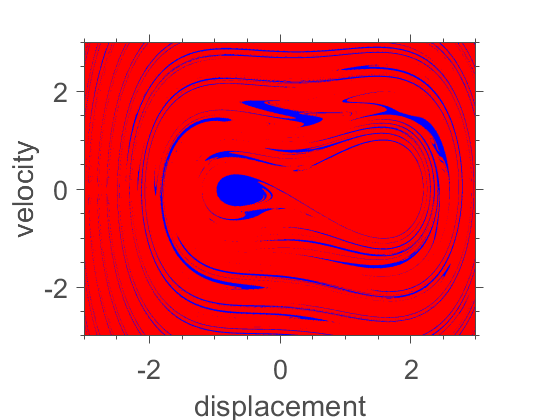}
         \caption{$f = 0.019$ and $\delta = 0.45$}
     \end{subfigure}
     \begin{subfigure}[b]{0.31\textwidth}
         \centering
         \includegraphics[width=\textwidth]{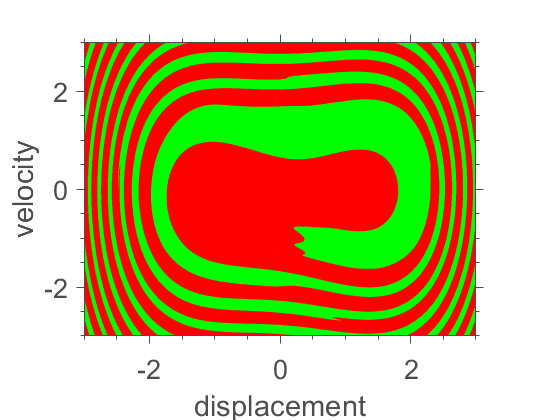}
         \caption{$f = 0.083$ and $\delta = 0.15$}
     \end{subfigure}
     \begin{subfigure}[b]{0.31\textwidth}
         \centering
         \includegraphics[width=\textwidth]{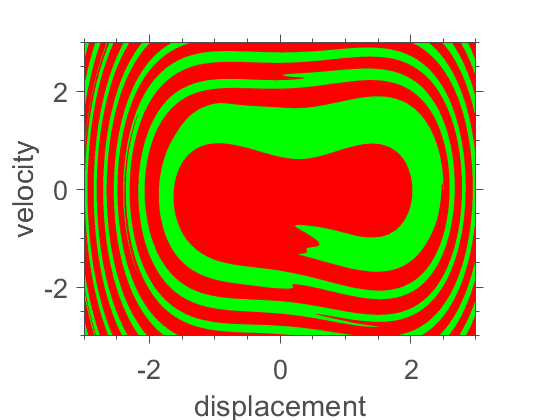}
         \caption{$f = 0.083$ and $\delta = 0.30$}
     \end{subfigure}
     \begin{subfigure}[b]{0.31\textwidth}
         \centering
         \includegraphics[width=\textwidth]{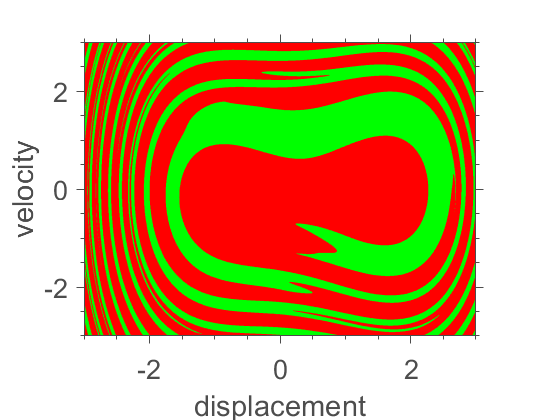}
         \caption{$f = 0.083$ and $\delta = 0.45$}
     \end{subfigure}
     \begin{subfigure}[b]{0.31\textwidth}
         \centering
         \includegraphics[width=\textwidth]{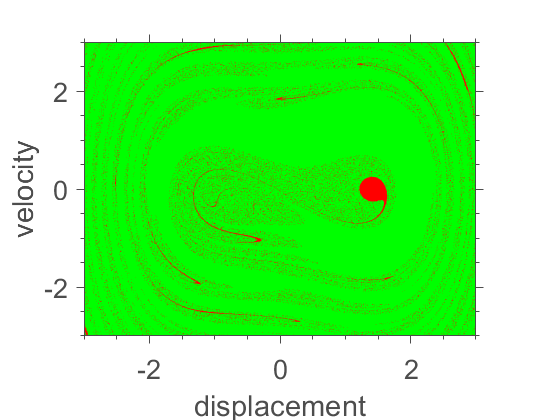}
         \caption{$f = 0.115$ and $\delta = 0.15$}
     \end{subfigure}
     \begin{subfigure}[b]{0.31\textwidth}
         \centering
         \includegraphics[width=\textwidth]{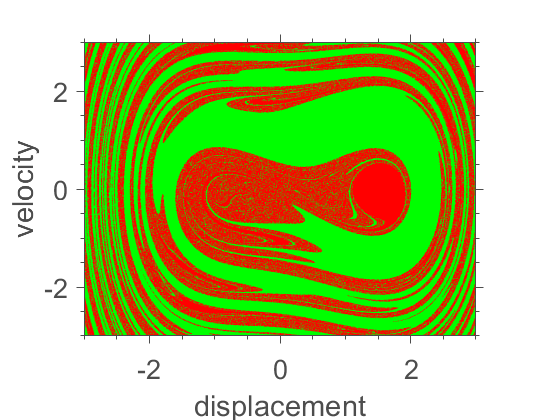}
         \caption{$f = 0.115$ and $\delta = 0.30$}
     \end{subfigure}
     \begin{subfigure}[b]{0.31\textwidth}
         \centering
         \includegraphics[width=\textwidth]{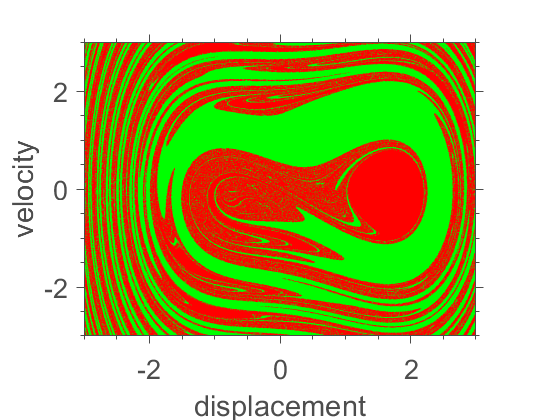}
         \caption{$f = 0.115$ and $\delta = 0.45$}
     \end{subfigure}
     
    \caption{Basins of attraction intersection with the plane $v=0$ for the asymmetric model with a forcing frequency of $\Omega=0.8$ and inclination $\phi=0^\circ$, within the bounds of $-3 \leq x_0 \leq 3$ and $-3 \leq \dot{x}_0 \leq 3$. Forcing amplitude of $0.019, 0.083$, and $0.115$ under quadratic terms of $0.15, 0.30$, and $0.45$ are investigated. The corresponding attractor colors are shown in Fig.\ref{fig:atractor_phi0}. }
    \label{fig:basins_phi0}
\end{figure}
\begin{figure}
    \centering
    \begin{subfigure}[b]{0.31\textwidth}
         \centering
         \includegraphics[width=\textwidth]{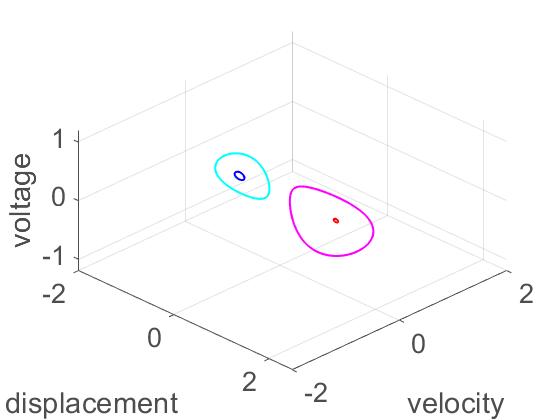}
         \caption{$f = 0.019$ and $\delta = 0.15$}
     \end{subfigure}
     \begin{subfigure}[b]{0.31\textwidth}
         \centering
         \includegraphics[width=\textwidth]{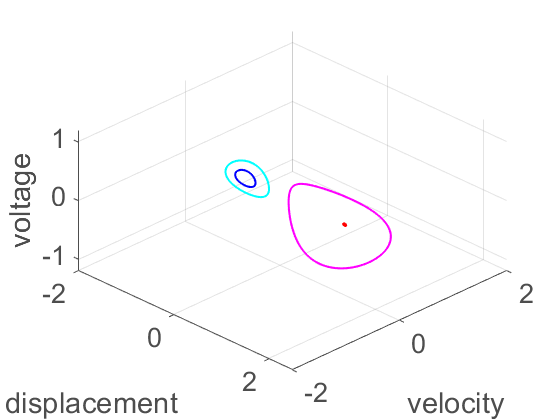}
         \caption{$f = 0.019$ and $\delta = 0.30$}
     \end{subfigure}
     \begin{subfigure}[b]{0.31\textwidth}
         \centering
         \includegraphics[width=\textwidth]{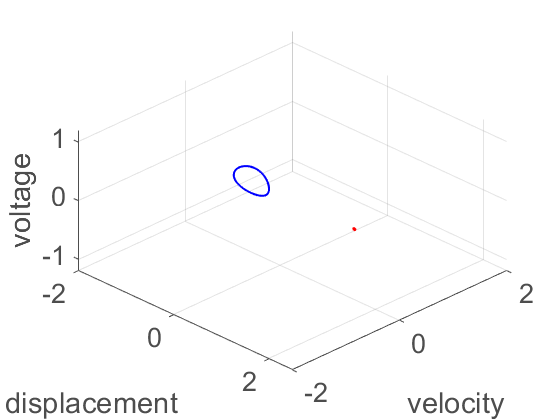}
         \caption{$f = 0.019$ and $\delta = 0.45$}
     \end{subfigure}
     \begin{subfigure}[b]{0.31\textwidth}
         \centering
         \includegraphics[width=\textwidth]{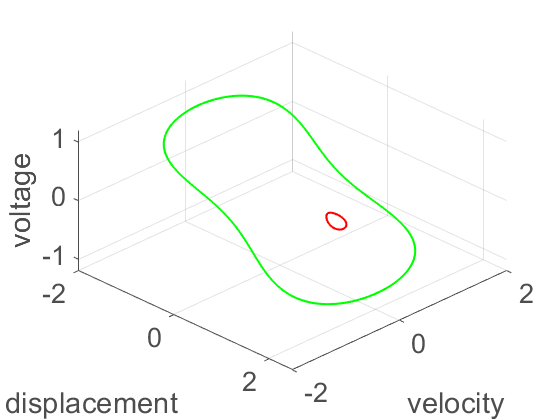}
         \caption{$f = 0.083$ and $\delta = 0.15$}
     \end{subfigure}
     \begin{subfigure}[b]{0.31\textwidth}
         \centering
         \includegraphics[width=\textwidth]{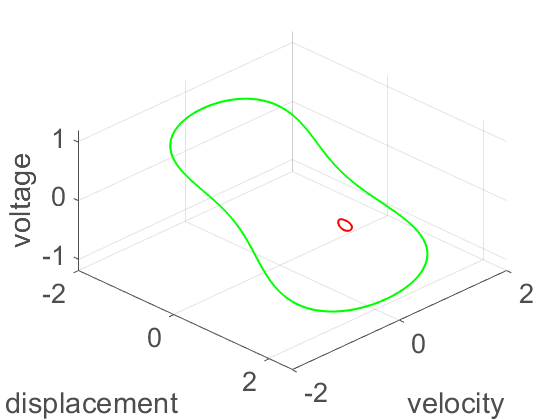}
         \caption{$f = 0.083$ and $\delta = 0.30$}
     \end{subfigure}
     \begin{subfigure}[b]{0.31\textwidth}
         \centering
         \includegraphics[width=\textwidth]{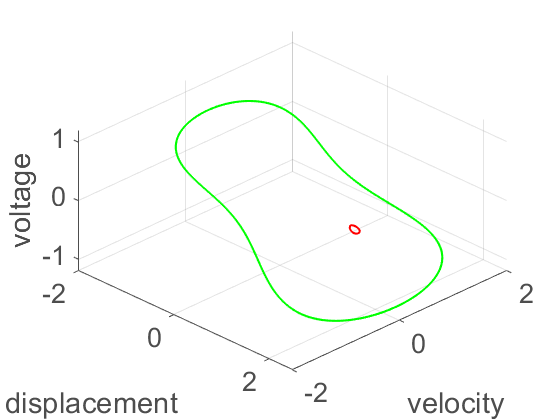}
         \caption{$f = 0.083$ and $\delta = 0.45$}
     \end{subfigure}
     \begin{subfigure}[b]{0.31\textwidth}
         \centering
         \includegraphics[width=\textwidth]{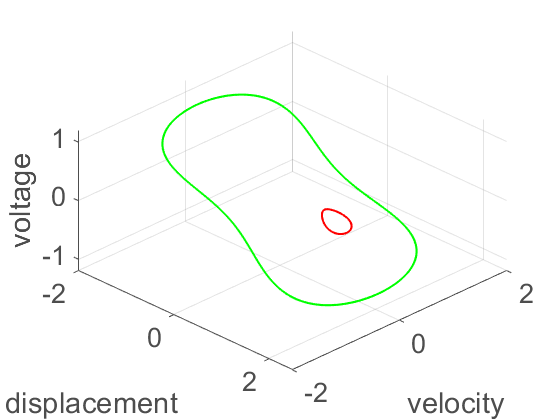}
         \caption{$f = 0.115$ and $\delta = 0.15$}
     \end{subfigure}
     \begin{subfigure}[b]{0.31\textwidth}
         \centering
         \includegraphics[width=\textwidth]{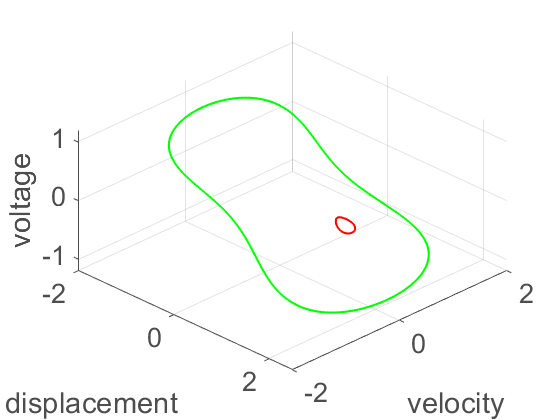}
         \caption{$f = 0.115$ and $\delta = 0.30$}
     \end{subfigure}
     \begin{subfigure}[b]{0.31\textwidth}
         \centering
         \includegraphics[width=\textwidth]{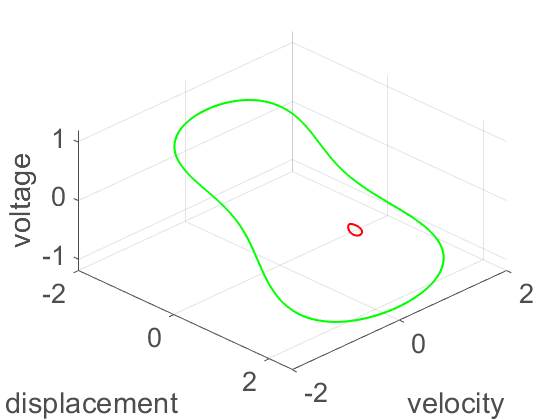}
         \caption{$f = 0.115$ and $\delta = 0.45$}
     \end{subfigure}

	\caption{Attractors of the asymmetric model for $\Omega=0.8$ and $\phi=0^\circ$ under different forcing amplitudes and quadratic terms. Attractors correspond to the basins of attraction presented in Fig.\ref{fig:basins_phi0}.}
	\label{fig:atractor_phi0}
\end{figure}

Fig.~\ref{fig:basins_phi0} and \ref{fig:atractor_phi0} show the basis of attraction and attractors for the asymmetric bistable model at a sloping angle of zero for various quadratic terms. Three quadratic terms are analyzed, $\delta = 0.15, 0.30$ and $0.45$, under different amplitudes of excitation, $f=0.019, 0.083$ and $0.115$. When $f=0.019$ and $\delta=0.15$ (Fig.~\ref{fig:basins_phi0}a), the red basin is dominant, but there are high-density regions with blue, cyan, and magenta basins. As the $\delta$ increases (Fig.~\ref{fig:basins_phi0}b and c), the red area also increases, almost covering the total area. The corresponding attractors in Fig.~\ref{fig:atractor_phi0}a,b, and c indicate that the energy of the red orbit decreases as $\delta$ increases. Fig.~\ref{fig:basins_phi0}d,e, and f show the basins of attraction when $f=0.083$ for $\delta=0.15,0.30$ and $0.45$, respectively. Red and green basins cover them, and their areas almost did not change as $\delta$ increased. The notable difference is verified at the red energy orbit in Fig.~\ref{fig:atractor_phi0}d,e, and f, which decreases as $\delta$ increases. Finally, Fig.~\ref{fig:basins_phi0}g, h, and i show basins of attraction for $f=0.115$ and $\delta=0.15, 0.30$ and $0.45$, respectively. The negative effect of the quadratic term is evident in this case. As $\delta$ increases, the red area dominates the region, and the green area decreases. The reduction of the red energy orbit can also be visualized in Fig.~\ref{fig:atractor_phi0}g,h, and i.

These results demonstrate that asymmetry from the quadratic term can harm the harvesting process, reducing the conditions for high-energy orbits and decreasing the amplitude of oscillations. Even when $\delta=0.15$, where the basins are similar to those of the symmetric model, the quadratic term still affects the energy of orbits, making its influence more visible when examining the attractors.

\begin{figure}
    \centering
    \begin{subfigure}[b]{0.31\textwidth}
         \centering
         \includegraphics[width=\textwidth]{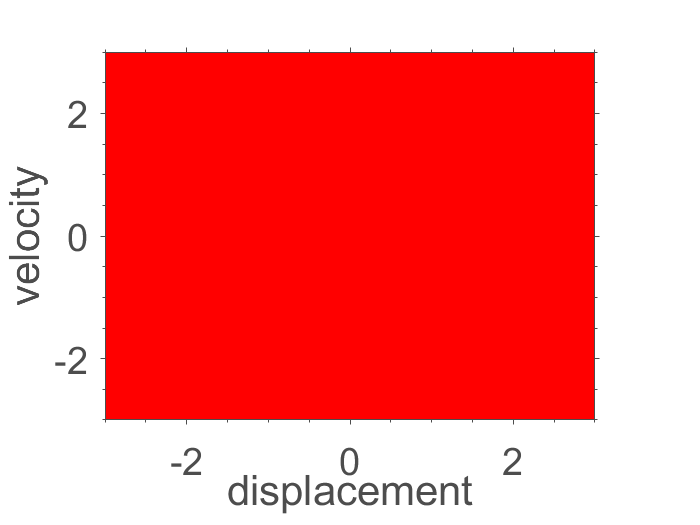}
         \caption{$f = 0.019$}
     \end{subfigure}
     \begin{subfigure}[b]{0.31\textwidth}
         \centering
         \includegraphics[width=\textwidth]{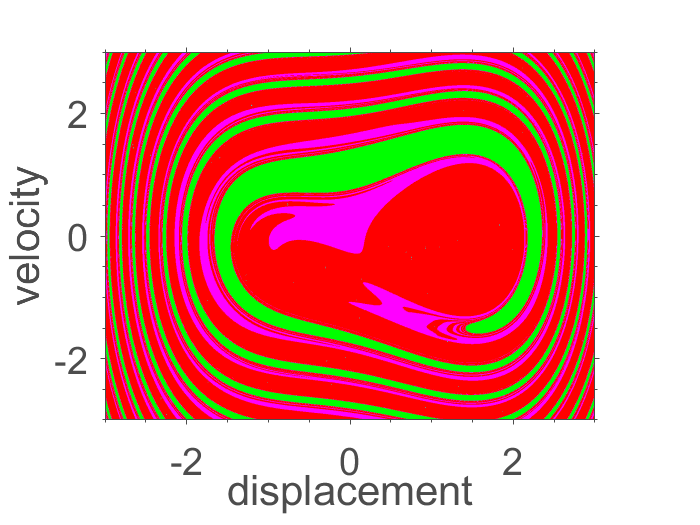}
         \caption{$f = 0.051$}
     \end{subfigure}
     \begin{subfigure}[b]{0.31\textwidth}
         \centering
         \includegraphics[width=\textwidth]{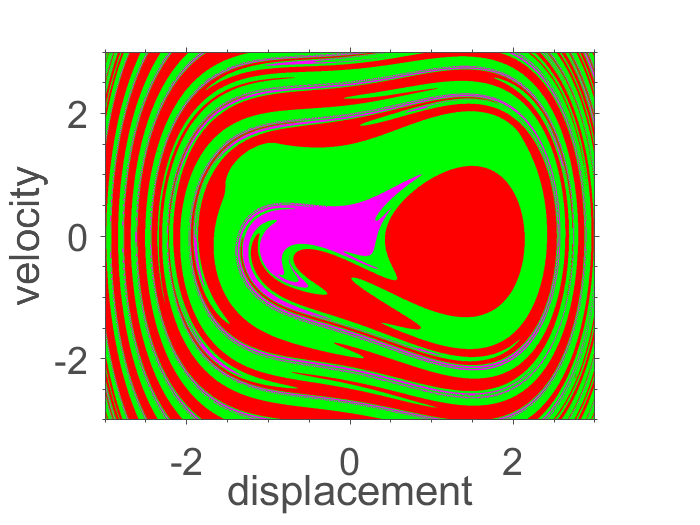}
         \caption{$f = 0.083$}
     \end{subfigure}
     \begin{subfigure}[b]{0.31\textwidth}
         \centering
         \includegraphics[width=\textwidth]{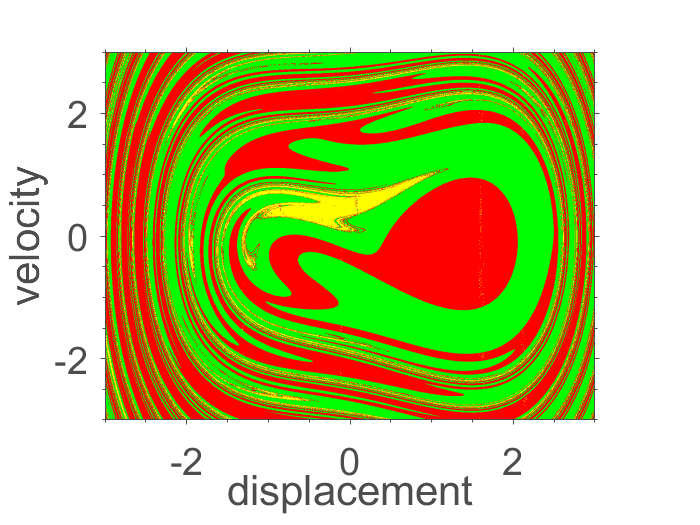}
         \caption{$f = 0.115$}
     \end{subfigure}
     \begin{subfigure}[b]{0.31\textwidth}
         \centering
         \includegraphics[width=\textwidth]{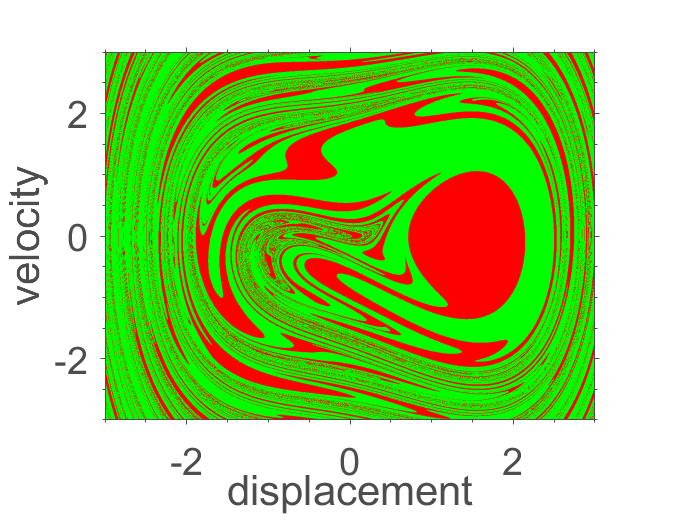}
         \caption{$f = 0.147$}
     \end{subfigure}
     \begin{subfigure}[b]{0.31\textwidth}
         \centering
         \includegraphics[width=\textwidth]{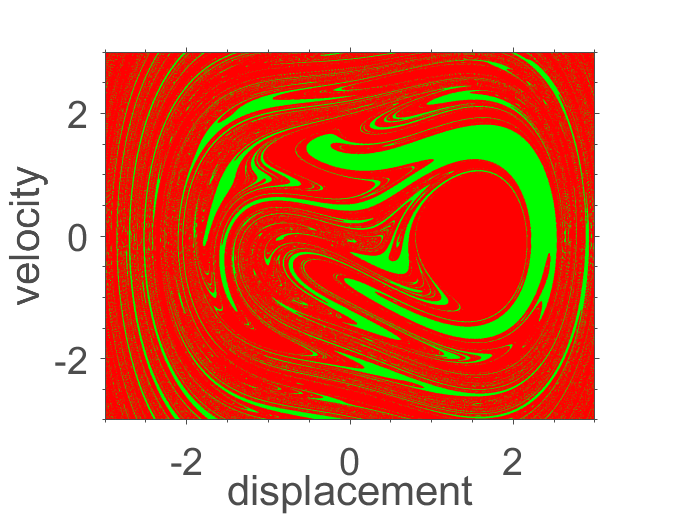}
         \caption{$f = 0.179$}
     \end{subfigure}
     \begin{subfigure}[b]{0.31\textwidth}
         \centering
         \includegraphics[width=\textwidth]{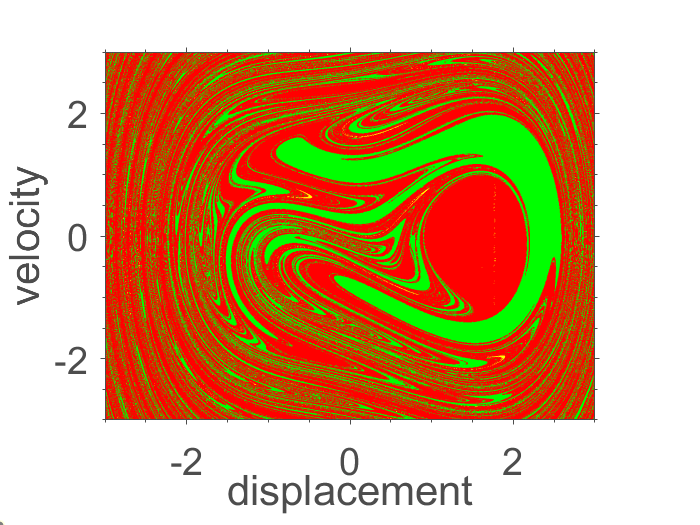}
         \caption{$f = 0.211$}
     \end{subfigure}
     \begin{subfigure}[b]{0.31\textwidth}
         \centering
         \includegraphics[width=\textwidth]{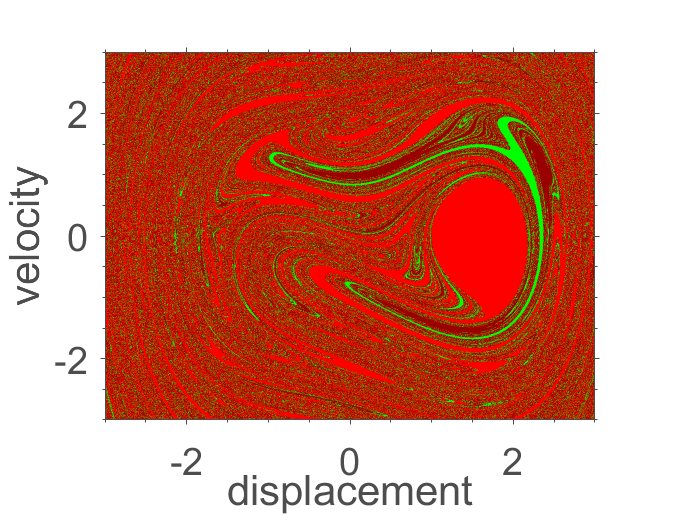}
         \caption{$f = 0.243$}
     \end{subfigure}
     \begin{subfigure}[b]{0.31\textwidth}
         \centering
         \includegraphics[width=\textwidth]{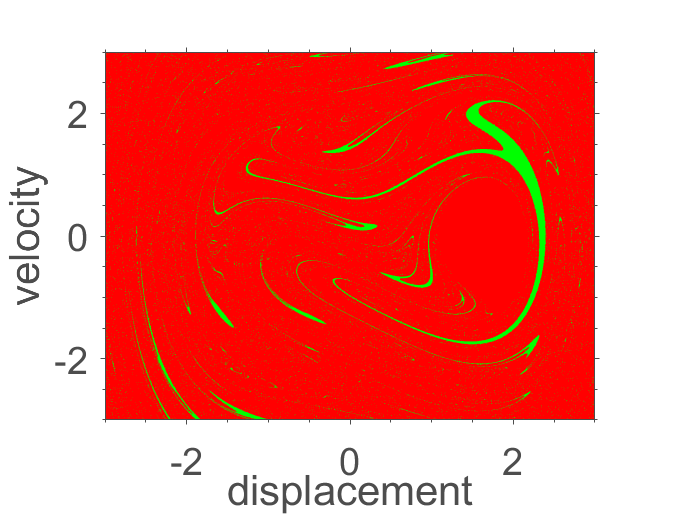}
         \caption{$f = 0.275$}
     \end{subfigure}

    \caption{Basins of attraction intersection with the plane $v=0$ for the asymmetric model with a forcing frequency of $\Omega=0.8$ and inclination $\phi=35^\circ$, within the bounds of $-3 \leq x_0 \leq 3$ and $-3 \leq \dot{x}_0 \leq 3$. Different forcing amplitudes are shown. The corresponding attractor colors are shown in Fig.\ref{fig:atractor_phi35}.}
    \label{fig:basins_phi35}
\end{figure}

\begin{figure}
    \centering
    \begin{subfigure}[b]{0.31\textwidth}
         \centering
         \includegraphics[width=\textwidth]{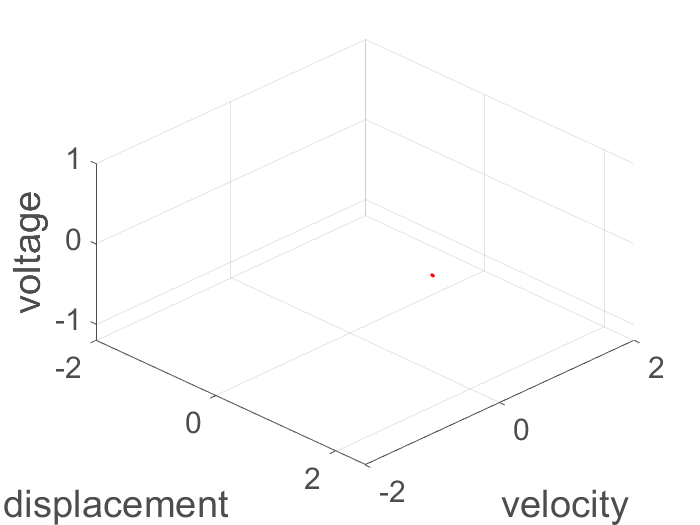}
         \caption{$f = 0.019$}
     \end{subfigure}
     \begin{subfigure}[b]{0.31\textwidth}
         \centering
         \includegraphics[width=\textwidth]{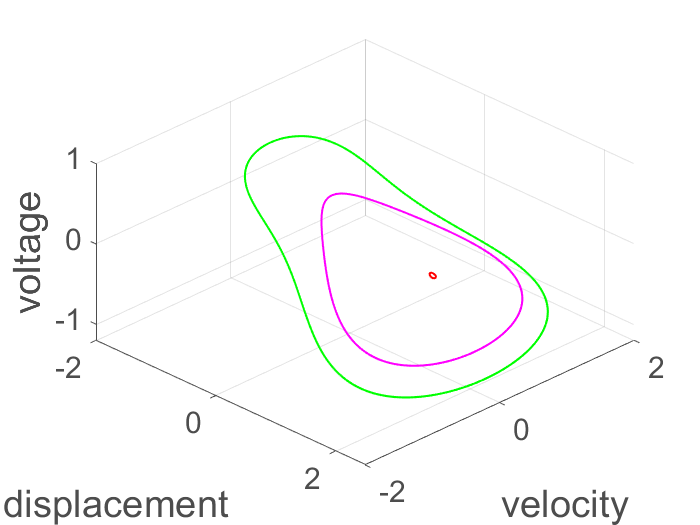}
         \caption{$f = 0.051$}
     \end{subfigure}
     \begin{subfigure}[b]{0.31\textwidth}
         \centering
         \includegraphics[width=\textwidth]{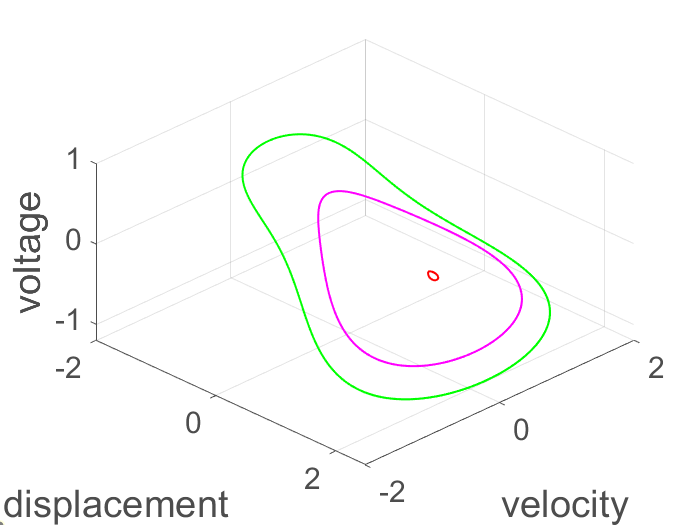}
         \caption{$f = 0.083$}
     \end{subfigure}
     \begin{subfigure}[b]{0.31\textwidth}
         \centering
         \includegraphics[width=\textwidth]{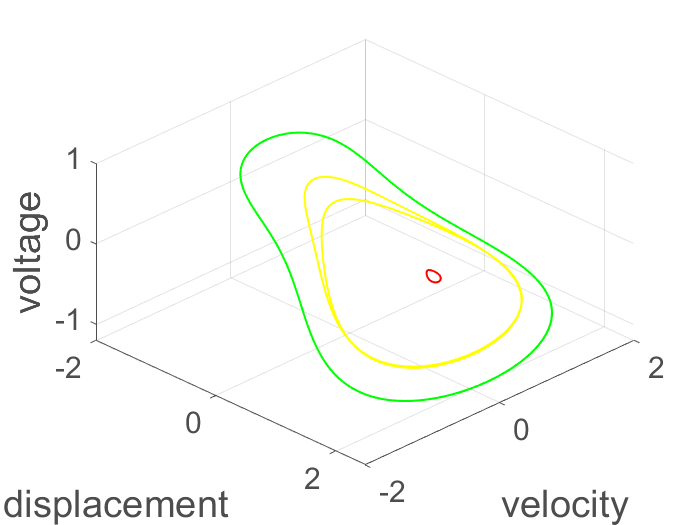}
         \caption{$f = 0.115$}
     \end{subfigure}
     \begin{subfigure}[b]{0.31\textwidth}
         \centering
         \includegraphics[width=\textwidth]{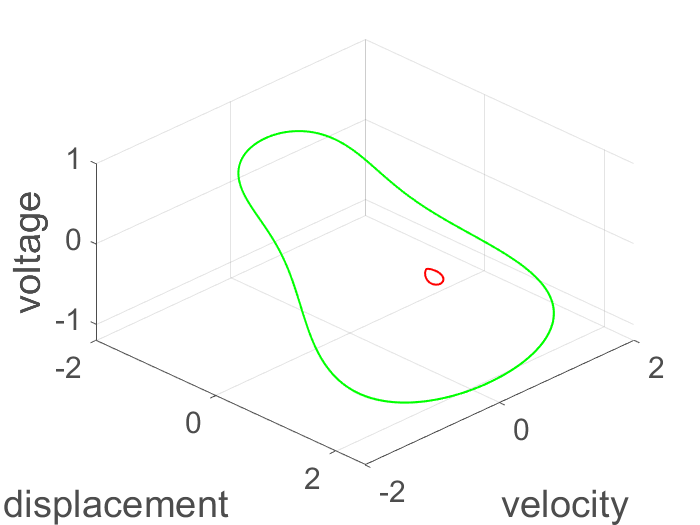}
         \caption{$f = 0.147$}
     \end{subfigure}
     \begin{subfigure}[b]{0.31\textwidth}
         \centering
         \includegraphics[width=\textwidth]{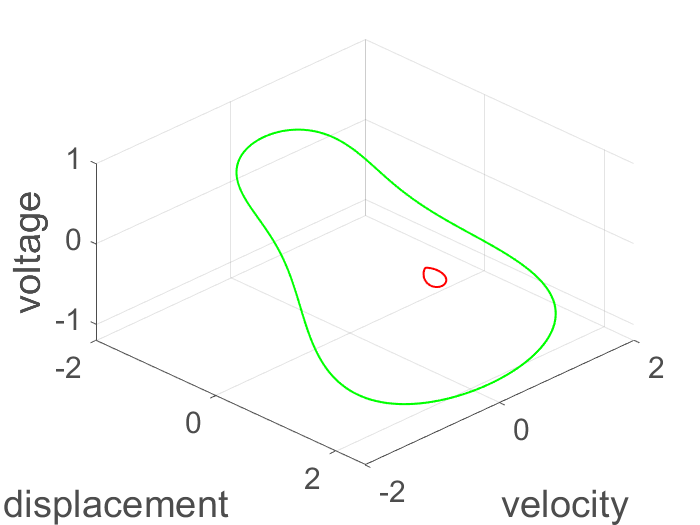}
         \caption{$f = 0.179$}
     \end{subfigure}
     \begin{subfigure}[b]{0.31\textwidth}
         \centering
         \includegraphics[width=\textwidth]{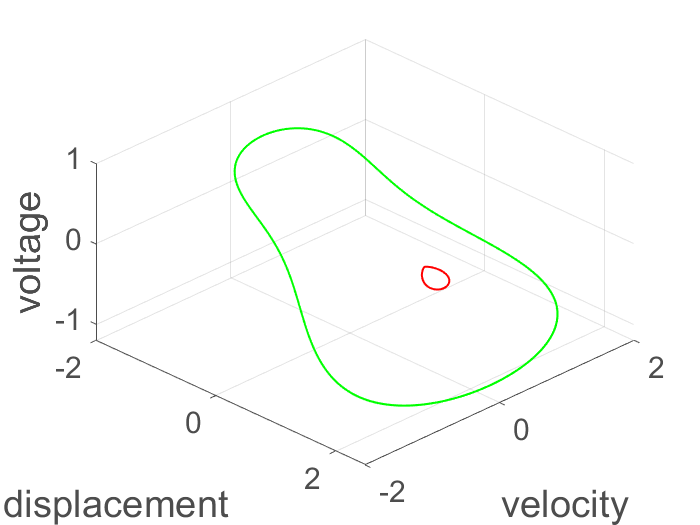}
         \caption{$f = 0.211$}
     \end{subfigure}
     \begin{subfigure}[b]{0.31\textwidth}
         \centering
         \includegraphics[width=\textwidth]{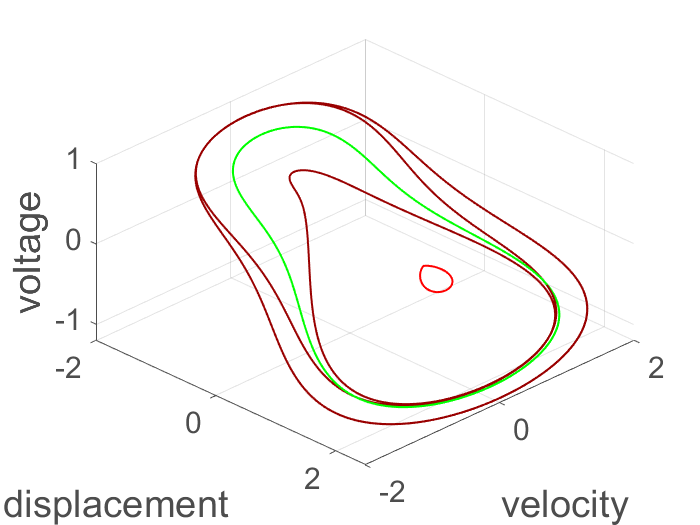}
         \caption{$f = 0.243$}
     \end{subfigure}
     \begin{subfigure}[b]{0.31\textwidth}
         \centering
         \includegraphics[width=\textwidth]{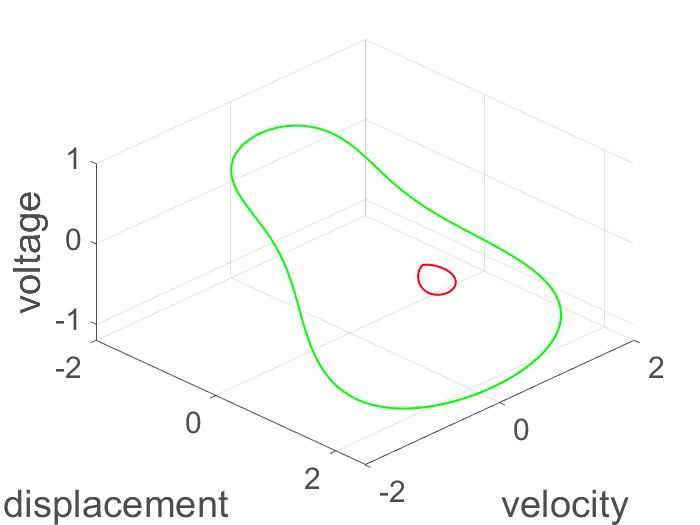}
         \caption{$f = 0.275$}
     \end{subfigure}
	\caption{Attractors of the asymmetric model for $\Omega=0.8$ and $\phi=35^\circ$ under different forcing amplitudes. Attractors correspond to basins of attraction presented in Fig.\ref{fig:basins_phi35}.}
	\label{fig:atractor_phi35}
\end{figure}

Basins of attraction and attractors for an asymmetric bistable model when $\phi = 35^\circ$, $\delta=0.15$ under different forcing amplitude values are shown in Fig.~\ref{fig:basins_phi35} and \ref{fig:atractor_phi35}. A high value of the sloping angle is chosen to investigate a strongly asymmetric condition. The excitation amplitude varies from $0.019\leq\mathnormal{f}\leq0.275$ with increments of $0.032$. The initial displacement and velocity plane is defined as $(x_0,\dot{x}_0) \in [-3 \leq x_0 \leq 3] \times [-3 \leq \dot{x}_0 \leq 3]$ with the restriction $v_0 = 0$.

For the lowest amplitude of excitation ($\mathnormal{f} = 0.019$), the plane is dominated by the red basins, representing low-energy attractors so small that they look like a point. This is the worst condition for energy harvesting. Three solutions are present when $\mathnormal{f} \in \{0.051, 0.083\}$: green, red, and magenta basins. The red attractor has increased energy, and its orbit is now visible. The green attractor is asymmetrical, with the positive equilibrium point having a deeper well than the negative one. The magenta basins appear near the negative equilibrium point, and its attractor is a homoclinic orbit. The magenta basins disappear when $\mathnormal{f} = 0.115$, and part of this solution becomes a yellow basin with a two-period solution. For $\mathnormal{f} = 0.147$, there are two solutions, and the green basin is dominant, but the other solution has low energy. As the excitation amplitude increases, the red basins become more frequent, and the green basins almost disappear. A purple basin with several frequencies appears when $\mathnormal{f} = 0.243$. Chaos is not observed in any of the conditions studied.

The blue basins, which represent low-energy orbits around the negative equilibrium point, do not appear. This suggests that the solution is attracted to the other equilibrium point, which is physically explained by the deep energy well of that equilibrium point. When the system exceeds the potential barrier, it performs a snap-through motion, where the oscillation goes from one equilibrium point to another. On the other hand, when the initial conditions are near the negative equilibrium point, even with low external energy, the system goes to the positive one, and the solution is attracted to a red orbit. In the basins of attraction for the symmetric model (Fig.~\ref{fig:basins_O08}), blue basins often appear near the negative equilibrium point. However, the asymmetric model fills this region with red basins (or basins with a snap-through motion).

The asymmetric model shows unusual behavior for $\phi = 35^\circ$. Low orbit response dominates as more external energy is inserted into the system, and the best energy harvesting performance is not achieved at the highest excitation amplitude. The best performance occurs at $\mathnormal{f}=0.147$, where the green basin dominates. However, the symmetric model performs better than the asymmetric model in all scenarios, making asymmetry an undesired configuration for energy harvesting.

\subsection{Optimal asymmetric bistable energy harvester}

The optimal angle, defined in Eq.~\ref{eq:optimal_angle}, is studied to illustrate how asymmetry from the quadratic term can be compensated for a sloping angle. Secondly, different angle values of sloping are studied to verify the robustness of this compensation technique. Finally, the relative area for each basin is calculated under different forcing amplitude values for symmetric and asymmetric cases and compared between the optimal asymmetric model.

\begin{figure}
    \centering
    \begin{subfigure}[b]{0.31\textwidth}
         \centering
         \includegraphics[width=\textwidth]{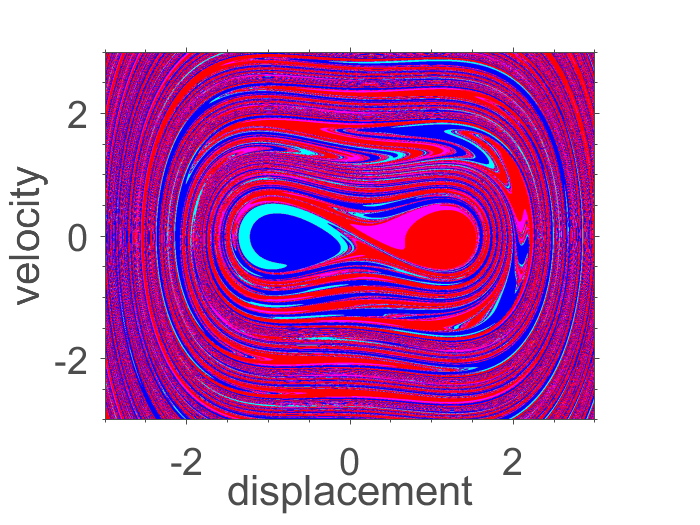}
         \caption{$f = 0.019$}
     \end{subfigure}
     \begin{subfigure}[b]{0.31\textwidth}
         \centering
         \includegraphics[width=\textwidth]{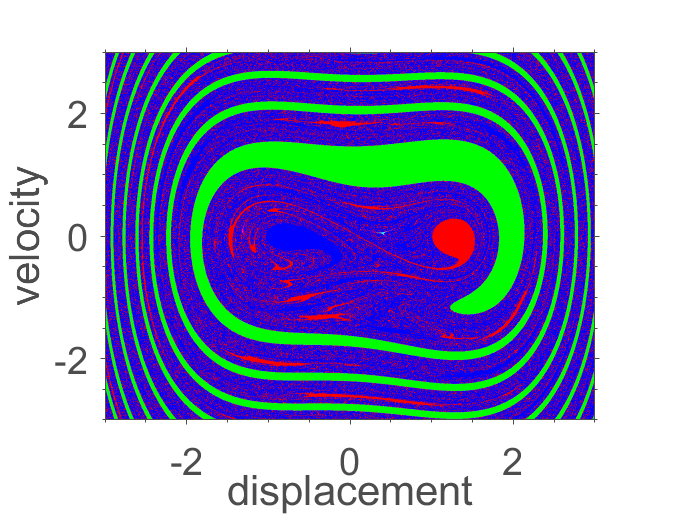}
         \caption{$f = 0.051$}
     \end{subfigure}
     \begin{subfigure}[b]{0.31\textwidth}
         \centering
         \includegraphics[width=\textwidth]{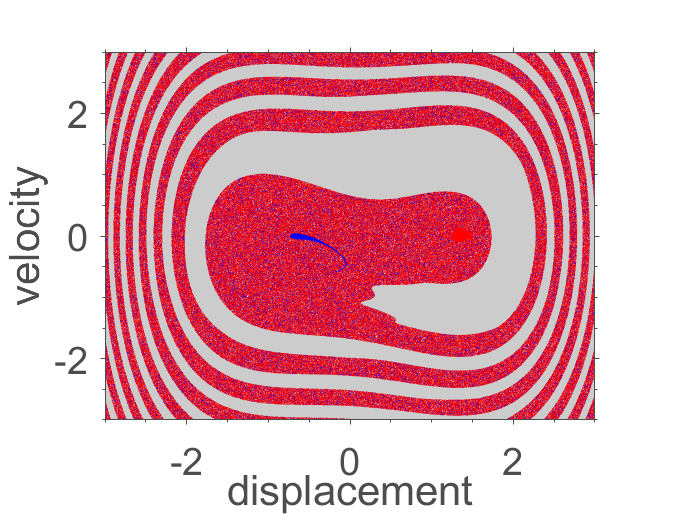}
         \caption{$f = 0.083$}
     \end{subfigure}
     \begin{subfigure}[b]{0.31\textwidth}
         \centering
         \includegraphics[width=\textwidth]{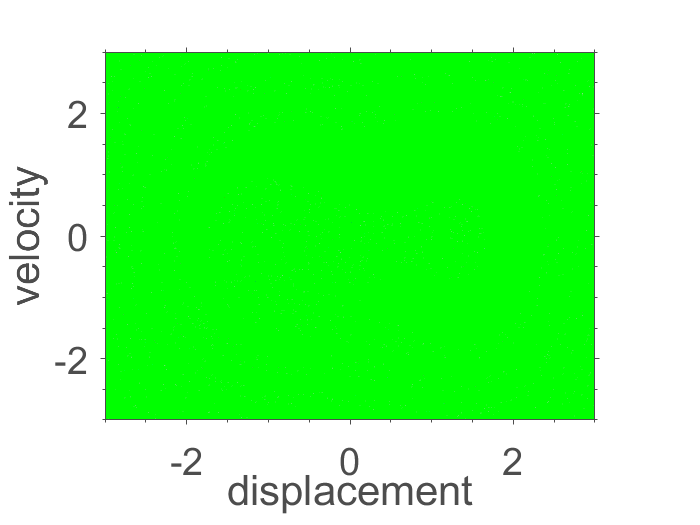}
         \caption{$f = 0.115$}
     \end{subfigure}
     \begin{subfigure}[b]{0.31\textwidth}
         \centering
         \includegraphics[width=\textwidth]{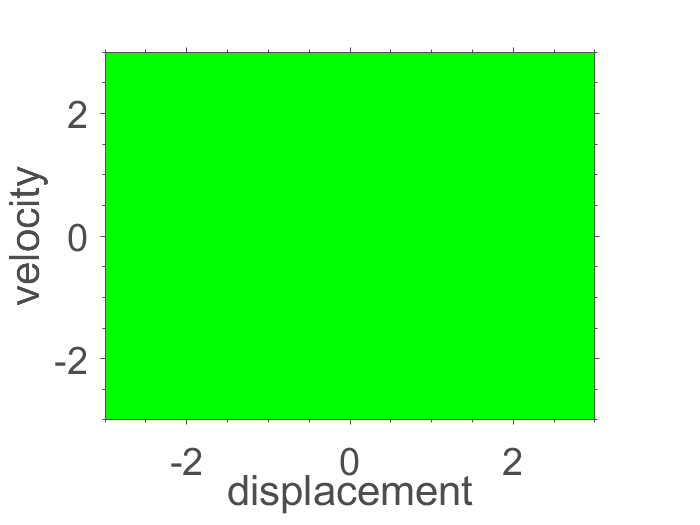}
         \caption{$f = 0.147$}
     \end{subfigure}
     \begin{subfigure}[b]{0.31\textwidth}
         \centering
         \includegraphics[width=\textwidth]{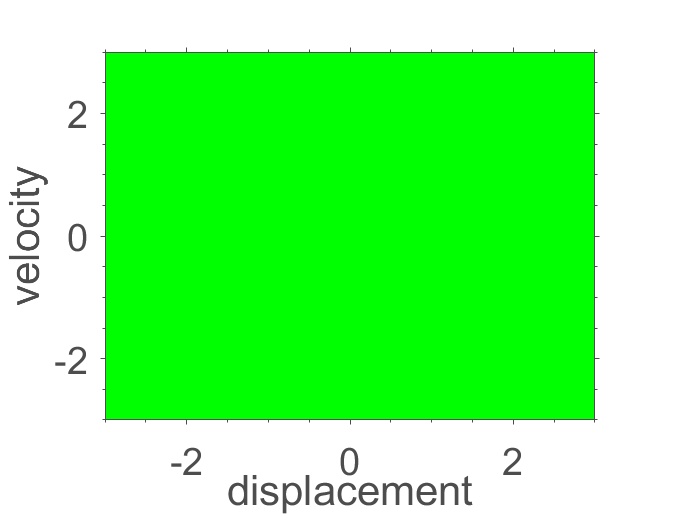}
         \caption{$f = 0.179$}
     \end{subfigure}
     \begin{subfigure}[b]{0.31\textwidth}
         \centering
         \includegraphics[width=\textwidth]{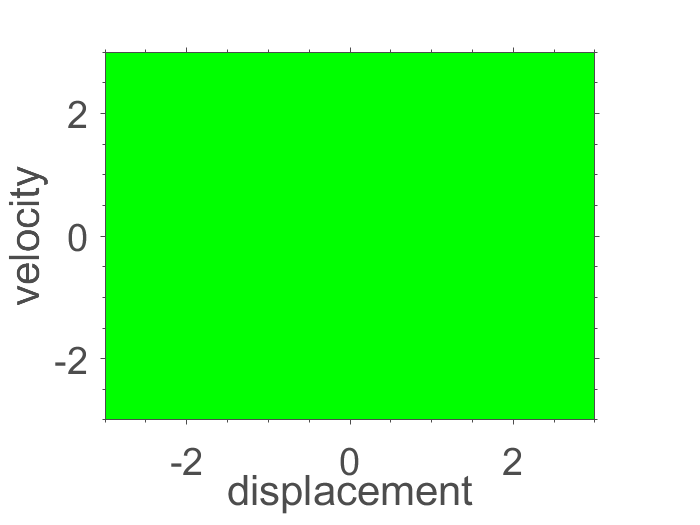}
         \caption{$f = 0.211$}
     \end{subfigure}
     \begin{subfigure}[b]{0.31\textwidth}
         \centering
         \includegraphics[width=\textwidth]{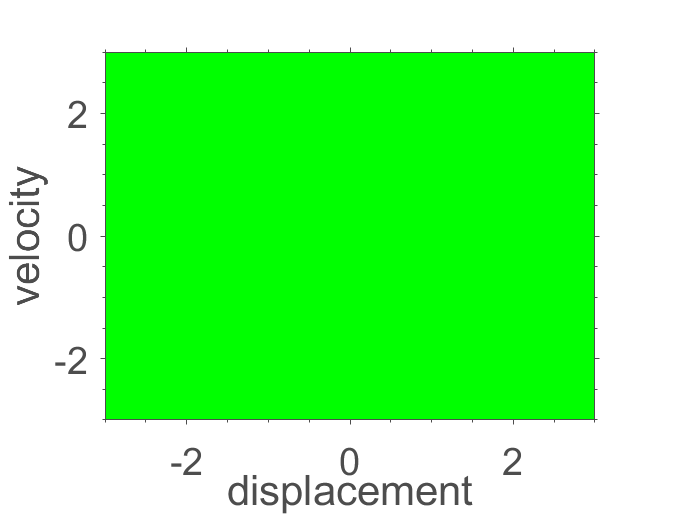}
         \caption{$f = 0.243$}
     \end{subfigure}
     \begin{subfigure}[b]{0.31\textwidth}
         \centering
         \includegraphics[width=\textwidth]{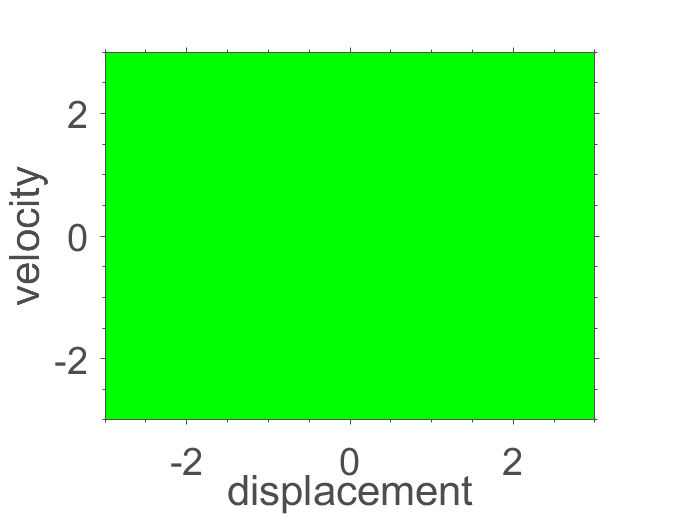}
         \caption{$f = 0.275$}
     \end{subfigure}
    \caption{Basins of attraction intersection with the plane $v=0$ for the asymmetric model with a forcing frequency of $\Omega=0.8$ and inclination $\phi=-4.95^\circ$, within the bounds of $-3 \leq x_0 \leq 3$ and $-3 \leq \dot{x}_0 \leq 3$. Different forcing amplitudes are investigated. The corresponding attractor colors are shown in Fig.\ref{fig:atractor_phi-4}.}
    \label{fig:basins_phi-4}
\end{figure}

\begin{figure}
    \centering
    \begin{subfigure}[b]{0.31\textwidth}
         \centering
         \includegraphics[width=\textwidth]{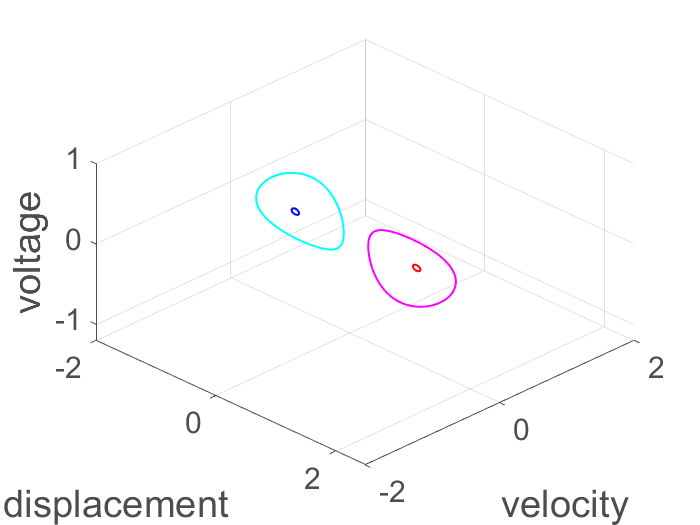}
         \caption{$f = 0.019$}
     \end{subfigure}
     \begin{subfigure}[b]{0.31\textwidth}
         \centering
         \includegraphics[width=\textwidth]{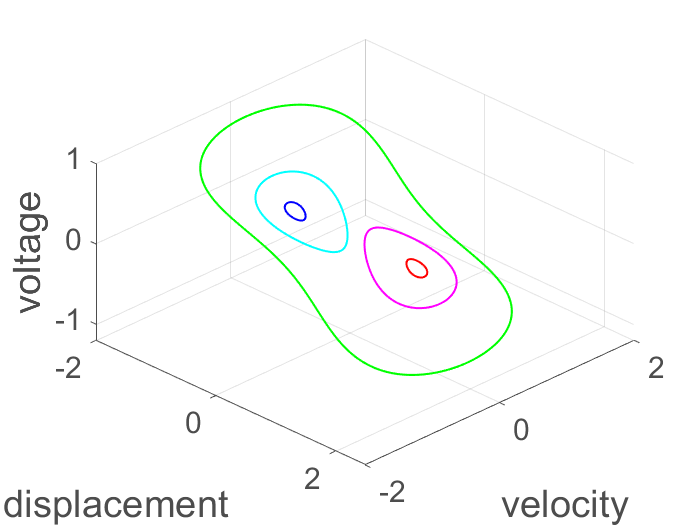}
         \caption{$f = 0.051$}
     \end{subfigure}
     \begin{subfigure}[b]{0.31\textwidth}
         \centering
         \includegraphics[width=\textwidth]{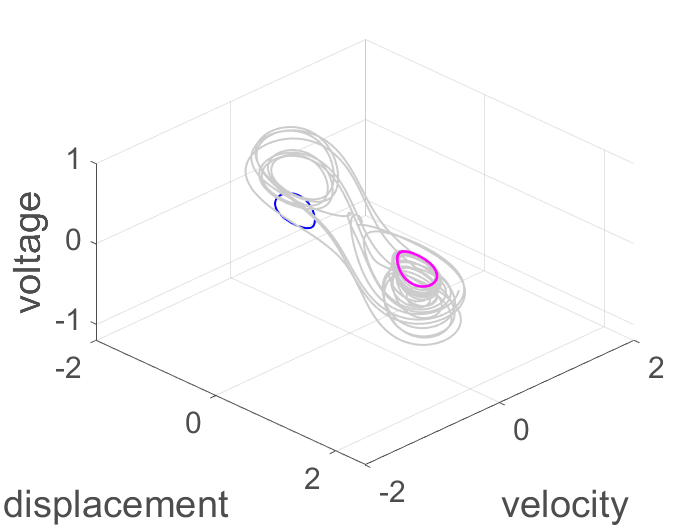}
         \caption{$f = 0.083$}
     \end{subfigure}
     \begin{subfigure}[b]{0.31\textwidth}
         \centering
         \includegraphics[width=\textwidth]{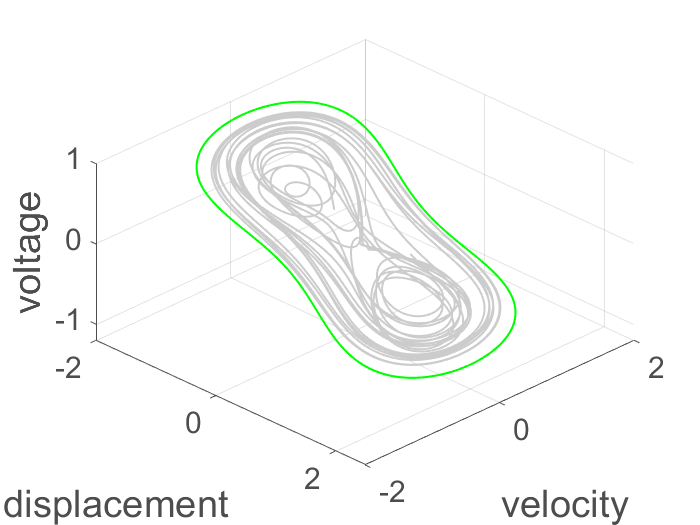}
         \caption{$f = 0.115$}
     \end{subfigure}
     \begin{subfigure}[b]{0.31\textwidth}
         \centering
         \includegraphics[width=\textwidth]{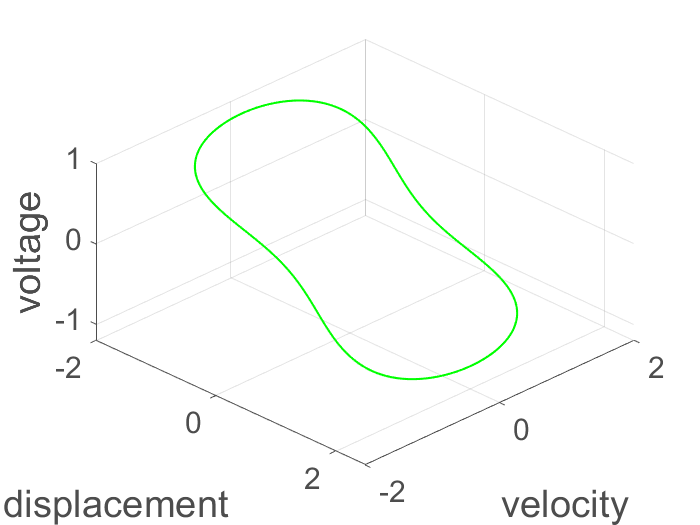}
         \caption{$f = 0.147$}
     \end{subfigure}
     \begin{subfigure}[b]{0.31\textwidth}
         \centering
         \includegraphics[width=\textwidth]{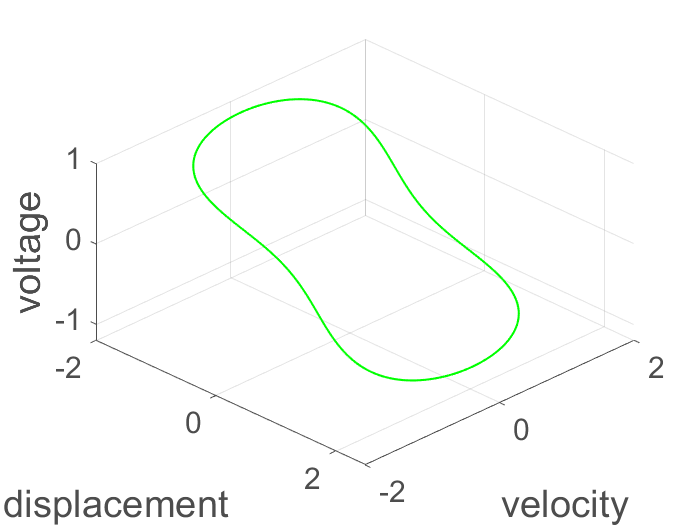}
         \caption{$f = 0.179$}
     \end{subfigure}
     \begin{subfigure}[b]{0.31\textwidth}
         \centering
         \includegraphics[width=\textwidth]{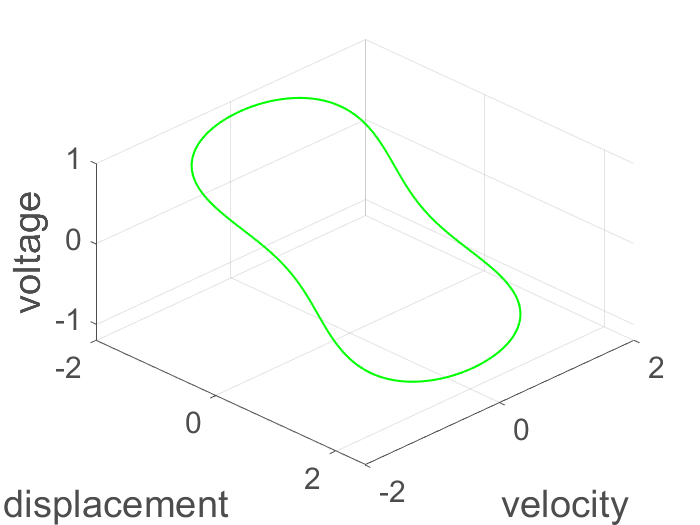}
         \caption{$f = 0.211$}
     \end{subfigure}
     \begin{subfigure}[b]{0.31\textwidth}
         \centering
         \includegraphics[width=\textwidth]{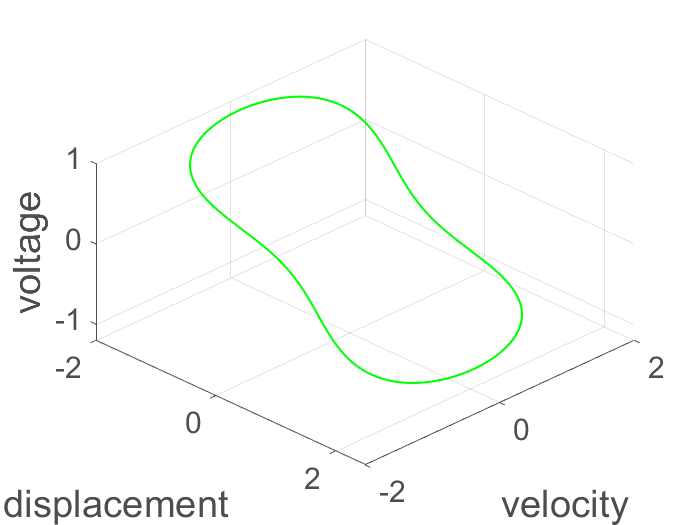}
         \caption{$f = 0.243$}
     \end{subfigure}
     \begin{subfigure}[b]{0.31\textwidth}
         \centering
         \includegraphics[width=\textwidth]{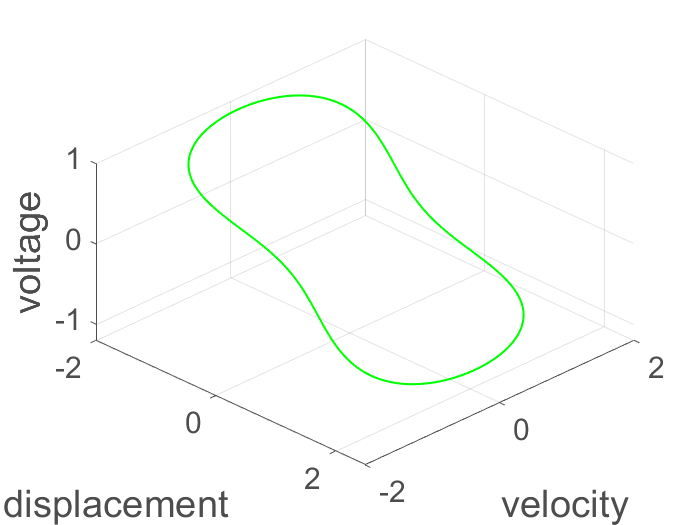}
         \caption{$f = 0.275$}
     \end{subfigure}

	\caption{Attractors of the asymmetric model for $\Omega=0.8$ and $\phi=-4.95^\circ$ under different forcing amplitudes. Attractors correspond to basins of attraction presented in Fig.\ref{fig:basins_phi-4}.}
	\label{fig:atractor_phi-4}
\end{figure}

The basins of attraction and attractors are calculated for the optimal asymmetric bistable model corresponding to quadratic term $\delta$ equal to 0.15. At this condition, the sloping plane is $\phi_{opt} = -4.95^\circ$. Fig.~\ref{fig:basins_phi-4} and \ref{fig:atractor_phi-4} depict the basins of attraction and their corresponding attractors, respectively. The excitation frequency was fixed at $\Omega = 0.8$ while the amplitude was varied from $0.019 \leq f \leq 0.275$ with intervals of $0.032$. The initial conditions were defined as $(x_0,\dot{x}_0) \in [-3 \leq x_0 \leq 3] \times [-3 \leq \dot{x}_0 \leq 3]$ on the restriction plane of $v_0 = 0$.

At an amplitude of $f=0.019$, four solutions are represented by blue, cyan, red, and magenta basins with low energy responses. The green basin emerged at $f=0.051$, resulting in five solutions. At $f=0.083$, besides the red and blue basins, there is an appearance of a chaotic solution. At this condition, the green attractor disappeared. At $f=0.115$, two solutions displayed inter-well motions, one chaotic and one regular, although the chaotic solution was barely noticeable, and the green basins filled the entire region. With increasing amplitude, the chaotic solution vanished, leaving only the green attractor as a solution.

It is important to note that the attractors were symmetrical for this sloping angle value, indicating that the depth of the potential energy wells is equal. The asymmetric effect is compensated despite the asymmetric terms added in the model, resulting in a similar behavior to the symmetric system. The negative effects of the asymmetric model previously seen were canceled, and the green basins dominated, particularly at high values of the excitation amplitude. Comparing the basins of attraction and attractors of the asymmetric system to the symmetric (Figs.~\ref{fig:basins_O08} and \ref{fig:atractors_O08}), they are found to be similar, with the solutions exhibiting similar behaviors.


\begin{figure}
    \centering
    \begin{subfigure}[b]{0.31\textwidth}
         \centering
         \includegraphics[width=\textwidth]{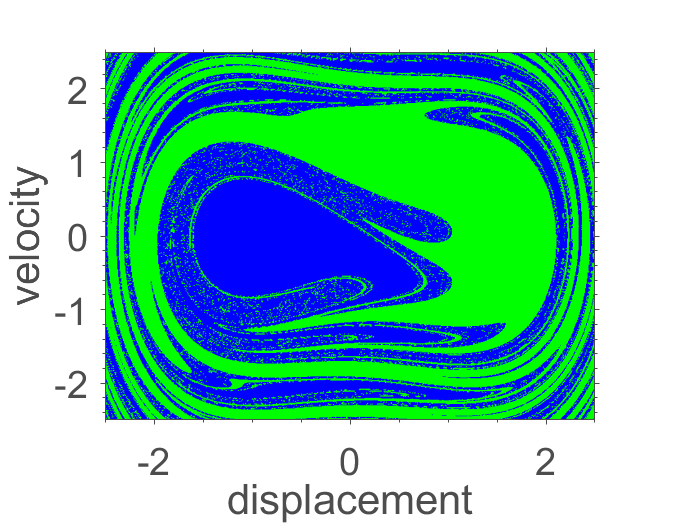}
         \caption{$\phi = -35^\circ$}
     \end{subfigure}
     \begin{subfigure}[b]{0.31\textwidth}
         \centering
         \includegraphics[width=\textwidth]{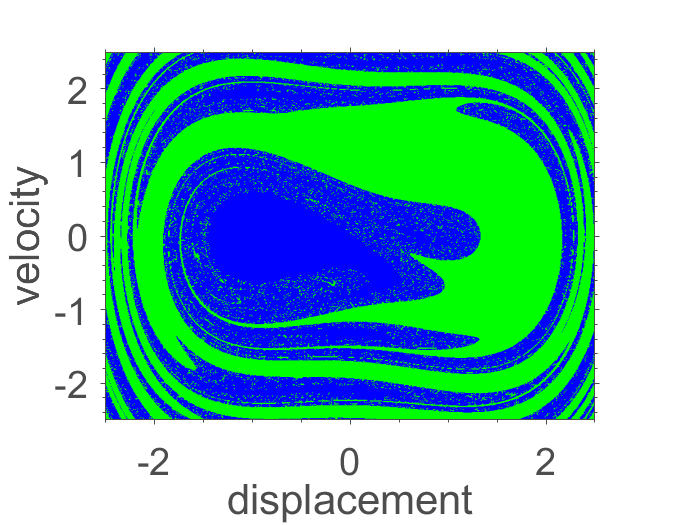}
         \caption{$\phi = -25^\circ$}
     \end{subfigure}
     \begin{subfigure}[b]{0.31\textwidth}
         \centering
         \includegraphics[width=\textwidth]{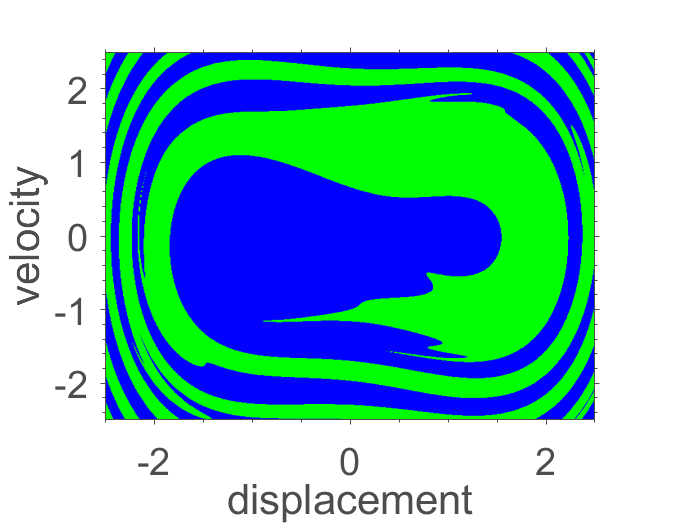}
         \caption{$\phi = -15^\circ$}
     \end{subfigure}
     \begin{subfigure}[b]{0.31\textwidth}
         \centering
         \includegraphics[width=\textwidth]{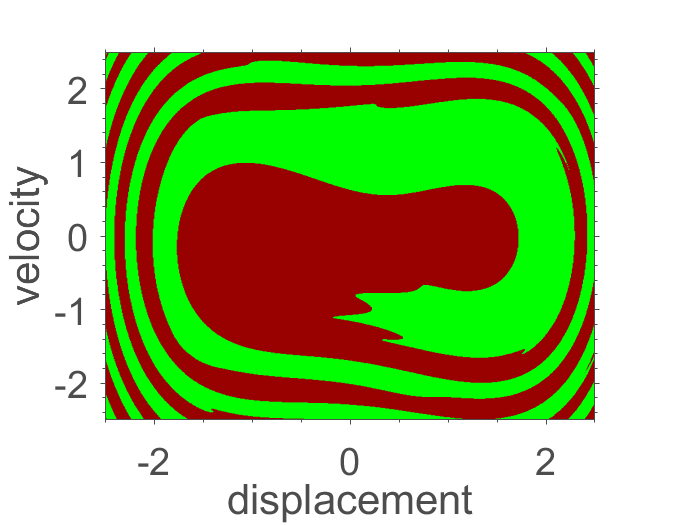}
         \caption{$\phi = -5^\circ$}
     \end{subfigure}
     \begin{subfigure}[b]{0.31\textwidth}
         \centering
         \includegraphics[width=\textwidth]{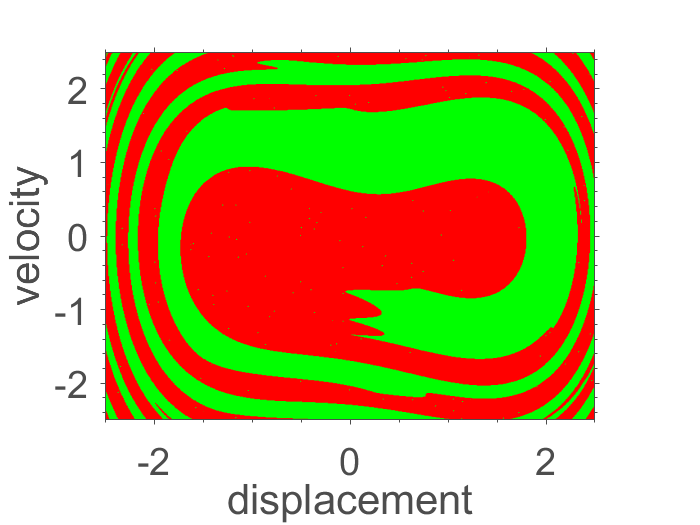}
         \caption{$\phi = 0^\circ$}
     \end{subfigure}
     \begin{subfigure}[b]{0.31\textwidth}
         \centering
         \includegraphics[width=\textwidth]{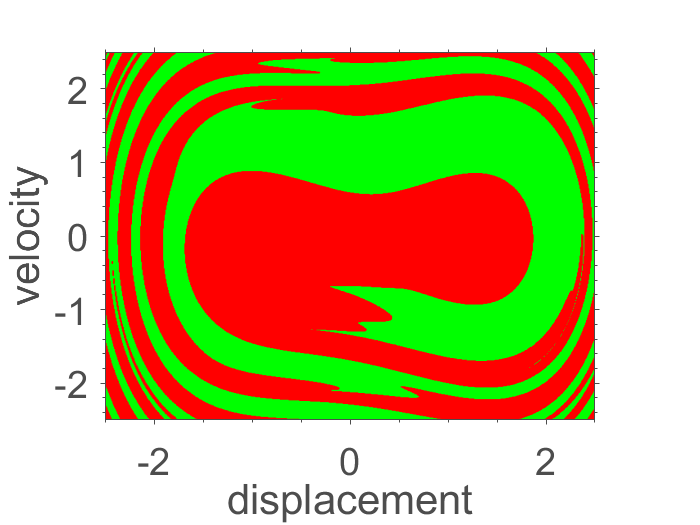}
         \caption{$\phi = 5^\circ$}
     \end{subfigure}
     \begin{subfigure}[b]{0.31\textwidth}
         \centering
         \includegraphics[width=\textwidth]{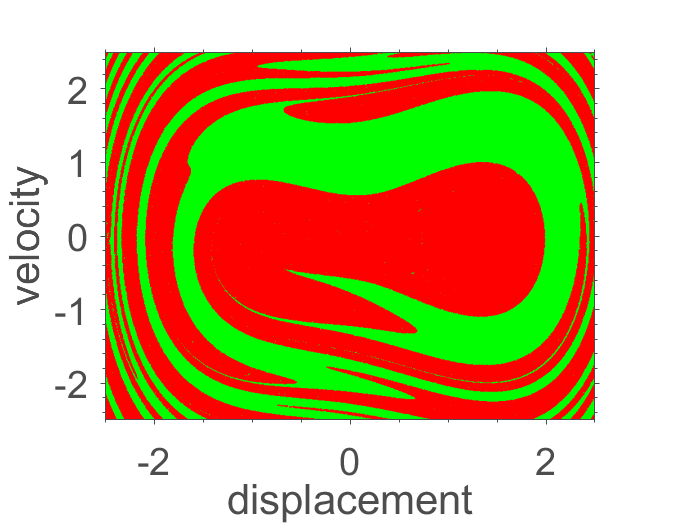}
         \caption{$\phi = 15^\circ$}
     \end{subfigure}
     \begin{subfigure}[b]{0.31\textwidth}
         \centering
         \includegraphics[width=\textwidth]{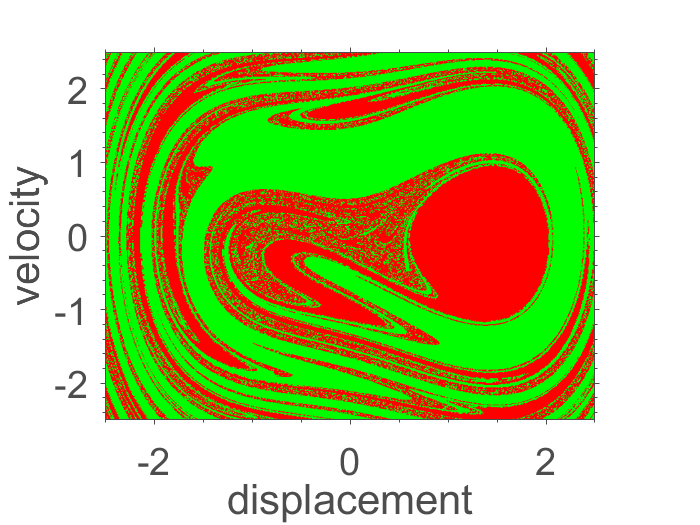}
         \caption{$\phi = 25^\circ$}
     \end{subfigure}
     \begin{subfigure}[b]{0.31\textwidth}
         \centering
         \includegraphics[width=\textwidth]{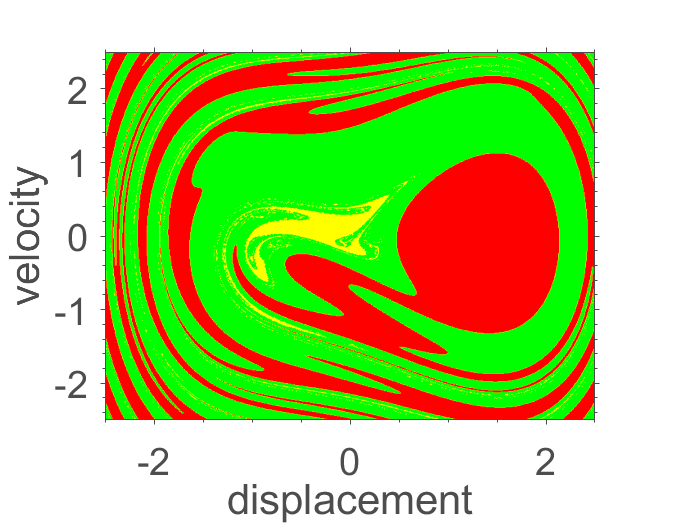}
         \caption{$\phi = 35^\circ$}
     \end{subfigure}

    \caption{Basin of attraction intersection with the plane $v=0$ for the asymmetric model with a forcing frequency of $\Omega=0.8$ and forcing amplitude $\mathnormal{f}=0.091$, within the bounds of $-2.5 \leq x_0 \leq 2.5$ and $-2.5 \leq \dot{x}_0 \leq 2.5$. Different sloping angles are investigated. The corresponding attractor colors are shown in Fig.\ref{fig:atractor_f091}.}
    \label{fig:basins_f091}
\end{figure}

\begin{figure}
    \centering
            \begin{subfigure}[b]{0.31\textwidth}
         \centering
         \includegraphics[width=\textwidth]{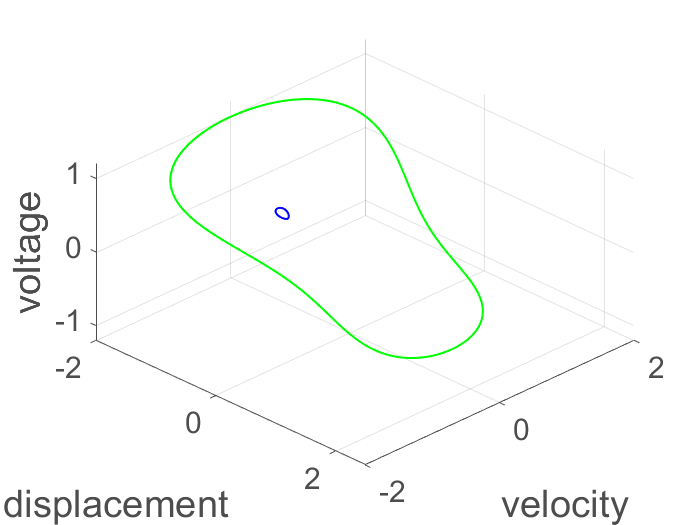}
         \caption{$\phi = -35^\circ$}
     \end{subfigure}
     \begin{subfigure}[b]{0.31\textwidth}
         \centering
         \includegraphics[width=\textwidth]{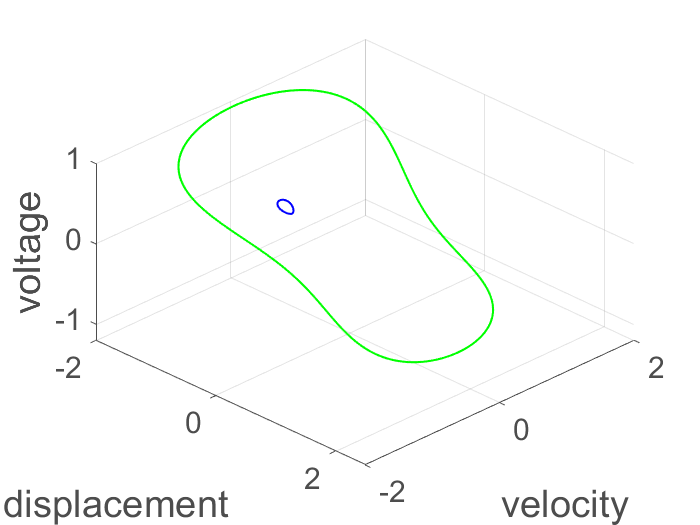}
         \caption{$\phi = -25^\circ$}
     \end{subfigure}
     \begin{subfigure}[b]{0.31\textwidth}
         \centering
         \includegraphics[width=\textwidth]{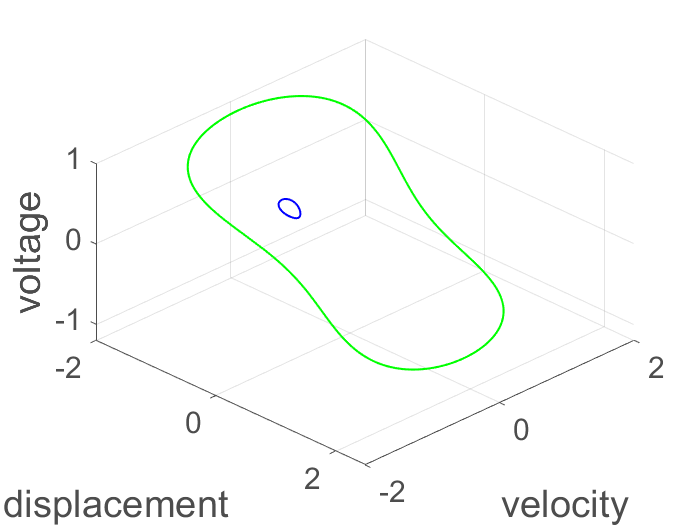}
         \caption{$\phi = -15^\circ$}
     \end{subfigure}
     \begin{subfigure}[b]{0.31\textwidth}
         \centering
         \includegraphics[width=\textwidth]{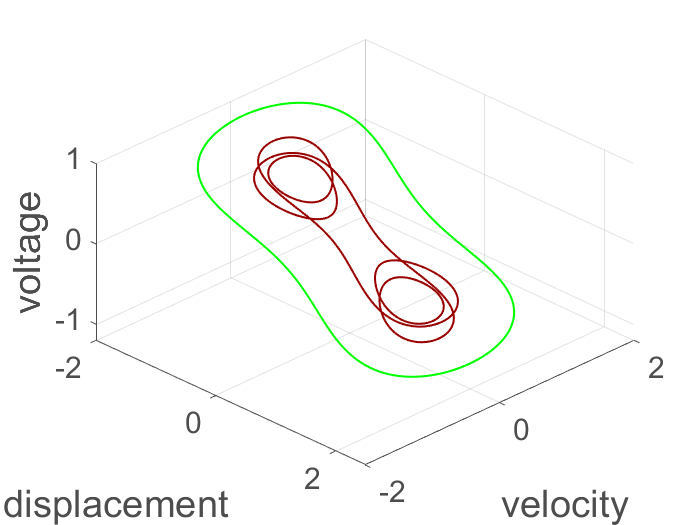}
         \caption{$\phi = -5^\circ$}
     \end{subfigure}
     \begin{subfigure}[b]{0.31\textwidth}
         \centering
         \includegraphics[width=\textwidth]{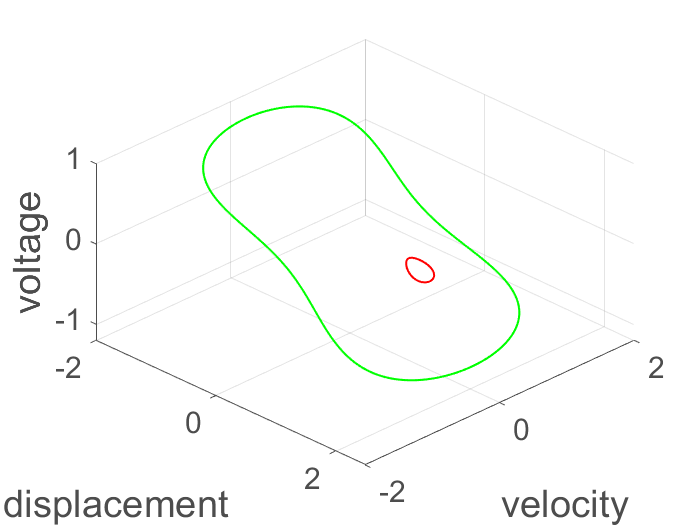}
         \caption{$\phi = 0^\circ$}
     \end{subfigure}
     \begin{subfigure}[b]{0.31\textwidth}
         \centering
         \includegraphics[width=\textwidth]{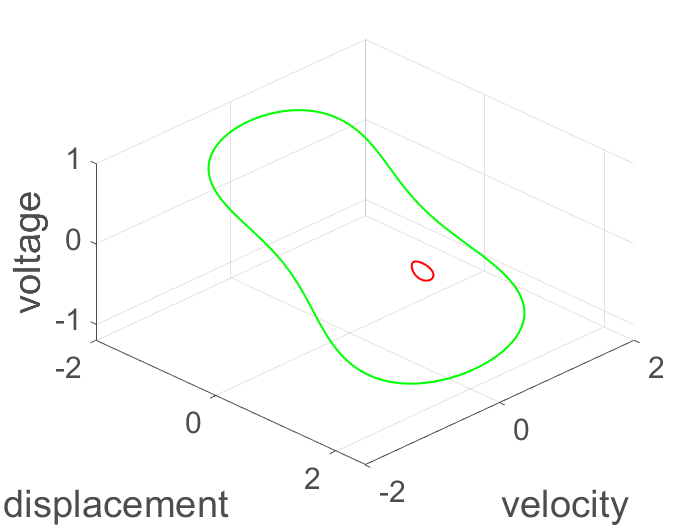}
         \caption{$\phi = 5^\circ$}
     \end{subfigure}
     \begin{subfigure}[b]{0.31\textwidth}
         \centering
         \includegraphics[width=\textwidth]{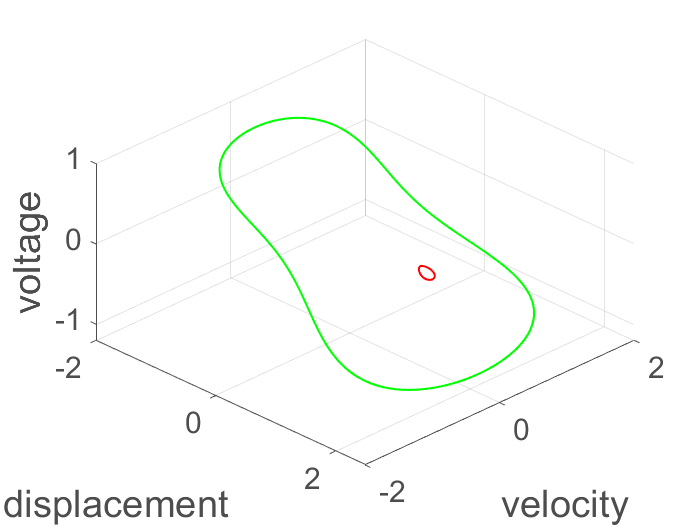}
         \caption{$\phi = 15^\circ$}
     \end{subfigure}
     \begin{subfigure}[b]{0.31\textwidth}
         \centering
         \includegraphics[width=\textwidth]{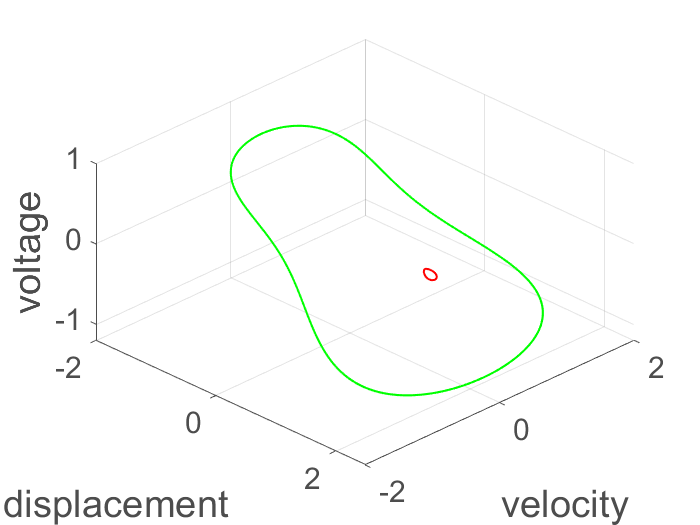}
         \caption{$\phi = 25^\circ$}
     \end{subfigure}
     \begin{subfigure}[b]{0.31\textwidth}
         \centering
         \includegraphics[width=\textwidth]{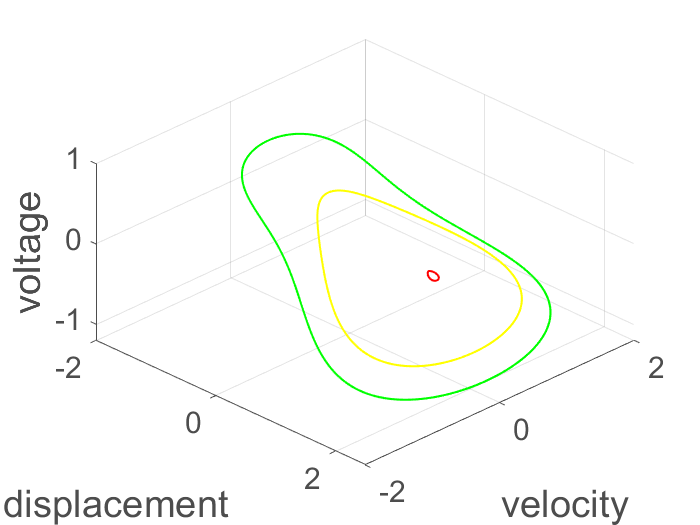}
         \caption{$\phi = 35^\circ$}
     \end{subfigure}
	\caption{Attractors of the asymmetric model for $\Omega=0.8$ and $\mathnormal{f}=0.091$ under different sloping angle values. Attractors correspond to basins of attraction presented in Fig.~\ref{fig:basins_f091}.}
	\label{fig:atractor_f091}
\end{figure}

 The Fig.~\ref{fig:basins_f091} and \ref{fig:atractor_f091} present the basins of attraction and their respective attractors for $\mathnormal{f} = 0.091$, $\Omega = 0.8$ and $\delta=0.15$ with the sloping angle $\phi$ varying between ${-35^\circ, -25^\circ, -15^\circ, -5^\circ, 0^\circ, 5^\circ, 25^\circ, 35^\circ}$. The initial conditions $(x_0,\dot{x}_0)$ are defined in the range of $[-2.5 \leq x_0 \leq 2.5] \times [-2.5 \leq \dot{x}_0 \leq 2.5]$ on the restriction plane of $v_0 = 0$.

For $\phi \in {-35^\circ, -25^\circ, -15^\circ}$, there are two attractors, blue and green, with the blue being a low-energy response around the negative equilibrium point. At $\phi = -5^\circ$, which is close to the optimal angle, the system has two basins, green and purple. The purple basin corresponds to multiple periods of inter-well motion. At $\phi = 0^\circ$, where the asymmetry is driven solely by the $\delta$ coefficient, there are two solutions, the green and red attractors. The red attractor corresponds to a low-energy orbit around the positive equilibrium point. As $\phi$ increases, so does the asymmetry, resulting in the presence of both the green and red basins. For $\phi = 35^\circ$, there is also a tiny yellow basin with a homoclinic orbit.

The asymmetry in the system is inverted after $\phi = -5^\circ$. For negative values of $\phi$, the blue basin is concentrated at the negative displacement values. As $\phi$ approaches the optimal value, a portion of the blue region moves towards positive displacement values. After the optimal value, the blue basin is replaced by the red basin and no longer attracts any responses. As $\phi$ increases, the red solution concentrates around the positive equilibrium point.

These results also show that the basins of attraction around the optimal angle $\phi_{opt}=-4.95^\circ$ (Fig. \ref{fig:basins_f091}c, d, and e) exhibit a similar green area. The only difference lies in the basin of attraction close to the equilibrium points, where for $\phi=0^\circ$ and $\phi=5^\circ$, the red basin dominates the area near equilibrium points. In contrast, for $\phi=-5^\circ$, the purple basin dominates this region. Despite the condition $\phi=-5^\circ$ being almost exactly at the optimal angle value and presenting better conditions compared to the other angle values in Fig. \ref{fig:basins_f091} as it only presents high energy orbits, a small difference in the angle did not damage the green region. Therefore, the sloping angle can mitigate the negative impact of asymmetry, and at its optimal value, it can completely cancel the asymmetry's impact.

\begin{figure}
    \centering
    \includegraphics[width=\textwidth]{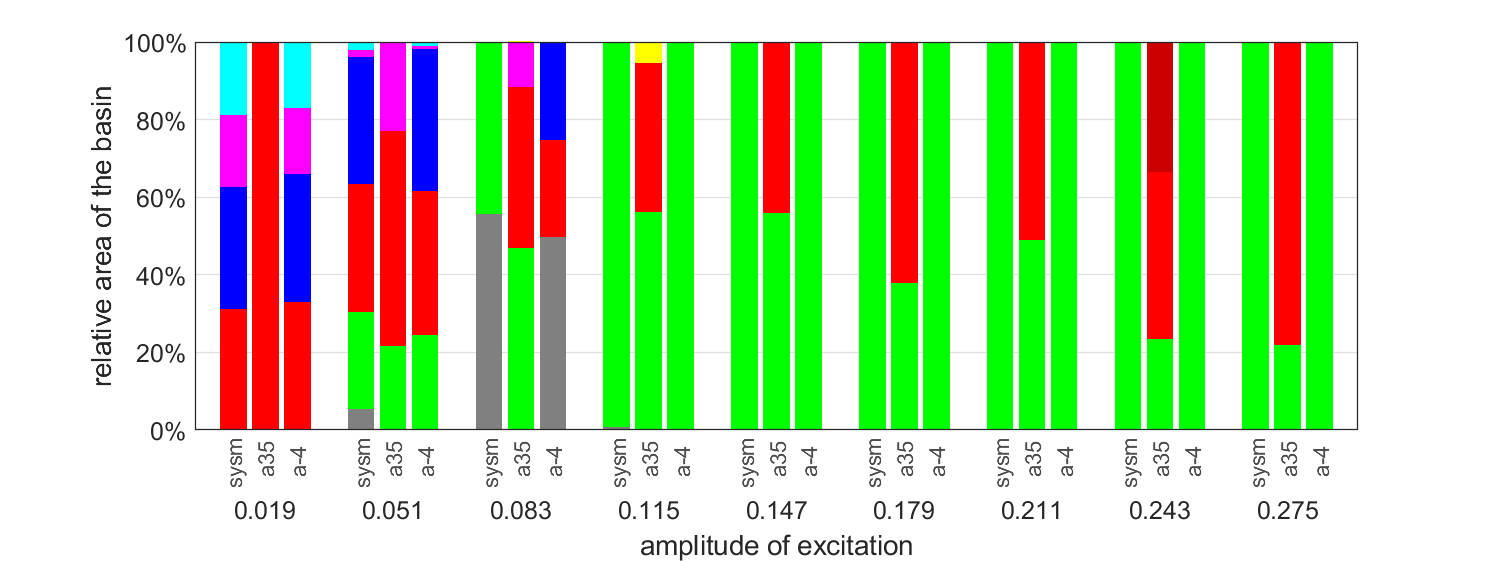}
	\caption{Comparison of the relative area of the basins of attraction for the symmetric model (sysm), asymmetric model with $\phi=35^\circ$ (a35), and optimal asymmetric model with $\phi=-4.95^\circ$ (a-4) for a fixed frequency of excitation $\Omega=0.8$ and varying amplitude of excitation. The area of the basins is measured relative to the full range of initial conditions. The red and blue bars indicate the red and blue basins, respectively, where the attractors have low-energy orbits. The green bar represents the basin where the attractor has a high-energy orbit, while the gray bar represents the basin where the attractor exhibits chaotic orbits. The pink, cyan, yellow, and purple bars correspond to basins where the attractor presents multi-periodic orbits.}
	\label{fig:basins_area}
\end{figure}

Finally, Fig.~\ref{fig:basins_area} compares the occupied area of each basin of attraction in a defined domain of initial conditions ($\left\lbrace (x_0,\dot{x}_0) ~ | ~ -3 \leq x_0 \leq 3 ~ \mbox{and} ~ -3 \leq \dot{x}_0 \leq 3 \right\rbrace$) for three different models: symmetric (sysm), asymmetric with $\phi=35^\circ$ (a35), and optimal asymmetric (a-4), under different excitation amplitudes with $\Omega = 0.8$. The occupied area ratio, the area of each colored basin divided by the area of the total domain, is presented using a stacked bar graph. Each color represents a specific basin as previously defined. The symmetric and optimal asymmetric models exhibit similar occupied area ratios for the green, red, blue, pink, cyan, and gray basins. Except for $f=0.083$, the symmetric system displays a green bar (indicating a high vibration amplitude). In contrast, the optimal asymmetric system shows red and blue bars (indicating a low vibration amplitude). At higher excitation amplitude, the green basin becomes dominant for both systems, which is the ideal scenario for energy harvesting since the green attractor has a high vibration amplitude. In contrast, in the asymmetric model with $\phi = 35^\circ$, the red basin dominates the area of the domain. Even as the excitation amplitude increases, the green basin only reaches 55\%, when $f=0.115$, of the occupied area before significantly dropping. For all excitation amplitudes, the asymmetric model with $\phi = 35^\circ$ has a smaller green basin area than the other models. Thus, this figure demonstrates that asymmetries harm energy harvesting. However, this effect can be mitigated when the asymmetry angle excitation is optimal.

\section{Conclusions}
\label{concl_sect}

This paper presents a numerical study of the dynamic behavior of bistable energy harvesters. Both symmetric and asymmetric models were analyzed, with the latter incorporating two sources of asymmetric, a quadratic term for the nonlinear restoring force and a sloping angle in the plane where the system is attached. Three cases of asymmetry were studied. One system with asymmetry only from the quadratic term, one with a strong influence from the quadratic term and a high value of sloping angle, and another that the sloping angle value canceled the quadratic term, mitigating the asymmetry effect. Dynamic responses for different excitation and initial conditions using bifurcation diagrams and basins of attractions were investigated.

A comprehensive analysis was conducted by varying the amplitude and frequency of excitation over a wide range, spanning from low to high energy. The findings revealed the existence of multiple solutions at forward and backward bifurcations in the amplitude and frequency diagrams. Chaotic response bands appeared for both the symmetric and optimal asymmetric models. In contrast, the asymmetric system did not exhibit such behavior. The strong asymmetric system resulted in low-amplitude motions. The symmetric and optimal asymmetric models displayed similar behaviors across all the studied diagrams. Moreover, the sloping angle was swept from $-35^\circ$ to $35^\circ$ for different forcing frequencies. Chaotic bands could emerge for angles near $0^\circ$ when the forcing amplitude was high. For all cases, multiple-period responses rarely appeared.

Basins of attraction based on the 0-1 test for chaos and their corresponding attractor were used to analyze the multiple solutions for initial conditions for symmetric and asymmetric models under different excitation conditions. The 0-1 test for chaos was used as cost-effective in classifying dynamic behavior. Our results showed that the basins of attraction for all cases exhibited diverse dynamics, including chaotic, one-periodic, and multi-periodic responses, with high and low-energy orbits. The symmetric model displayed periodic oscillation with high amplitude, which is the most favorable scenario for energy harvesting. However, an intensified density region for the low-energy orbit was observed when we introduced asymmetry through the quadratic term. For a strong asymmetric model with a high sloping angle value, the low-energy orbit dominated the basins of attraction, even with high excitation amplitude. In the optimal value for the sloping angle, we found that it can cancel the negative effect of asymmetry, leading to similar basins of attraction as in the symmetric model. Finally, the basins of attraction for different sloping angle values demonstrated that angles close to the optimal angle presented the same area of high orbits. It indicated that the sloping angle could mitigate the drastically harmful impact of the asymmetry.


For future works, it would be valuable to construct an experimental apparatus to investigate the dynamic behavior of the system under various excitation conditions and asymmetries. These experiments would serve to validate the findings of this study and provide empirical confirmation. Additionally, the experimental investigations could enhance our understanding of the system's dynamics, verify the robustness of the optimal angle, explore the sensitivity of physical parameters on the system's behavior, and contribute to the advancement of knowledge in this field.

\section*{Acknowledgements}

The authors express their gratitude to Prof. Marcelo Savi (UFRJ), Prof. Paulo Paupitz Gonçalves (UNESP), and Prof. Marcus Varanis (UFMS) for their valuable contributions to the discussions regarding the results presented in this work. They also extend their appreciation to Dr. Helder Casa Grande (UERJ) for his thorough review of the manuscript and to Mr. Michel Tosin (UERJ) for creating Figure~\ref{fig:test01}.

\section*{Funding}

This research received financial support from the Brazilian agencies Coordena\c{c}\~{a}o de Aperfei\c{c}oamento de Pessoal de N\'{\i}vel Superior - Brazil (CAPES) with the Finance Code 001, Conselho Nacional de Desenvolvimento Cient\'{\i}fico e Tecnol\'{o}gico (CNPq) grant 305476/2022-0, and Funda\c{c}\~{a}ão Carlos Chagas Filho de Amparo \`{a} Pesquisa do Estado do Rio de Janeiro (FAPERJ) under the grants: 211.304/2015, 210.021/2018, 210.167/2019, 211.037/2019 and 201.294/2021.

\section*{Code availability}

The simulations presented in this paper were carried out using the {\bf STONEHENGE -- Suite for Nonlinear Analysis of Energy Harvesting Systems} software \cite{STONEHENGE_paper}. The authors make the code accessible to the scientific community via a free, open-source repository on GitHub \cite{STONEHENGE}. This allows for the reproducibility of the results.


\section*{Conflict of Interest}
The authors declare they have no conflict of interest.

\section*{Disclaimer}

This manuscript underwent a thorough grammatical revision and improvement with the assistance of artificial intelligence-powered tools such as Grammarly and ChatGPT. Nevertheless, the authors retain full responsibility for the original language and phrasing.

\bibliographystyle{elsarticle-num}


\end{document}